# Biprops

Volodymyr Lyubashenko*

April 16, 2026

**Abstract**

We define biprops as a generalization of coloured props and of symmetric weak multicategories. These are bicategories whose objects form a free monoid. They are equipped with some structure resembling a symmetric strict tensor product. We prove that a symmetric weak multicategory gives rise to a biprop and a symmetric weak multifunctor gives rise to a morphism of biprops. This is a functor from the category of symmetric weak multicategories to the category of biprops. [1]

### 0.1 Conventions

We use the symmetric strict monoidal category $\mathcal{S}_{\mathsf{sk}}$ skeletal for the category of finite totally ordered sets and arbitrary maps. Thus, $\operatorname{Ob}\mathcal{S}_{\mathsf{sk}} \cong \mathbb{N} = \{0, 1, 2, \dots\}$. The monoidal product of $I, J \in \operatorname{Ob}\mathcal{S}_{\mathsf{sk}}$ is their lexicographic product, which is the cartesian product of sets $I \times J$, equipped with the ordering $(i_1, j_1) <'' (i_2, j_2)$ iff $(j_1 < j_2$ or $(j_1 = j_2$ and $i_1 < i_2))$. We denote this monoidal product by $I \times J$ although it is not a cartesian one. The symmetry is that of the category of sets. The strict monoidal subcategory $\mathcal{O}_{\mathsf{sk}} \subset \mathcal{S}_{\mathsf{sk}}$ has the same objects, but only order-preserving maps.

### 0.2 Introduction

The notion of a biprop generalizes that of coloured prop (see e.g. [HR15, Section 1.1]). It arises also as a generalization of symmetric weak multicategories [Lyu25]. It is a kind of symmetric monoidal bicategory $\mathcal{P}$. Objects of this bicategory form a free monoid on the set of colours. So an object is a tuple of colours $(X_i)_{i\in I}$, $I \in \mathcal{S}_{\mathsf{sk}}$, $X_i \in \operatorname{Col}\mathcal{P}$. $\operatorname{Col}\mathcal{P}$, the set of colours of $\mathcal{P}$, is an analogue of the set of objects of a symmetric weak multicategory.

Monoidal structure of bicategory $\mathcal{P}$ is not a conventional one, for tensor product functors are indexed by cospans $I \xrightarrow{f} L \xleftarrow{g} J$ in $\mathcal{S}_{\mathsf{sk}}$. Heuristically, such a tensor product combines a conventional tensor product with a symmetry. Besides, the tensor product is strictly associative in appropriate sense (1.1.1) and (1.1.6). On the other hand, the composition in $\mathcal{P}$ is not strictly associative, so $\mathcal{P}$ is not a 2-category. This notion or its version with conventional tensor products may indicate an alternative for the usual biased definition of a symmetric monoidal bicategory (cf. [Sta16, Definition 4.8], [SP11, Definition 2.3] and references therein).

The usual definition (e.g. [HR15, Section 1.1]) of props uses explicit action of symmetric groups permuting the inputs and outputs. In Leinster's definition of symmetric multicategory explicit action of symmetric groups is not used. Instead he constructs this action in [Lei03, Lemma A.2.2] departing from the structure data of a symmetric multicategory. We generalize this action to biprops so that a permutation of inputs (or outputs) acts on the left (or on the right) by an equivalence of categories (Corollary 2.1.3), the identity permutation acts by an equivalence isomorphic to the identity functor (Remark 2.1.1), the composition of two permutations acts by an equivalence isomorphic to the composition of the actions of two permutations-factors (Proposition 2.1.2). Moreover, these isomorphisms satisfy the non-abelian cocycle identities (Propositions 2.1.5 and 2.1.6). In this sense the action of a symmetric group on hom-categories of a biprop is weak. Left and right actions commute up to an isomorphism (Proposition 2.1.4). They agree with the composition in an expected way (Remarks 2.1.7–2.1.9).

* *Email address:* lub@imath.kiev.ua

The main reason for the existence of biprops are symmetric weak multicategories. We construct a functor which sends a symmetric weak multicategory to a biprop, a symmetric weak multifunctor to a morphism of biprops.

The guiding principle in writing definitions of biprops, their morphisms, their natural transformations and their modifications is contractibility. In our case this means that all natural transformations constructed from the structure data, must be equal.

### 0.3 Notations

Instead of pasting diagrams in bicategories (=extendable pasting schemes) of [JY21, §3] we use Bartlett's notation [Bar14] for computations in the symmetric cartesian 2-category $\mathcal{C}at$. Its objects are locally small categories.

**Acknowledgement.** I thank Paul Taylor for writing a package 'diagrams.sty', which I use a lot. I am grateful for clarifying discussions to Prof. Dr. Anna Beliakova and Prof. Dr. Giovanni Felder. Excellent environment for work in the University of Zurich and of the Swiss Federal Institute of Technology in Zurich were due to the generosity of the National Centre of Competence in Research SwissMAP of the Swiss National Science Foundation (grant number 205607), to whom I express my deep gratitude. The author thanks Prof. Dr. Christoph Schweigert and Faculty of Mathematics, Computer Science and Natural Sciences of the University of Hamburg for interest in the subject and excellent conditions of work. I thank Prof. Tom Leinster at the University of Edinburgh for fruitful discussions. The author wishes to thank the Isaac Newton Institute and London Mathematical Society for the financial support and the School of Mathematics at the University of Edinburgh for their hospitality. I thank Prof. Paolo Papi from Sapienza Università di Roma and Prof. Paolo Stellari from Università degli studi di Milano for creating nice conditions of work. I am really grateful to the Armed Forces of Ukraine who gave me the possibility to work quietly on the subject.

# 1 Biprops

## 1.1 Main properties of biprops

**1.1.1 Definition.** A *biprop* is a bicategory $(\mathcal{P}, 1_X, \mathsf{a}, \mathsf{r}, \mathsf{l})$ together with

— The set of objects $\operatorname{Ob}\mathcal{P} = S^*$, the free monoid generated by a set $S = \operatorname{Col}\mathcal{P}$. The concatenation operation in $S^*$ is denoted $\bigsqcup_{l\in L}: (S^*)^L \to S^*$, $(A_l)_{l\in L} \mapsto \bigsqcup_{l\in L} A_l$.

— Monoidal products – the system of maps

$$\otimes^{I\xrightarrow{f}L}: (\operatorname{Ob}\mathcal{P})^L \to \operatorname{Ob}\mathcal{P}: \big(\bigsqcup_{i\in f^{-1}l} A_i\big)_{l\in L} \mapsto \bigsqcup_{i\in I} A_i, \quad A_i \in \operatorname{Ob}\mathcal{P},$$

indexed by a map $f: I \to L \in \mathcal{S}_{\mathsf{sk}}$, and of functors

$$\otimes^{I\xrightarrow{f}L\xleftarrow{g}J}: \prod_{l\in L}\mathcal{P}\big(\bigsqcup_{i\in f^{-1}l} A_i; \bigsqcup_{j\in g^{-1}l} B_j\big) \to \mathcal{P}\big(\bigsqcup_{i\in I} A_i; \bigsqcup_{j\in J} B_j\big),$$

indexed by a pair of maps $I \xrightarrow{f} L \xleftarrow{g} J \in \mathcal{S}_{\mathsf{sk}}$. This system is supposed to be strict, that is (compare with strict homomorphisms of Bénabou [Bén67, Remark 4.2], see also e.g. [Lei98, § 1.1]), there is a commutative diagram of functors

$$\begin{array}{ccc}
\prod_{l\in L}\mathcal{P}\big(\bigsqcup_{i\in f^{-1}l} A_i; \bigsqcup_{j\in g^{-1}l} B_j\big) \times \prod_{l\in L}\mathcal{P}\big(\bigsqcup_{j\in g^{-1}l} B_j; \bigsqcup_{k\in h^{-1}l} C_k\big) & \xrightarrow{\prod_L m} & \prod_{l\in L}\mathcal{P}\big(\bigsqcup_{i\in f^{-1}l} A_i; \bigsqcup_{k\in h^{-1}l} C_k\big) \\
\Big\downarrow {\scriptstyle \otimes^{I\xrightarrow{f}L\xleftarrow{g}J} \times \otimes^{J\xrightarrow{g}L\xleftarrow{h}K}} & = & \Big\downarrow {\scriptstyle \otimes^{I\xrightarrow{f}L\xleftarrow{h}K}} \\
\mathcal{P}\big(\bigsqcup_{i\in I} A_i; \bigsqcup_{j\in J} B_j\big) \times \mathcal{P}\big(\bigsqcup_{j\in J} B_j; \bigsqcup_{k\in K} C_k\big) & \xrightarrow{m} & \mathcal{P}\big(\bigsqcup_{i\in I} A_i; \bigsqcup_{k\in K} C_k\big)
\end{array} \tag{1.1.1}$$

As follows from axioms (1.1.2), (1.1.4), (1.1.5) (below) the collection $\otimes^{I\xrightarrow{f}L\xleftarrow{f}I}$ is a strict homomorphism of bicategories in the sense of Bénabou [Bén67, Remark 4.2] (in unbiased form).

Assume given maps in $\mathcal{S}_{\mathsf{sk}}$

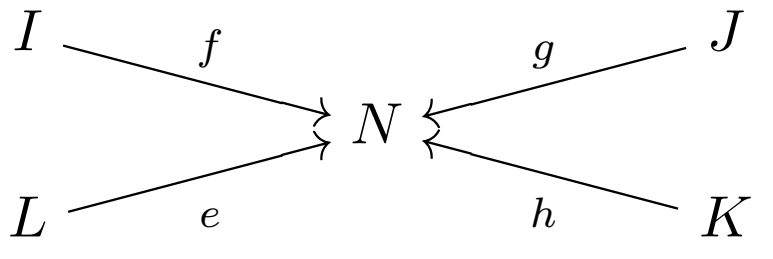

Let $(A_i)_{i\in I}$, $(B_j)_{j\in J}$, $(C_k)_{k\in K}$, $(D_l)_{l\in L}$ be families of objects of $\mathcal{P}$. It is required that there is a commutative cube (or its remnants)

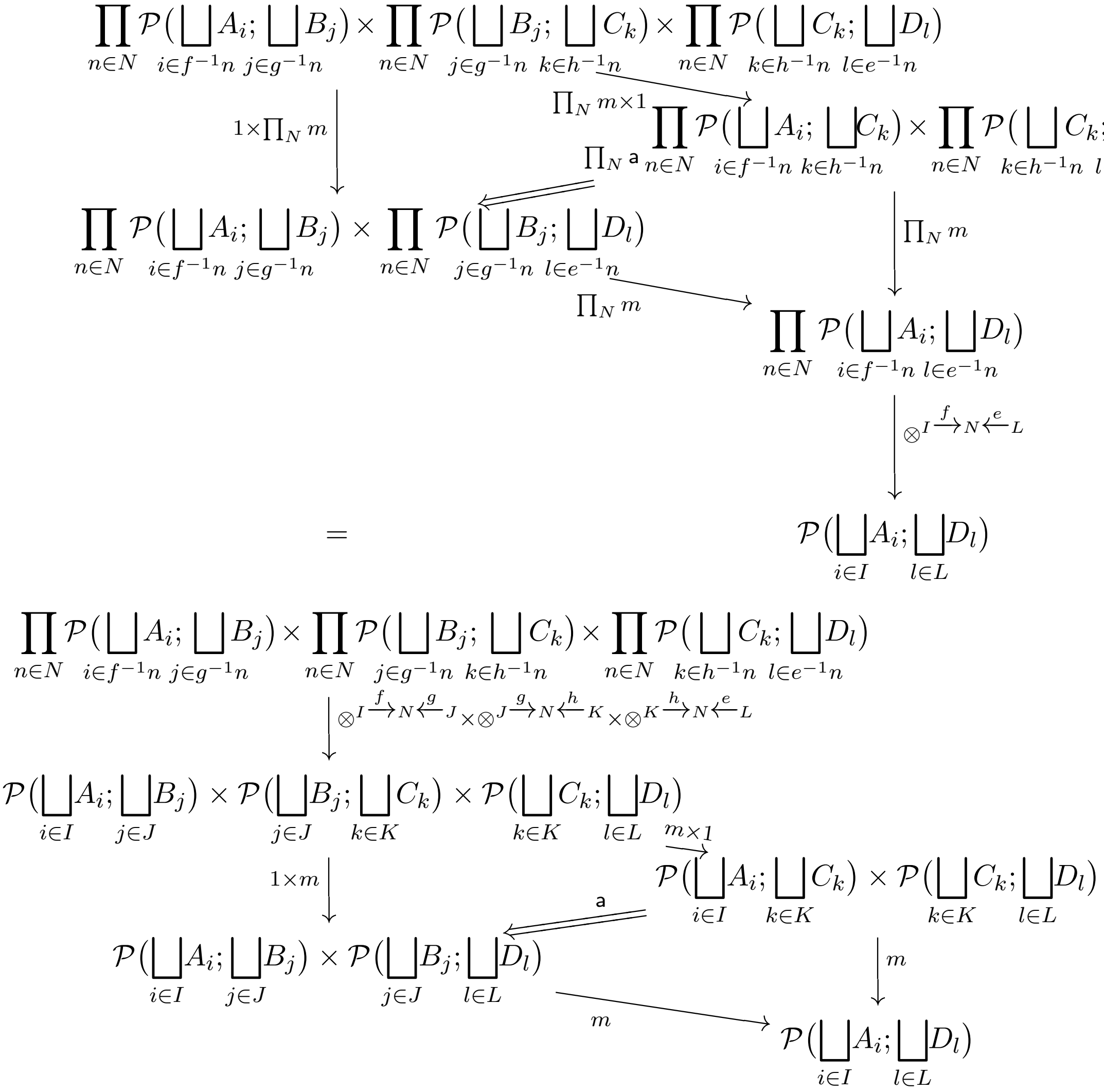

(1.1.2)

Equality of source and target functors follows from diagram (1.1.1).

It is required that $\otimes^{I\xrightarrow{\triangledown}\mathbf{1}\xleftarrow{\triangledown}J}: \mathcal{P}\big(\bigsqcup_{i\in I} A_i; \bigsqcup_{j\in J} B_j\big) \to \mathcal{P}\big(\bigsqcup_{i\in I} A_i; \bigsqcup_{j\in J} B_j\big)$ is the identity functor. Furthermore,

$$\otimes^{I\xrightarrow{f}L\xleftarrow{f}I}: \prod_{l\in L} \mathcal{P}\Big(\bigsqcup_{i\in f^{-1}l} A_i; \bigsqcup_{i\in f^{-1}l} A_i\Big) \to \mathcal{P}\Big(\bigsqcup_{i\in I} A_i; \bigsqcup_{i\in I} A_i\Big), \quad (1)_{l\in L} \mapsto 1, \tag{1.1.3}$$

and

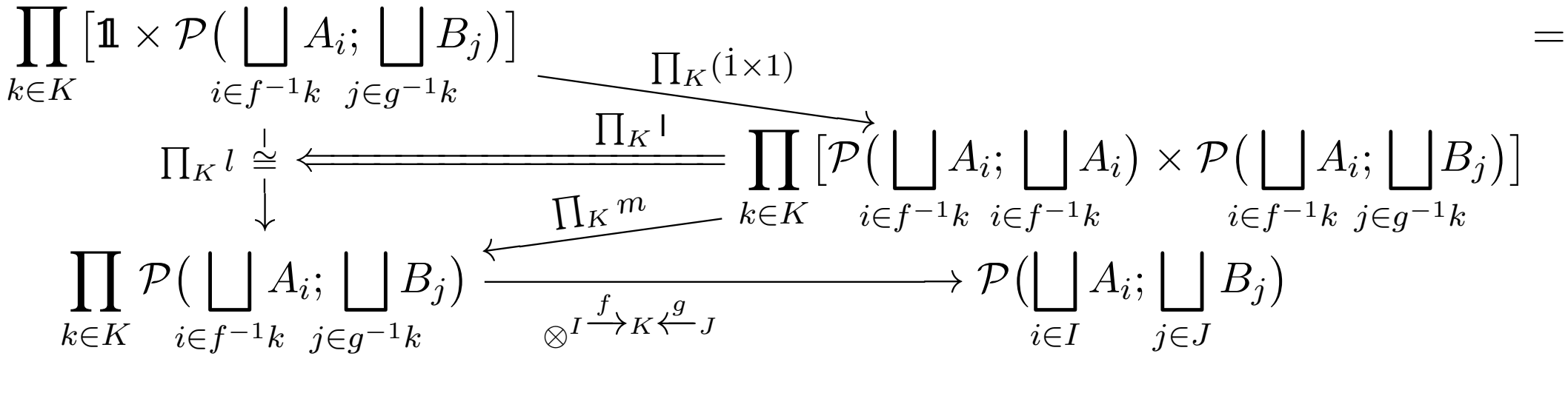

$$
\begin{array}{ccc}
\mathbb{1}\times\mathcal{P}(\bigsqcup_{i\in I}A_i;\bigsqcup_{j\in J}B_j) & \xleftarrow{\;l^{-1}\circ\otimes^{I\xrightarrow{f}K\xleftarrow{g}J}\circ\prod_K l\;} & \prod_{k\in K}[\mathbb{1}\times\mathcal{P}(\bigsqcup_{i\in f^{-1}k}A_i;\bigsqcup_{j\in g^{-1}k}B_j)] \\
 & \searrow^{\mathrm{i}\times 1} & \downarrow \\
l \downarrow \cong & \xLeftarrow{\quad\mathsf{l}\quad} & \mathcal{P}(\bigsqcup_{i\in I}A_i;\bigsqcup_{i\in I}A_i)\times\mathcal{P}(\bigsqcup_{i\in I}A_i;\bigsqcup_{j\in J}B_j) \\
 & \swarrow_{m} & \\
\mathcal{P}(\bigsqcup_{i\in I}A_i;\bigsqcup_{j\in J}B_j) & &
\end{array}
\tag{1.1.4}
$$

Here equality of source functors and of target functors follows from (1.1.3) and (1.1.1). Similarly, an equation for right unitors has to hold:

$$
\begin{array}{ccc}
\prod_{k\in K}[\mathcal{P}(\bigsqcup_{i\in f^{-1}k}A_i;\bigsqcup_{j\in g^{-1}k}B_j)\times\mathbb{1}] & & \\
 & \searrow^{\prod_K(1\times\mathrm{i})} & \\
\prod_K r \downarrow \cong & \xLeftarrow{\quad\prod_K\mathsf{r}\quad} & \prod_{k\in K}[\mathcal{P}(\bigsqcup_{i\in f^{-1}k}A_i;\bigsqcup_{j\in g^{-1}k}B_j)\times\mathcal{P}(\bigsqcup_{j\in g^{-1}k}B_j;\bigsqcup_{j\in g^{-1}k}B_j)] \\
 & \swarrow_{\prod_K m} & \\
\prod_{k\in K}\mathcal{P}(\bigsqcup_{i\in f^{-1}k}A_i;\bigsqcup_{j\in g^{-1}k}B_j) & \xrightarrow[\otimes^{I\xrightarrow{f}K\xleftarrow{g}J}]{} & \mathcal{P}(\bigsqcup_{i\in I}A_i;\bigsqcup_{j\in J}B_j)
\end{array}
=
$$

$$
\begin{array}{ccc}
\mathcal{P}(\bigsqcup_{i\in I}A_i;\bigsqcup_{j\in J}B_j)\times\mathbb{1} & \xleftarrow{\;r^{-1}\circ\otimes^{I\xrightarrow{f}K\xleftarrow{g}J}\circ\prod_K r\;} & \prod_{k\in K}[\mathcal{P}(\bigsqcup_{i\in f^{-1}k}A_i;\bigsqcup_{j\in g^{-1}k}B_j)\times\mathbb{1}] \\
 & \searrow^{1\times\mathrm{i}} & \\
r \downarrow \cong & \xLeftarrow{\quad\mathsf{r}\quad} & \mathcal{P}(\bigsqcup_{i\in I}A_i;\bigsqcup_{j\in J}B_j)\times\mathcal{P}(\bigsqcup_{j\in J}B_j;\bigsqcup_{j\in J}B_j) \\
 & \swarrow_{m} & \\
\mathcal{P}(\bigsqcup_{i\in I}A_i;\bigsqcup_{j\in J}B_j) & &
\end{array}
\tag{1.1.5}
$$

— for any commutative diagram in $\mathcal{S}_{\mathsf{sk}}$

$$
\begin{array}{ccccc}
I & \xrightarrow{\;f\;} & K & \xleftarrow{\;g\;} & J \\
 & {}_{\alpha}\searrow \;=\; & \downarrow h & \;=\; \swarrow_{\beta} & \\
 & & L & &
\end{array}
$$

there is an equality (or, rather, a symmetry in $\mathcal{S}et$) of functors

$$
\begin{array}{ccc}
\prod_{k\in K}\mathcal{P}(\bigsqcup_{i\in f^{-1}k}A_i;\bigsqcup_{j\in g^{-1}k}B_j) & \xrightarrow{\;\otimes^{I\xrightarrow{f}K\xleftarrow{g}J}\;} & \mathcal{P}(\bigsqcup_{i\in I}A_i;\bigsqcup_{j\in J}B_j) \\
\cong\downarrow & = & \uparrow\otimes^{I\xrightarrow{\alpha}L\xleftarrow{\beta}J} \\
\prod_{l\in L}\prod_{k\in h^{-1}l}\mathcal{P}(\bigsqcup_{i\in f^{-1}k}A_i;\bigsqcup_{j\in g^{-1}k}B_j) & \xrightarrow[\prod_{l\in L}\otimes^{\alpha^{-1}l\xrightarrow{f|}h^{-1}l\xleftarrow{g|}\beta^{-1}l}]{} & \prod_{l\in L}\mathcal{P}(\bigsqcup_{i\in\alpha^{-1}l}A_i;\bigsqcup_{j\in\beta^{-1}l}B_j)
\end{array}
\tag{1.1.6}
$$

We refer to this property as to strictness of tensor multiplication.

Let $\alpha,\beta,\gamma$ from $I\xrightarrow{\alpha}J\xrightarrow{\beta}K\xrightarrow{\gamma}L\in\mathcal{S}_{\mathsf{sk}}$ be bijections. Let $(A'_l)_{l\in L}$, $(B'_l)_{l\in L}$, $(C'_l)_{l\in L}$, $(D'_l)_{l\in L}$ be families of objects of $\mathcal{P}$. Apply axiom (1.1.2) to maps

$$
\begin{array}{ccccc}
I & \searrow^{\alpha\beta\gamma} & & {}^{\beta\gamma}\swarrow & J \\
 & & L & & \\
L & \nearrow\!\!\!= & & \nwarrow_{\gamma} & K
\end{array}
$$

and to objects $A_i=A'_{\gamma\beta\alpha i}$, $i\in I$, $B_j=B'_{\gamma\beta j}$, $j\in J$, $C_k=C'_{\gamma k}$, $k\in K$, $D_l=D'_l$, $l\in L$. Then (1.1.2)

implies the identity

$$\prod_{l\in L}\mathcal{P}(A_l;B_l)\times\prod_{l\in L}\mathcal{P}(B_l;C_l)\times\prod_{l\in L}\mathcal{P}(C_l;D_l) \tag{1.1.7}$$

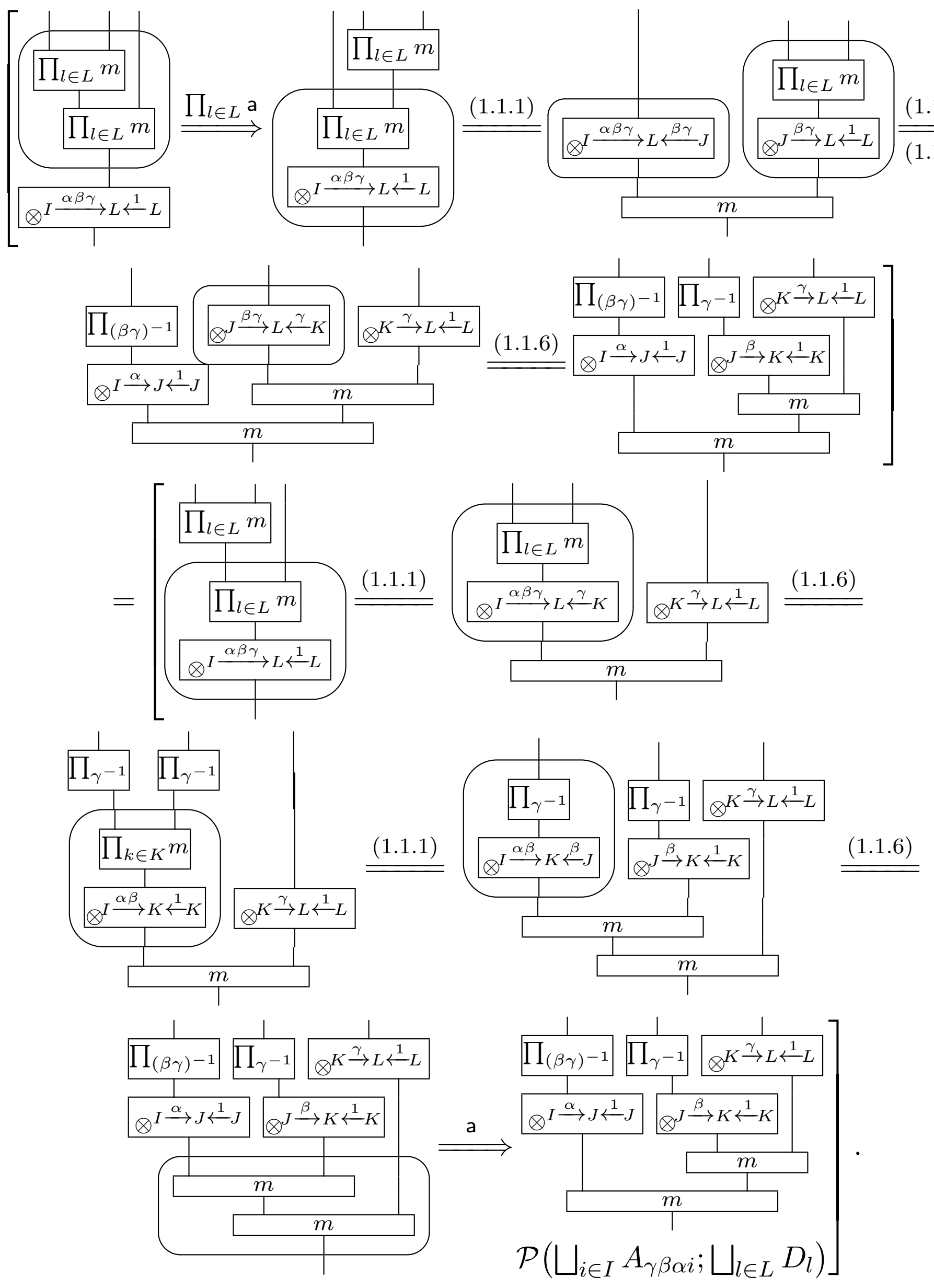

**1.1.2 Proposition.** *Let $\alpha,\beta,\gamma$ from $I\xrightarrow{\alpha}J\xrightarrow{\beta}K\xrightarrow{\gamma}L\in\mathcal{S}_{\mathsf{sk}}$ be bijections. Let $(X_i)_{i\in I}$, $(Y_j)_{j\in J}$, $(Z_k)_{k\in K}$, $(W_l)_{l\in L}$ be families of colours of $\mathcal{P}$ such that $X_i=Y_{\alpha i}$, $Y_j=Z_{\beta j}$, $Z_k=W_{\gamma k}$ for all $i\in I$, $j\in J$, $k\in K$. Then*

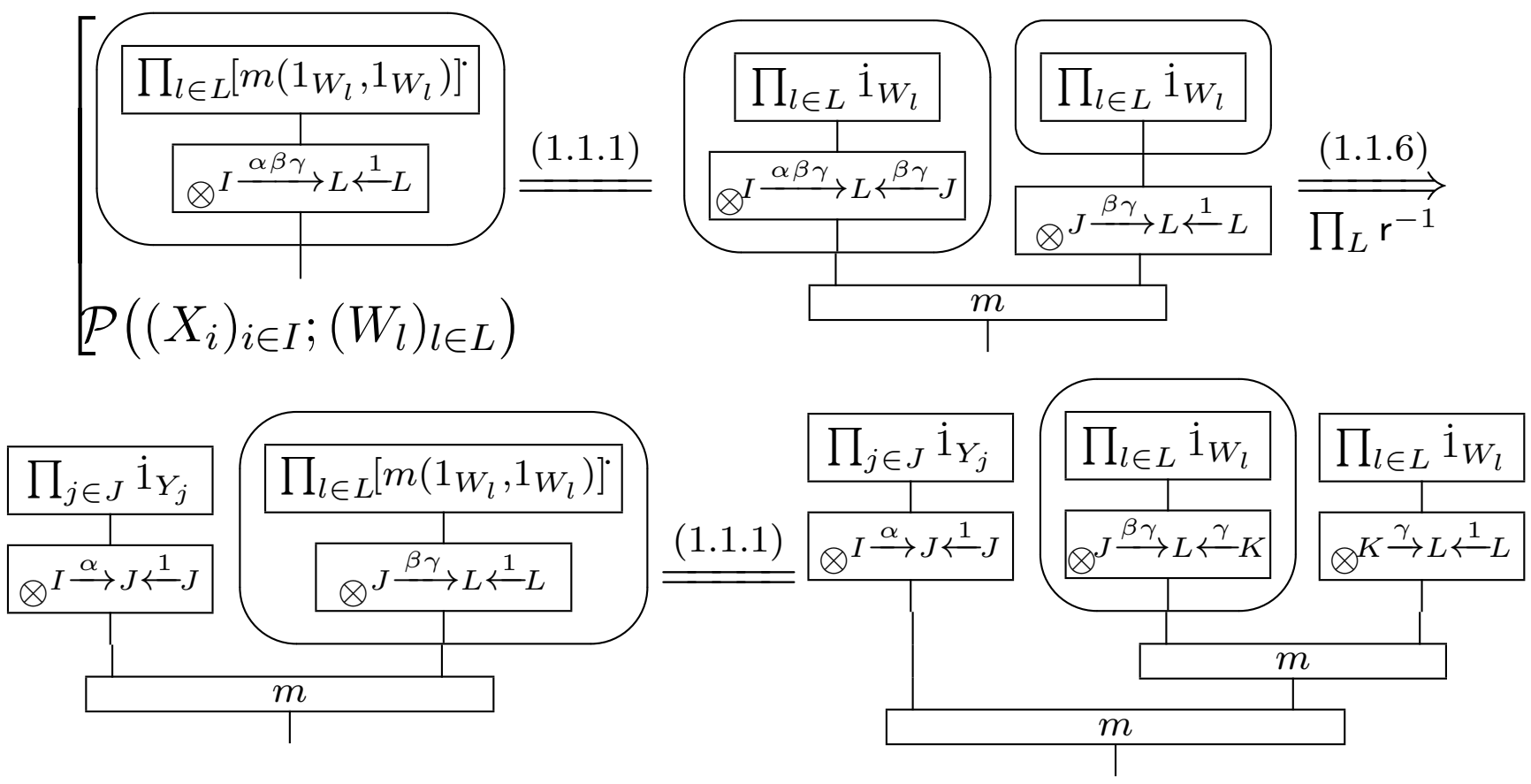

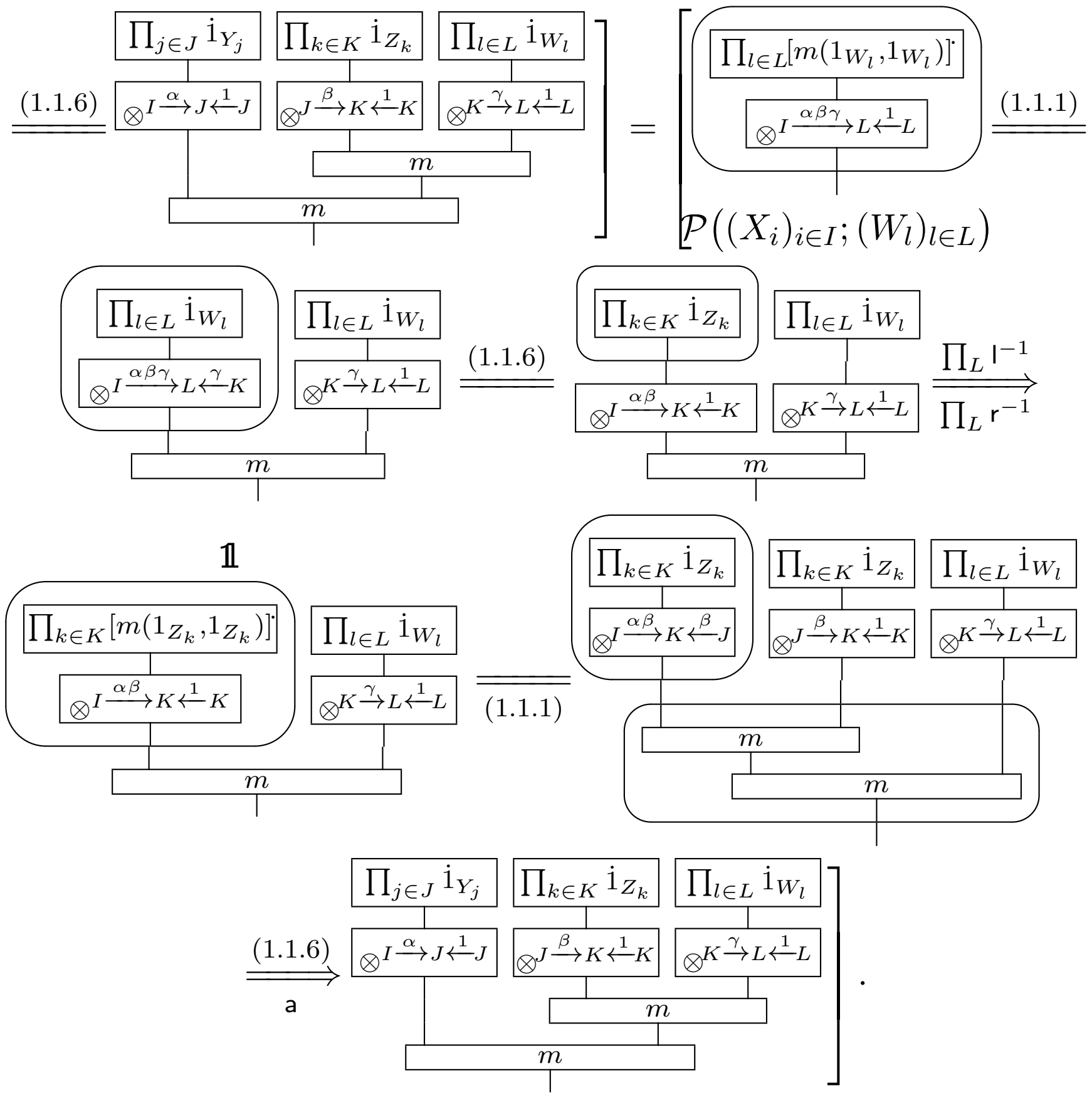

(1.1.8)

*Proof.* Follows from identity (1.1.7) precomposed with

$$\boxed{\textstyle\prod_{l\in L}\mathrm{i}_{W_l}}\;\boxed{\textstyle\prod_{l\in L}\mathrm{i}_{W_l}}\;\boxed{\textstyle\prod_{l\in L}\mathrm{i}_{W_l}}\colon \mathbb{1}\to\prod_{l\in L}\mathcal{P}(W_l;W_l)\times\prod_{l\in L}\mathcal{P}(W_l;W_l)\times\prod_{l\in L}\mathcal{P}(W_l;W_l).$$

□

**1.1.3 Example.** An example of biprop $\mathfrak{C}$ comes from the symmetric monoidal 2-category $(\mathcal{C}at,\times)$.

— $\operatorname{Col}\mathfrak{C}$ = class of small categories. Hence, objects of $\mathfrak{C}$ are tuples $(\mathcal{C}_i)_{i\in I}$ of small categories;

— $\mathfrak{C}\big((\mathcal{C}_i)_{i\in I};(\mathcal{D}_j)_{j\in J}\big)=\underline{\mathcal{C}at}(\prod_{i\in I}\mathcal{C}_i,\prod_{j\in J}\mathcal{D}_j)$ = category of functors $\prod_{i\in I}\mathcal{C}_i\to\prod_{j\in J}\mathcal{D}_j$ and of natural transformations of such for each tuple $\big((\mathcal{C}_i)_{i\in I},(\mathcal{D}_j)_{j\in J}\big)$ of small categories. Note that $\mathcal{C}at$ is a closed monoidal category and the internal homs $\underline{\mathcal{C}at}$ form a $\mathcal{C}at$-category;

— the composition in $\mathfrak{C}$

$$m\colon\underline{\mathcal{C}at}(\prod_{i\in I}\mathcal{C}_i,\prod_{j\in J}\mathcal{D}_j)\times\underline{\mathcal{C}at}(\prod_{j\in J}\mathcal{D}_j,\prod_{k\in K}\mathcal{E}_k)\to\underline{\mathcal{C}at}(\prod_{i\in I}\mathcal{C}_i,\prod_{k\in K}\mathcal{E}_k)$$

is that of $\underline{\mathcal{C}at}$. It is strictly associative, so we choose the associator $\mathsf{a}$ to be the identity transformation;

— The unit 1-morphism of an object $(\mathcal{C}_i)_{i\in I}$ is the identity functor

$$1_{(\mathcal{C}_i)_{i\in I}}=\mathrm{Id}_{\prod_{i\in I}\mathcal{C}_i}\in\operatorname{Ob}\underline{\mathcal{C}at}(\prod_{i\in I}\mathcal{C}_i,\prod_{i\in I}\mathcal{C}_i)=\mathcal{C}at(\prod_{i\in I}\mathcal{C}_i,\prod_{i\in I}\mathcal{C}_i).$$

It is a strict unit, so the unitors $\mathsf{r}$, $\mathsf{l}$ are chosen as the identity transformation. Summing up, $\mathfrak{C}$ is also a $\mathcal{C}at$-category, that is, a 2-category;

— let $A_i=(\mathcal{C}_i^s)_{s\in S_i}$ be an object of $\mathfrak{C}$, $i\in I$, and let $f\colon I\to L$ be a map. For each $l\in L$ we have $\bigsqcup_{i\in f^{-1}l}A_i=\big((\mathcal{C}_i^s)_{s\in S_i}\big)_{i\in f^{-1}l}$. So the map $\otimes^{I\xrightarrow{f}L}\colon(\operatorname{Ob}\mathfrak{C})^L\to\operatorname{Ob}\mathfrak{C}$ on objects is chosen as $\big(\big((\mathcal{C}_i^s)_{s\in S_i}\big)_{i\in f^{-1}l}\big)_{l\in L}\mapsto\big((\mathcal{C}_i^s)_{s\in S_i}\big)_{i\in I}$. Let $B_j=(\mathcal{D}_j^q)_{q\in Q_j}$ be an object of $\mathfrak{C}$, $j\in J$, and let $g\colon J \to L$ be a map. Introduce the notation $\prod A_i = \prod_{s\in S_i} \mathcal{C}_i^s$. Similarly, $\prod B_j = \prod_{q\in Q_j} \mathcal{D}_j^q$, etc. The functor

$$\otimes^{I \xrightarrow{f} L \xleftarrow{g} J}\colon \prod_{l\in L} \mathcal{P}\Big(\bigsqcup_{i\in f^{-1}l} A_i; \bigsqcup_{j\in g^{-1}l} B_j\Big) \to \mathcal{P}\Big(\bigsqcup_{i\in I} A_i; \bigsqcup_{j\in J} B_j\Big),$$

indexed by a pair of maps $I \xrightarrow{f} L \xleftarrow{g} J \in \mathcal{S}_{\mathsf{sk}}$ is

$$\prod_{l\in L} \underline{\mathcal{C}at}\Big(\prod\big((\mathcal{C}_i^s)_{s\in S_i}\big)_{i\in f^{-1}l}, \prod\big((\mathcal{D}_j^q)_{q\in Q_j}\big)_{j\in g^{-1}l}\Big)$$
$$\to \underline{\mathcal{C}at}\Big(\prod\big((\mathcal{C}_i^s)_{s\in S_i}\big)_{i\in I}, \prod\big((\mathcal{D}_j^q)_{q\in Q_j}\big)_{j\in J}\Big), \quad (F_l)_{l\in L} \mapsto$$
$$\Big(\prod_{i\in I}\prod A_i \xrightarrow{\prod_{\sigma(f)}\prod A_i} \prod_{l\in L}\prod_{i\in f^{-1}l}\prod\prod A_i \xrightarrow{\prod_{l\in L} F_l} \prod_{l\in L}\prod_{j\in g^{-1}l}\prod\prod B_j \xrightarrow{\prod_{\sigma(g)^{-1}}\prod B_j} \prod_{j\in J}\prod B_j\Big), \tag{1.1.9}$$

where the bijection $\sigma(f)\colon I \to \bigsqcup_{l\in L} f^{-1}l, i \mapsto (i, fi) \in \bigsqcup_L I = I \times L$ is the graph of $f$ (see [Lyu23, (2.1.1) and above] and [Lyu25, Example 1.1.2]). Of course, the bijection $\sigma(g)\colon J \to \bigsqcup_{l\in L} g^{-1}l, j \mapsto (j, gj) \in \bigsqcup_L J = J \times L$ is the graph of $g$. When $f\colon I \to L$ is a bijection we may identify $\bigsqcup_{l\in L} f^{-1}l$ with $L$ and the bijection $\sigma(f)$ identifies with $f$.

Suppose that $(F_l)_{l\in L}$, $(G_l)_{l\in L}$ are objects (or morphisms) from the left top corner of (1.1.1):

$$\prod_{l\in L}\underline{\mathcal{C}at}\Big(\prod_{i\in f^{-1}l}\prod A_i, \prod_{j\in g^{-1}l}\prod B_j\Big) \times \prod_{l\in L}\underline{\mathcal{C}at}\Big(\prod_{j\in g^{-1}l}\prod B_j, \prod_{k\in h^{-1}l}\prod C_k\Big) \xrightarrow{\prod_L m} \prod_{l\in L}\underline{\mathcal{C}at}\Big(\prod_{i\in f^{-1}l}\prod A_i, \prod_{k\in h^{-1}l}\prod C_k\Big)$$
$$\Big\downarrow \otimes^{I \xrightarrow{f} L \xleftarrow{g} J} \times \otimes^{J \xrightarrow{g} L \xleftarrow{h} K} \qquad\qquad \Big\downarrow \otimes^{I \xrightarrow{f} L \xleftarrow{h} K}$$
$$\underline{\mathcal{C}at}\Big(\prod_{i\in I}\prod A_i, \prod_{j\in J}\prod B_j\Big) \times \underline{\mathcal{C}at}\Big(\prod_{j\in J}\prod B_j, \prod_{k\in K}\prod C_k\Big) \xrightarrow{m} \underline{\mathcal{C}at}\Big(\prod_{i\in I}\prod A_i, \prod_{k\in K}\prod C_k\Big)$$

Both sides of the above applied to $(F_l)_{l\in L}, (G_l)_{l\in L}$ give

$$\Big(\prod_{i\in I}\prod A_i \xrightarrow{\prod_{\sigma(f)}\prod A_i} \prod_{l\in L}\prod_{i\in f^{-1}l}\prod\prod A_i \xrightarrow{\prod_{l\in L}(F_l\cdot G_l)} \prod_{l\in L}\prod_{k\in h^{-1}l}\prod\prod C_k \xrightarrow{\prod_{\sigma(h)^{-1}}\prod C_k} \prod_{k\in K}\prod\prod C_k\Big),$$

which proves equation (1.1.1).

Equation (1.1.2) holds for $\mathfrak{C}$ because both sides are identity transformations.

For $\triangledown\colon I \to \mathbf{1}$ the bijection $\big(I \xrightarrow{\sigma(\triangledown)} \bigsqcup_{\mathbf{1}} I \xrightarrow{\cong} I\big)$ is the identity map, therefore,

$$\otimes^{I \xrightarrow{\triangledown} \mathbf{1} \xleftarrow{\triangledown} J}\colon \mathfrak{C}\Big(\bigsqcup_{i\in I} A_i; \bigsqcup_{j\in J} B_j\Big) \to \mathfrak{C}\Big(\bigsqcup_{i\in I} A_i; \bigsqcup_{j\in J} B_j\Big)$$

is the identity functor.

The functor

$$\otimes^{I \xrightarrow{f} L \xleftarrow{f} I}\colon \prod_{l\in L}\underline{\mathcal{C}at}\Big(\prod\big((\mathcal{C}_i^s)_{s\in S_i}\big)_{i\in f^{-1}l}, \prod\big((\mathcal{C}_i^s)_{s\in S_i}\big)_{i\in f^{-1}l}\Big)$$
$$\to \underline{\mathcal{C}at}\Big(\prod\big((\mathcal{C}_i^s)_{s\in S_i}\big)_{i\in I}, \prod\big((\mathcal{C}_i^s)_{s\in S_i}\big)_{i\in I}\Big)$$

takes $(\mathrm{Id})_L$ to Id as (1.1.9) shows. This proves (1.1.3) for $\mathfrak{C}$.

Equations (1.1.4), (1.1.5) hold for $\mathfrak{C}$ because both sides are identity transformations.

Let us prove now identity (1.1.6) for $\mathfrak{C}$. The left-bottom-right path of (1.1.6) is

$$\prod_{k\in K}\underline{\mathcal{C}at}\Big(\prod_{i\in f^{-1}k}\prod A_i,\prod_{j\in g^{-1}k}\prod B_j\Big)\xrightarrow{\cong}\prod_{l\in L}\prod_{k\in h^{-1}l}\underline{\mathcal{C}at}\Big(\prod_{i\in f^{-1}k}\prod A_i,\prod_{j\in g^{-1}k}\prod B_j\Big)$$
$$\xrightarrow{\prod_{l\in L}\otimes^{\alpha^{-1}l\xrightarrow{f|}h^{-1}l\xleftarrow{g|}\beta^{-1}l}}$$
$$\prod_{l\in L}\underline{\mathcal{C}at}\Big(\prod_{i\in\alpha^{-1}l}\prod A_i,\prod_{j\in\beta^{-1}l}\prod B_j\Big)\xrightarrow{\otimes^{I\xrightarrow{\alpha}L\xleftarrow{\beta}J}}\underline{\mathcal{C}at}\Big(\prod_{i\in I}\prod A_i,\prod_{j\in J}\prod B_j\Big).$$

It sends (an object or morphism) $(F_k)_{k\in K}$ to

$$\big((F_k)_{k\in h^{-1}l}\big)_{l\in L}\mapsto\Big(\prod_{i\in\alpha^{-1}l}\prod A_i\xrightarrow{\prod_{\sigma(f|_{\alpha^{-1}l})}\prod A_i}\prod_{k\in h^{-1}l}\prod_{i\in f^{-1}k}\prod A_i\xrightarrow{\prod_{k\in h^{-1}l}F_k}$$
$$\prod_{k\in h^{-1}l}\prod_{j\in g^{-1}k}\prod B_j\xrightarrow{\prod_{\sigma(g|_{\beta^{-1}l})^{-1}}\prod B_j}\prod_{j\in\beta^{-1}l}\prod B_j\Big)_{l\in L}\mapsto\Big(\prod_{i\in I}\prod A_i\xrightarrow{\prod_{\sigma(\alpha)}\prod A_i}$$
$$\prod_{l\in L}\prod_{i\in\alpha^{-1}l}\prod A_i\xrightarrow{\prod_{l\in L}\prod_{\sigma(f|_{\alpha^{-1}l})}\prod A_i}\prod_{l\in L}\prod_{k\in h^{-1}l}\prod_{i\in f^{-1}k}\prod A_i\xrightarrow{\prod_{l\in L}\prod_{k\in h^{-1}l}F_k}$$
$$\prod_{l\in L}\prod_{k\in h^{-1}l}\prod_{j\in g^{-1}k}\prod B_j\xrightarrow{\prod_{l\in L}\prod_{\sigma(g|_{\beta^{-1}l})^{-1}}\prod B_j}\prod_{l\in L}\prod_{j\in\beta^{-1}l}\prod B_j\xrightarrow{\prod_{\sigma(\beta)^{-1}}\prod B_j}\prod_{j\in J}\prod B_j\Big).\tag{1.1.10}$$

The top arrow of (1.1.6) is

$$\otimes^{I\xrightarrow{f}K\xleftarrow{g}J}:\prod_{k\in K}\underline{\mathcal{C}at}\Big(\prod_{i\in f^{-1}k}\prod A_i,\prod_{j\in g^{-1}k}\prod B_j\Big)\to\underline{\mathcal{C}at}\Big(\prod_{i\in I}\prod A_i,\prod_{j\in J}\prod B_j\Big).$$

It sends $(F_k)_{k\in K}$ to

$$\Big(\prod_{i\in I}\prod A_i\xrightarrow{\prod_{\sigma(f)}\prod A_i}\prod_{k\in K}\prod_{i\in f^{-1}k}\prod A_i\xrightarrow{\prod_{k\in K}F_k}\prod_{k\in K}\prod_{j\in g^{-1}k}\prod B_j\xrightarrow{\prod_{\sigma(g)^{-1}}\prod B_j}\prod_{j\in J}\prod B_j\Big),$$

This expression coincides with the result of (1.1.10) due to the following equality in $\mathcal{S}_{\mathsf{sk}}$ for arbitrary two composable maps $\alpha=\big(I\xrightarrow{f}K\xrightarrow{h}L\big)\in\mathcal{S}_{\mathsf{sk}}$

$$\begin{array}{ccc} I & \xrightarrow{\sigma(\alpha)} & \bigsqcup_{l\in L}\alpha^{-1}l \\ {\scriptstyle\sigma(f)}\downarrow & = & \downarrow{\scriptstyle\sqcup_{l\in L}\sigma(f|_{\alpha^{-1}l})} \\ \bigsqcup_{k\in K}f^{-1}k & \xrightarrow[\sqcup_{\sigma(h)}f^{-1}k]{} & \bigsqcup_{l\in L}\bigsqcup_{k\in h^{-1}l}f^{-1}k \end{array}\tag{1.1.11}$$

Equivalently we may substitute $f\centerdot h$ for $\alpha$:

$$\begin{array}{ccc} I & \xrightarrow{\sigma(f\centerdot h)} & \bigsqcup_{l\in L}f^{-1}h^{-1}l \\ {\scriptstyle\sigma(f)}\downarrow & = & \downarrow{\scriptstyle\sqcup_{l\in L}\sigma(f|_{f^{-1}h^{-1}l})} \\ \bigsqcup_{k\in K}f^{-1}k & \xrightarrow[\sqcup_{\sigma(h)}f^{-1}k]{} & \bigsqcup_{l\in L}\bigsqcup_{k\in h^{-1}l}f^{-1}k \end{array}\tag{1.1.12}$$

On elements we have the equality

$$\begin{array}{ccc} i & \longmapsto & i\text{ indexed by }hfi \\ \downarrow & = & \downarrow \\ i\text{ indexed by }fi & \longmapsto & (i\text{ indexed by }fi)\text{ indexed by }hfi \end{array}$$

hence, (1.1.12) holds true, equivalently, (1.1.11) holds true. By this proof a similar equality holds for $g$ and $\beta$.

Note that equality (1.1.7) holds true because the both hand sides are identity transformations.

## 1.2 Biprop envelope of a symmetric weak multicategory

Although enrichment over a symmetric weak multicategory generalises enrichment over a biprop, one can work over a symmetric weak multicategory as one does over a biprop with the help of the following construction. The *biprop envelope* of a symmetric weak multicategory $\mathsf{V}$ is the biprop $\mathbb{F}\mathsf{V}$ whose objects $(X_1,\dots,X_n)$ are words of objects of $\mathsf{V}$ of length $n \geqslant 0$, and in which the category $(\mathbb{F}\mathsf{V})\big((X_1,\dots,X_n);(Y_1,\dots,Y_m)\big)$ is the category

$$\coprod_{\varphi\colon \{1,\dots,n\}\to\{1,\dots,m\}} \prod_{1\leqslant j\leqslant m} \mathsf{V}\big((X_i)_{i\in\varphi^{-1}(j)};Y_j\big)$$

composition in $\mathbb{F}\mathsf{V}$ is defined using the symmetric weak multicategory structure of $\mathsf{V}$, and the tensor product and symmetry are given by concatenation and permutation of words. This construction defines a functor from the category of symmetric weak multicategories and symmetric multifunctors to the category of biprops and their morphisms.

The non-symmetric version of this construction, namely the monoidal envelope of a multicategory, was defined in [Lin71], where it was used (among other things) to define the notions of category, functor, and natural transformation enriched over a multicategory in terms of the same notions enriched over a strict monoidal category. We make the analogous definitions in the symmetric weak case.

Let $\mathsf{C}$ be a symmetric weak multicategory [Lyu25, Definition 1.1.1]. The biprop $\mathbb{F}\mathsf{C}$ has

— $\operatorname{Ob}\mathbb{F}\mathsf{C} = (\operatorname{Ob}\mathsf{C})^*$, the free monoid on the set of colours $\operatorname{Col}\mathbb{F}\mathsf{C} = \operatorname{Ob}\mathsf{C}$;

— $\mathbb{F}\mathsf{C}\big((X_i)_{i\in I},(Y_j)_{j\in J}\big) = \coprod_{\phi\colon I\to J\in\mathcal{S}_{\mathsf{sk}}}\prod_{j\in J}\mathsf{C}\big((X_i)_{i\in\phi^{-1}j};Y_j\big)$;

— the composition $m$ is

$$\begin{aligned}
&\mathbb{F}\mathsf{C}\big((X_i)_{i\in I},(Y_j)_{j\in J}\big)\times\mathbb{F}\mathsf{C}\big((Y_j)_{j\in J},(Z_k)_{k\in K}\big)\\
&\qquad\cong \coprod_{I\xrightarrow{\phi}J\xrightarrow{\psi}K\in\mathcal{S}_{\mathsf{sk}}}\prod_{k\in K}\Big[\Big(\prod_{j\in\psi^{-1}k}\mathsf{C}\big((X_i)_{i\in\phi^{-1}j};Y_j\big)\Big)\times\mathsf{C}\big((Y_j)_{j\in\psi^{-1}k};Z_k\big)\Big]\\
&\qquad\xrightarrow{\coprod_{(\phi,\psi)\mapsto\phi\cdot\psi}\prod_{k\in K}\mu_{\phi|\colon\phi^{-1}\psi^{-1}k\to\psi^{-1}k}} \coprod_{\xi\colon I\to K\in\mathcal{S}_{\mathsf{sk}}}\prod_{k\in K}\mathsf{C}\big((X_i)_{i\in\xi^{-1}k};Z_k\big)\\
&\qquad= \mathbb{F}\mathsf{C}\big((X_i)_{i\in I},(Z_k)_{k\in K}\big).
\end{aligned}$$

— the identity morphism in $\mathbb{F}\mathsf{C}\big((X_i)_{i\in I},(X_i)_{i\in I}\big)$ is the object

$$(1_{X_i})_{i\in I}\in\prod_{i\in I}\operatorname{Ob}\mathsf{C}(X_i;X_i),$$

indexed by the identity map $\mathrm{id}_I$;

— the tensor multiplication on objects is the concatenation. Consider objects $A_i = (X_i^s)_{s\in S_i}$, $B_j = (Y_j^q)_{q\in Q_j}$. Introduce ordered sets $S = \bigsqcup_{i\in I} S_i$, $Q = \bigsqcup_{j\in J} Q_j$. A map $\phi : S \to Q \in \mathcal{S}_{\mathsf{sk}}$ such that the diagram

$$\begin{array}{ccccc}
S & & \xrightarrow{\quad\phi\quad} & & Q\\
{\scriptstyle \mathrm{pr}_I}\big\downarrow & & & & \big\downarrow{\scriptstyle \mathrm{pr}_J}\\
I & \xrightarrow{\ f\ } & K & \xleftarrow{\ g\ } & J
\end{array}\tag{1.2.1}$$

commutes amounts to family of maps $\phi_k\colon \bigsqcup_{i\in f^{-1}k} S_i \to \bigsqcup_{j\in g^{-1}k} Q_j \in \mathcal{S}_{\mathsf{sk}}$, $k\in K$. On categories of morphisms the tensor multiplication $\otimes^{I\xrightarrow{f}K\xleftarrow{g}J}$ is the functor

$$\otimes^{I\xrightarrow{f}K\xleftarrow{g}J}\colon \prod_{k\in K} \mathbb{F}\mathsf{C}\Big(\bigsqcup_{i\in f^{-1}k} A_i;\ \bigsqcup_{j\in g^{-1}k} B_j\Big)$$

$$= \prod_{k\in K} \mathbb{F}\mathsf{C}\big((X^s_{\mathrm{pr}_I s})_{s\in\sqcup_{i\in f^{-1}k}S_i};(Y^q_{\mathrm{pr}_J q})_{q\in\sqcup_{j\in g^{-1}k}Q_j}\big)$$

$$= \prod_{k\in K}\ \coprod_{\phi_k\colon \sqcup_{i\in f^{-1}k}S_i\to\sqcup_{j\in g^{-1}k}Q_j}\ \prod_{q\in\sqcup_{j\in g^{-1}k}Q_j} \mathsf{C}\big((X^s_{\mathrm{pr}_I s})_{s\in\phi_k^{-1}q};Y^q_{\mathrm{pr}_J q}\big)$$

$$\cong \coprod_{(\phi_k\colon \sqcup_{i\in f^{-1}k}S_i\to\sqcup_{j\in g^{-1}k}Q_j)_{k\in K}}\ \prod_{k\in K}\ \prod_{q\in\sqcup_{j\in g^{-1}k}Q_j} \mathsf{C}\big((X^s_{\mathrm{pr}_I s})_{s\in\phi_k^{-1}q};Y^q_{\mathrm{pr}_J q}\big)$$

$$= \coprod_{(\phi_k\colon \sqcup_{i\in f^{-1}k}S_i\to\sqcup_{j\in g^{-1}k}Q_j)_{k\in K}}\ \prod_{k\in K}\ \prod_{j\in g^{-1}k}\ \prod_{q\in Q_j} \mathsf{C}\big((X^s_{\mathrm{pr}_I s})_{s\in\phi_k^{-1}q};Y^q_j\big) \xhookrightarrow{\coprod_{(\phi_k)\mapsto\phi} 1}$$

$$\coprod_{\xi\colon S\to Q}\ \prod_{j\in J}\ \prod_{q\in Q_j} \mathsf{C}\big((X^s_{\mathrm{pr}_I s})_{s\in\xi^{-1}q};Y^q_j\big) = \mathbb{F}\mathsf{C}\big((X^s_{\mathrm{pr}_I s})_{s\in S},(Y^q_{\mathrm{pr}_J q})_{q\in Q}\big)$$

$$= \mathbb{F}\mathsf{C}\Big(\bigsqcup_{i\in I} A_i,\bigsqcup_{j\in J} B_j\Big), \quad (1.2.2)$$

where $A_i = (X^s_i)_{s\in S_i}$, $B_j = (Y^q_j)_{q\in Q_j}$, and $\phi\colon S\to Q$ is the only map, which makes (1.2.1) commutative and satisfies the condition $\big(\sqcup_{i\in f^{-1}k}S_i \hookrightarrow S \xrightarrow{\phi} Q\big) = \big(\sqcup_{i\in f^{-1}k}S_i \xrightarrow{\phi_k} \bigsqcup_{j\in g^{-1}k}Q_j \hookrightarrow Q\big)$. We shall see that the tensor multiplication is strictly associative.

Let us prove that diagram (1.1.1) with $\mathcal{P}=\mathbb{F}\mathsf{C}$, $A_i = (X^s_i)_{s\in S_i}$, $B_j = (Y^q_j)_{q\in Q_j}$, $C_k = (Z^p_k)_{p\in P_k}$ commutes. Denote $S=\bigsqcup_{i\in I}S_i$, $Q=\bigsqcup_{j\in J}Q_j$, $P=\bigsqcup_{k\in K}P_k$. We consider pairs of maps $S\xrightarrow{\phi}Q\xrightarrow{\psi}P$ such that the bottom of the diagram

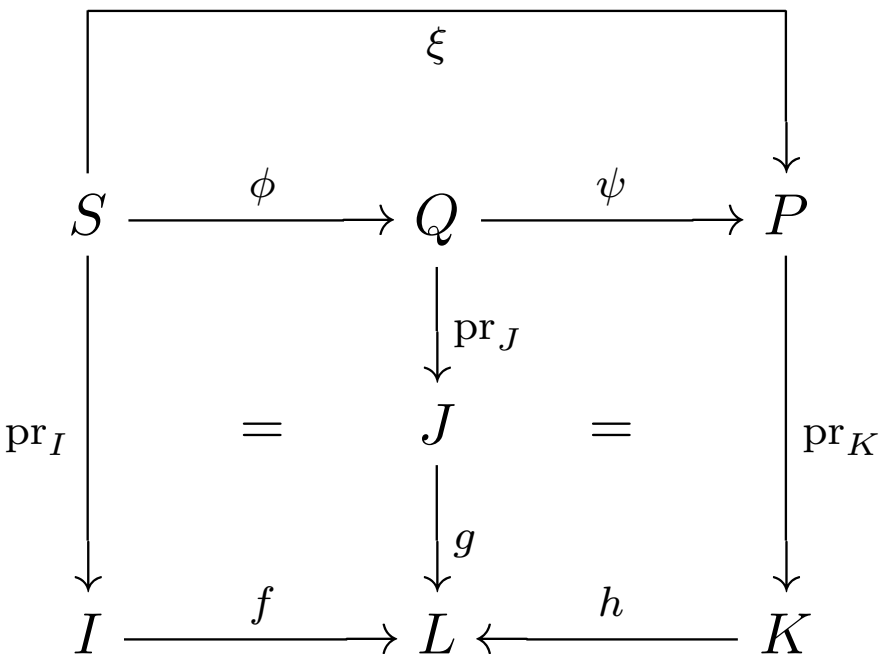

(1.2.3)

commutes.

The top–right path of diagram (1.1.1) is

$$\prod_{l\in L}\mathbb{F}\mathsf{C}\big(((X^s_i)_{s\in S_i})_{i\in f^{-1}l};((Y^q_j)_{q\in Q_j})_{j\in g^{-1}l}\big)$$

$$\times\prod_{l\in L}\mathbb{F}\mathsf{C}\big(((Y^q_j)_{q\in Q_j})_{j\in g^{-1}l};((Z^p_k)_{p\in P_k})_{k\in h^{-1}l}\big)\cong$$

$$\prod_{l\in L}\big[\mathbb{F}\mathsf{C}\big((X^s_{\mathrm{pr}_I s})_{s\in\sqcup_{i\in f^{-1}l}S_i};(Y^q_{\mathrm{pr}_J q})_{q\in\sqcup_{j\in g^{-1}l}Q_j}\big)$$

$$\times\mathbb{F}\mathsf{C}\big((Y^q_{\mathrm{pr}_J q})_{q\in\sqcup_{j\in g^{-1}l}Q_j};(Z^p_{\mathrm{pr}_K p})_{p\in\sqcup_{p\in h^{-1}l}P_k}\big)\big]$$

$$\cong\prod_{l\in L}\ \coprod_{\sqcup_{i\in f^{-1}l}S_i\xrightarrow{\phi_l}\sqcup_{j\in g^{-1}l}Q_j\xrightarrow{\psi_l}\sqcup_{k\in h^{-1}l}P_k}\ \prod_{p\in\sqcup_{k\in h^{-1}l}P_k}\Big[\Big(\prod_{q\in\psi_l^{-1}p}\mathsf{C}\big((X^s_{\mathrm{pr}_I s})_{s\in\phi_l^{-1}q};Y^q_{\mathrm{pr}_J q}\big)\Big)$$

$$\times\mathsf{C}\big((Y^q_{\mathrm{pr}_J q})_{q\in\psi_l^{-1}p};Z^p_{\mathrm{pr}_K p}\big)\Big]$$

$$\xrightarrow{\prod_{l\in L}\coprod_{(\phi_l,\psi_l)\mapsto\phi_l\centerdot\psi_l}\prod_{p\in\sqcup_{k\in h^{-1}l}P_k}\mu_{\phi_l|:\,\phi_l^{-1}\psi_l^{-1}p\to\psi_l^{-1}p}}$$

$$\prod_{l\in L}\coprod_{\xi_l:\,\sqcup_{i\in f^{-1}l}S_i\to\sqcup_{k\in h^{-1}l}P_k}\prod_{p\in\sqcup_{k\in h^{-1}l}P_k}\mathsf{C}\big((X^s_{\mathrm{pr}_I s})_{s\in\xi_l^{-1}p};Z^p_{\mathrm{pr}_K p}\big)\cong$$

$$\coprod_{(\xi_l:\,\sqcup_{i\in f^{-1}l}S_i\to\sqcup_{k\in h^{-1}l}P_k)_{l\in L}}\prod_{l\in L}\prod_{p\in\sqcup_{k\in h^{-1}l}P_k}\mathsf{C}\big((X^s_{\mathrm{pr}_I s})_{s\in\xi_l^{-1}p};Z^p_{\mathrm{pr}_K p}\big)$$

$$\xhookrightarrow{\coprod_{(\xi_l)\mapsto\xi}1}\coprod_{\xi:\,S\to P}\prod_{p\in P}\mathsf{C}\big((X^s_{\mathrm{pr}_I s})_{s\in\xi^{-1}p};Z^p_{\mathrm{pr}_K p}\big)=\mathbb{F}\mathsf{C}\big((X^s_{\mathrm{pr}_I s})_{s\in S};(Z^p_{\mathrm{pr}_K p})_{p\in P}\big). \quad (1.2.4)$$

Here $\phi\colon S\to Q$ (resp. $\psi\colon Q\to P$, $\xi\colon S\to P$) is the only map, which satisfies the condition $\phi|_{\sqcup_{i\in f^{-1}l}S_i}=\phi_l$ (resp. $\psi|_{\sqcup_{j\in g^{-1}l}Q_j}=\psi_l$, $\xi|_{\sqcup_{i\in f^{-1}l}S_i}=\xi_l$). This calculation can be written in this form since for $p\in\bigsqcup_{k\in h^{-1}l}P_k$ we have $(\phi|\colon\phi^{-1}\psi^{-1}p\to\psi^{-1}p)=(\phi_l|\colon\phi_l^{-1}\psi_l^{-1}p\to\psi_l^{-1}p)$.

The left–bottom path of diagram (1.1.1) is

$$\prod_{l\in L}\mathbb{F}\mathsf{C}\big(\big((X^s_i)_{s\in S_i}\big)_{i\in f^{-1}l};\big((Y^q_j)_{q\in Q_j}\big)_{j\in g^{-1}l}\big)$$

$$\times\prod_{l\in L}\mathbb{F}\mathsf{C}\big(\big((Y^q_j)_{q\in Q_j}\big)_{j\in g^{-1}l};\big((Z^p_k)_{p\in P_k}\big)_{k\in h^{-1}l}\big)\cong$$

$$\Big[\coprod_{(\phi_l:\,\sqcup_{i\in f^{-1}l}S_i\to\sqcup_{j\in g^{-1}l}Q_j)_{l\in L}}\prod_{l\in L}\prod_{q\in\sqcup_{j\in g^{-1}l}Q_j}\mathsf{C}\big((X^s_{\mathrm{pr}_I s})_{s\in\phi_l^{-1}q};Y^q_{\mathrm{pr}_J q}\big)\Big]\times$$

$$\coprod_{(\psi_l:\,\sqcup_{j\in g^{-1}l}Q_j\to\sqcup_{k\in h^{-1}l}P_k)_{l\in L}}\prod_{l\in L}\prod_{p\in\sqcup_{k\in h^{-1}l}P_k}\mathsf{C}\big((Y^q_{\mathrm{pr}_J q})_{q\in\psi_l^{-1}p};Z^p_{\mathrm{pr}_K p}\big)\xhookrightarrow{\coprod_{(\phi_l)\mapsto\phi}1\times\coprod_{(\psi_l)\mapsto\psi}1}$$

$$\Big[\coprod_{\chi:\,S\to Q}\prod_{q\in Q}\mathsf{C}\big((X^s_{\mathrm{pr}_I s})_{s\in\chi^{-1}q};Y^q_{\mathrm{pr}_J q}\big)\Big]\times\coprod_{\zeta:\,Q\to P}\prod_{p\in P}\mathsf{C}\big((Y^q_{\mathrm{pr}_J q})_{q\in\zeta^{-1}p};Z^p_{\mathrm{pr}_K p}\big)$$

$$\cong\coprod_{S\xrightarrow{\chi}Q\xrightarrow{\zeta}P}\prod_{p\in P}\Big[\Big(\prod_{q\in\zeta^{-1}p}\mathsf{C}\big((X^s_{\mathrm{pr}_I s})_{s\in\chi^{-1}q};Y^q_{\mathrm{pr}_J q}\big)\Big)\times\mathsf{C}\big((Y^q_{\mathrm{pr}_J q})_{q\in\zeta^{-1}p};Z^p_{\mathrm{pr}_K p}\big)\Big]$$

$$\xrightarrow{\coprod_{(\chi,\zeta)\mapsto\chi\centerdot\zeta}\prod_{p\in P}\mu_{\chi|:\,\chi^{-1}\zeta^{-1}p\to\zeta^{-1}p}}\coprod_{\xi:\,S\to P}\prod_{p\in P}\mathsf{C}\big((X^s_{\mathrm{pr}_I s})_{s\in\xi^{-1}p};Z^p_{\mathrm{pr}_K p}\big)$$

$$=\mathbb{F}\mathsf{C}\big((X^s_{\mathrm{pr}_I s})_{s\in S};(Z^p_{\mathrm{pr}_K p})_{p\in P}\big).$$

This equals (1.2.4).

— the unit object is the empty sequence $()=()_\varnothing$;

The associator for composition

$$\mathsf{a}\colon(m\times1)\centerdot m\to(1\times m)\centerdot m\colon$$
$$\mathbb{F}\mathsf{C}\big((X_i)_{i\in I},(Y_j)_{j\in J}\big)\times\mathbb{F}\mathsf{C}\big((Y_j)_{j\in J},(Z_k)_{k\in K}\big)\times\mathbb{F}\mathsf{C}\big((Z_k)_{k\in K},(W_l)_{l\in L}\big)$$
$$\to\mathbb{F}\mathsf{C}\big((X_i)_{i\in I},(W_l)_{l\in L}\big)$$

is given by

$$\coprod_{I\xrightarrow{\phi}J\xrightarrow{\psi}K\xrightarrow{\xi}L}\prod_{l\in L}\nu_{\phi^{-1}\psi^{-1}\xi^{-1}l\xrightarrow{\phi|}\psi^{-1}\xi^{-1}l\xrightarrow{\psi|}\xi^{-1}l}$$

In detail, after canonical identification of the source of $(1\times m)\centerdot m$ with the following this functor

becomes the composition

$$\coprod_{I\xrightarrow{\phi}J\xrightarrow{\psi}K\xrightarrow{\xi}L}\prod_{l\in L}\Big[\Big(\prod_{j\in\psi^{-1}\xi^{-1}l}\mathsf{C}\big((X_i)_{i\in\phi^{-1}j};Y_j\big)\Big)\times\Big(\prod_{k\in\xi^{-1}l}\mathsf{C}\big((Y_j)_{j\in\psi^{-1}k};Z_k\big)\Big)$$

$$\times\,\mathsf{C}\big((Z_k)_{k\in\xi^{-1}l};W_l\big)\Big]$$

$$\xrightarrow{\coprod_{(\phi,\psi,\xi)\mapsto(\phi,\psi\centerdot\xi)}\prod_{l\in L}(1\times\mu_{\psi|:\,\psi^{-1}\xi^{-1}l\to\xi^{-1}l})}$$

$$\coprod_{I\xrightarrow{\phi}J\xrightarrow{\zeta}L}\prod_{l\in L}\Big[\Big(\prod_{j\in\zeta^{-1}l}\mathsf{C}\big((X_i)_{i\in\phi^{-1}j};Y_j\big)\Big)\times\mathsf{C}\big((Y_j)_{j\in\zeta^{-1}l};W_l\big)\Big]$$

$$\xrightarrow{\coprod_{(\phi,\zeta)\mapsto(\phi\centerdot\zeta)}\prod_{l\in L}\mu_{\phi|:\,\phi^{-1}\zeta^{-1}l\to\zeta^{-1}l}}\coprod_{\chi:\,I\to L}\prod_{l\in L}\mathsf{C}\big((X_i)_{i\in\chi^{-1}l};W_l\big).$$

The source of $(m\times 1)\centerdot m$ is identified canonically with another expression, after which this functor becomes the composition

$$\coprod_{I\xrightarrow{\phi}J\xrightarrow{\psi}K\xrightarrow{\xi}L}\prod_{l\in L}\Big\{\prod_{k\in\xi^{-1}l}\Big[\prod_{j\in\psi^{-1}k}\mathsf{C}\big((X_i)_{i\in\phi^{-1}j};Y_j\big)\times\mathsf{C}\big((Y_j)_{j\in\psi^{-1}k};Z_k\big)\Big]$$

$$\times\,\mathsf{C}\big((Z_k)_{k\in\xi^{-1}l};W_l\big)\Big\}$$

$$\xrightarrow{\coprod_{(\phi,\psi,\xi)\mapsto(\phi\centerdot\psi,\xi)}\prod_{l\in L}\{\prod_{k\in\xi^{-1}l}[\mu_{\phi|:\,\phi^{-1}\psi^{-1}k\to\psi^{-1}k}]\times 1\}}$$

$$\coprod_{I\xrightarrow{\sigma}K\xrightarrow{\xi}L}\prod_{l\in L}\Big\{\Big[\prod_{k\in\xi^{-1}l}\mathsf{C}\big((X_i)_{i\in\sigma^{-1}k};Z_k\big)\Big]\times\mathsf{C}\big((Z_k)_{k\in\xi^{-1}l};W_l\big)\Big\}$$

$$\xrightarrow{\coprod_{(\sigma,\xi)\mapsto(\sigma\centerdot\xi)}\prod_{l\in L}\mu_{\sigma|:\,\sigma^{-1}\xi^{-1}l\to\xi^{-1}l}}\coprod_{\chi:\,I\to L}\prod_{l\in L}\mathsf{C}\big((X_i)_{i\in\chi^{-1}l};W_l\big).$$

The pentagon equation for $\mathsf{a}$ holds true being on a summand indexed by $I\xrightarrow{\phi}J\xrightarrow{\psi}K\xrightarrow{\xi}L\xrightarrow{\kappa}M$ of the product over $m\in M$ of pentagons [Lyu25, (1.1.2)] written for maps

$$\phi^{-1}\psi^{-1}\xi^{-1}\kappa^{-1}m\xrightarrow{\phi|}\psi^{-1}\xi^{-1}\kappa^{-1}m\xrightarrow{\psi|}\xi^{-1}\kappa^{-1}m\xrightarrow{\xi|}\kappa^{-1}m.$$

Let us prove equation (1.1.2) for $\mathcal{P}=\mathbb{F}\mathsf{C}$. Take arbitrary objects $A_i=(X_i^s)_{s\in S_i}$, $i\in I$, $B_j=(Y_j^q)_{q\in Q_j}$, $j\in J$, $C_k=(Z_k^p)_{p\in P_k}$, $k\in K$, $D_l=(W_l^r)_{r\in R_l}$, $l\in L$. Denote $S=\bigsqcup_{i\in I}S_i$, $Q=\bigsqcup_{j\in J}Q_j$, $P=\bigsqcup_{k\in K}P_k$, $R=\bigsqcup_{l\in L}R_l$. Let $f$, $g$, $h$, $e$ be maps from $\mathcal{S}_{\mathsf{sk}}$ as on the following diagram

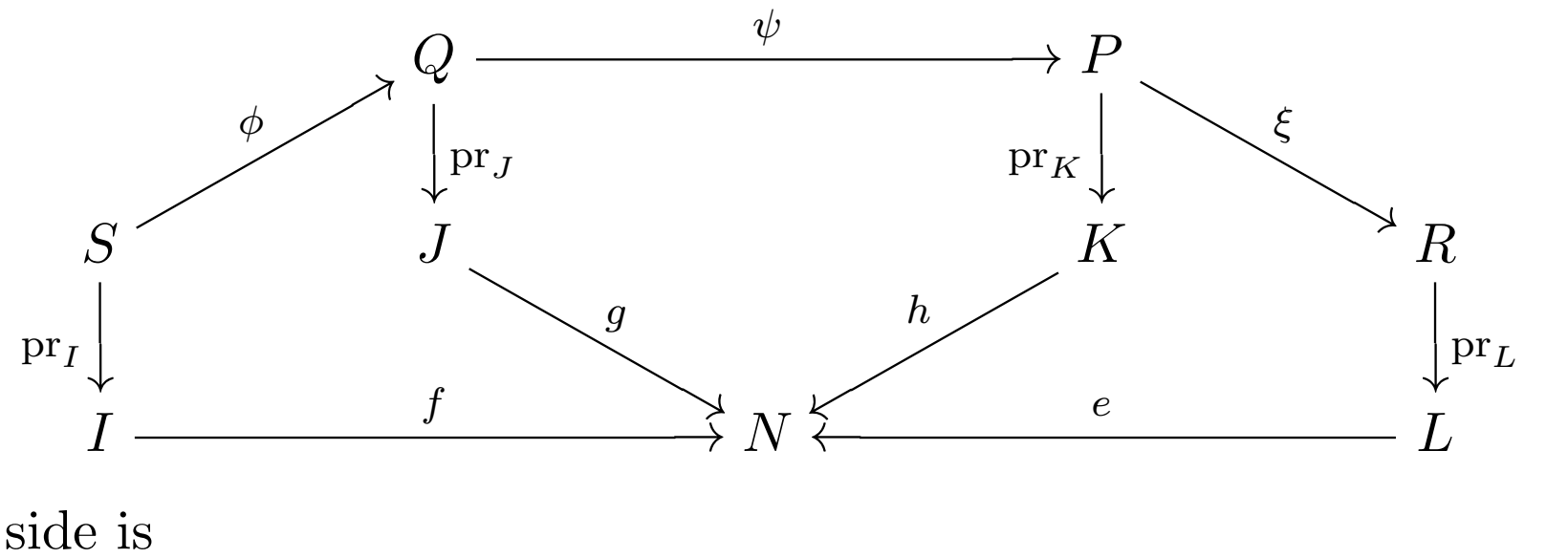

(1.2.5)

The left hand side is

$$\bigotimes^{I\xrightarrow{f}N\xleftarrow{e}L}\prod_{n\in N}\coprod_{\sqcup_{i\in f^{-1}n}S_i\xrightarrow{\phi_n}\sqcup_{j\in g^{-1}n}Q_j\xrightarrow{\psi_n}\sqcup_{j\in h^{-1}n}P_k\xrightarrow{\xi_n}\sqcup_{l\in e^{-1}n}R_l}\prod_{r\in\sqcup_{l\in e^{-1}n}R_l}\nu_{\phi_n^{-1}\psi_n^{-1}\xi_n^{-1}r\xrightarrow{\phi_n|}\psi_n^{-1}\xi_n^{-1}r\xrightarrow{\psi_n|}\xi_n^{-1}r}. \tag{1.2.6}$$

The collection $(\phi_n)_{n\in N}$ (resp. $(\psi_n)_{n\in N}$, $(\xi_n)_{n\in N}$) is encoded by a single map $\phi\colon S\to Q$ (resp. $\psi\colon Q\to P$, $\xi\colon P\to R$) which makes the corresponding pentagon of (1.2.5) commutative and $\phi|_{\sqcup_{i\in f^{-1}n}S_i} = \phi_n$ (resp. $\psi|_{\sqcup_{j\in g^{-1}n}Q_j} = \psi_n$, $\xi|_{\sqcup_{k\in h^{-1}n}P_k} = \xi_n$). The right hand side of (1.1.2) is

$$\coprod_{S\xrightarrow{\phi}Q\xrightarrow{\psi}P\xrightarrow{\xi}R} \prod_{r\in R} \nu_{\phi^{-1}\psi^{-1}\xi^{-1}r \xrightarrow{\phi|} \psi^{-1}\xi^{-1}r \xrightarrow{\psi|} \xi^{-1}r}$$

which coincides with (1.2.6) due to the following simple observation. Consider an element $r \in R$. Let $n \in N$ be such that $r \in \sqcup_{l\in e^{-1}n} R_l$. Then

$$(\phi^{-1}\psi^{-1}\xi^{-1}r \xrightarrow{\phi|} \psi^{-1}\xi^{-1}r \xrightarrow{\psi|} \xi^{-1}r) = (\phi_n^{-1}\psi_n^{-1}\xi_n^{-1}r \xrightarrow{\phi_n|} \psi_n^{-1}\xi_n^{-1}r \xrightarrow{\psi_n|} \xi_n^{-1}r).$$

The left unitor

$$\begin{array}{ccc}
\mathbb{1} \times \mathbb{F}\mathsf{C}\big((X_i)_{i\in I}, (Y_j)_{j\in J}\big) & \xrightarrow{\;\mathrm{i}\times 1\;} & \\
l \downarrow \cong \quad \overset{\mathsf{l}}{\Longleftarrow} & & \mathbb{F}\mathsf{C}\big((X_i)_{i\in I}, (X_i)_{i\in I}\big) \times \mathbb{F}\mathsf{C}\big((X_i)_{i\in I}, (Y_j)_{j\in J}\big) \\
\mathbb{F}\mathsf{C}\big((X_i)_{i\in I}, (Y_j)_{j\in J}\big) & \xleftarrow{\;m\;} &
\end{array}$$

is the isomorphism $\mathsf{l} = \coprod_{\phi\colon I\to J} \prod_{j\in J} \theta_{\phi^{-1}j}$ presented in diagram

$$\begin{array}{ccc}
\coprod_{\phi\colon I\to J} \prod_{j\in J} \big[\mathbb{1} \times \mathsf{C}\big((X_i)_{i\in\phi^{-1}j}; Y_j\big)\big] & \xrightarrow{\;\coprod_\phi \prod_J [(\mathrm{i})\times 1]\;} & \\
\coprod_\phi \prod_J l \downarrow \cong \quad \overset{\mathsf{l}}{\Longleftarrow} & & \coprod_{\phi\colon I\to J} \prod_{j\in J} \Big[\big(\prod_{i\in\phi^{-1}j} \mathsf{C}(X_i; X_i)\big) \times \mathsf{C}\big((X_i)_{i\in\phi^{-1}j}; Y_j\big)\Big] \\
\coprod_{\phi\colon I\to J} \prod_{j\in J} \big[\mathsf{C}\big((X_i)_{i\in\phi^{-1}j}; Y_j\big)\big] & \xleftarrow{\;\coprod_{\phi\colon I\to J} \prod_{j\in J} \mu_{\mathrm{id}_{\phi^{-1}j}}\;} &
\end{array} \tag{1.2.7}$$

The right unitor

$$\begin{array}{ccc}
\mathbb{F}\mathsf{C}\big((X_i)_{i\in I}, (Y_j)_{j\in J}\big) \times \mathbb{1} & \xrightarrow{\;1\times\mathrm{i}\;} & \\
r \downarrow \cong \quad \overset{\mathsf{r}}{\Longleftarrow} & & \mathbb{F}\mathsf{C}\big((X_i)_{i\in I}, (Y_j)_{j\in J}\big) \times \mathbb{F}\mathsf{C}\big((Y_j)_{j\in J}, (Y_j)_{j\in J}\big) \\
\mathbb{F}\mathsf{C}\big((X_i)_{i\in I}, (Y_j)_{j\in J}\big) & \xleftarrow{\;m\;} &
\end{array}$$

is the isomorphism $\mathsf{r} = \coprod_{\phi\colon I\to J} \prod_{j\in J} \zeta_{\phi^{-1}j}$ presented in diagram

$$\begin{array}{ccc}
\coprod_{\phi\colon I\to J} \prod_{j\in J} \big[\mathsf{C}\big((X_i)_{i\in\phi^{-1}j}; Y_j\big) \times \mathbb{1}\big] & \xrightarrow{\;\coprod_\phi \prod_J [1\times\mathrm{i}]\;} & \\
\coprod_\phi \prod_J r \downarrow \cong \quad \overset{\mathsf{r}}{\Longleftarrow} & & \coprod_{\phi\colon I\to J} \prod_{j\in J} \big[\mathsf{C}\big((X_i)_{i\in\phi^{-1}j}; Y_j\big) \times \mathsf{C}(Y_j; Y_j)\big] \\
\coprod_{\phi\colon I\to J} \prod_{j\in J} \big[\mathsf{C}\big((X_i)_{i\in\phi^{-1}j}; Y_j\big)\big] & \xleftarrow{\;\coprod_{\phi\colon I\to J} \prod_{j\in J} \mu_{\triangledown\colon \phi^{-1}j\to\mathbf{1}}\;} &
\end{array} \tag{1.2.8}$$

The relation between

$$\mathsf{a} = \coprod_{I\xrightarrow{\phi}J\xrightarrow{\psi}K} \prod_{k\in K} \nu_{\phi^{-1}\psi^{-1}k \xrightarrow{\phi|} \psi^{-1}k \xrightarrow{\mathrm{id}} \psi^{-1}k}, \quad \mathsf{r} = \coprod_{\phi\colon I\to J} \prod_{k\in K} \prod_{j\in\psi^{-1}k} \zeta_{\phi^{-1}j},$$

$$\text{and } \mathsf{l} = \coprod_{\psi\colon J\to K} \prod_{k\in K} \theta_{\psi^{-1}k}$$

on the summand indexed by $I \xrightarrow{\phi} J \xrightarrow{\psi} K$ is obtained as the product over $k \in K$ of relations [Lyu25, (1.1.4)] between $\nu_{\phi^{-1}\psi^{-1}k \xrightarrow{\phi|} \psi^{-1}k \xrightarrow{\mathrm{id}} \psi^{-1}k}$, $\prod_{j\in\psi^{-1}k} \zeta_{\phi^{-1}j}$ and $\theta_{\psi^{-1}k}$, see diagram at Figure 1 on the following page.

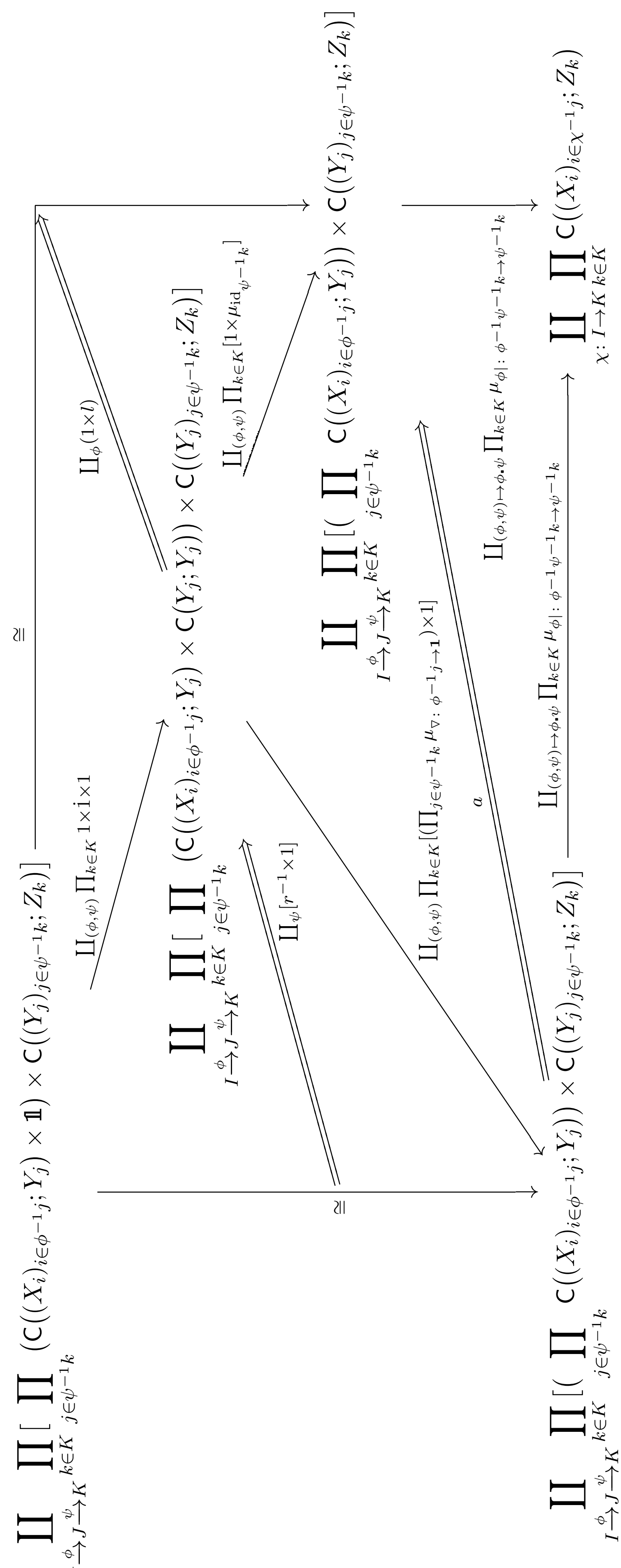

Figure 1: A triangle relation between $\mathsf{a}$, $\mathsf{r}$ and $\mathsf{l}$

Property (1.1.3) is checked below. Let $A_i = (X_i^s)_{s\in S_i}$, $i \in I$. Denote $S = \bigsqcup_{i\in I} S_i$. Then

$$\prod_{l\in L} \mathbb{F}\mathsf{C}\big((X^s_{\mathrm{pr}_I s})_{s\in\sqcup_{i\in f^{-1}l}S_i};(X^s_{\mathrm{pr}_I s})_{s\in\sqcup_{i\in f^{-1}l}S_i}\big) \hookleftarrow \prod_{l\in L}\prod_{i\in f^{-1}l}\prod_{s\in S_i}\mathsf{C}(X_i^s;X_i^s)$$

$$\cong \prod_{i\in I}\prod_{s\in S_i}\mathsf{C}(X_i^s;X_i^s) \hookrightarrow \mathbb{F}\mathsf{C}\big((X^s_{\mathrm{pr}_I s})_{s\in S};(X^s_{\mathrm{pr}_I s})_{s\in S}\big)$$

$$\big(1_{(X^s_{\mathrm{pr}_I s})_{s\in\sqcup_{i\in f^{-1}l}S_i}}\big)_{l\in L} \longleftarrow\!\shortmid \big(((1_{X_i^s})_{s\in S_i})_{i\in f^{-1}l}\big)_{l\in L} = ((1_{X_i^s})_{s\in S_i})_{i\in I} \mapsto 1_{(X^s_{\mathrm{pr}_I s})_{s\in S}}.$$

Equation (1.1.4) takes the following form. Look at diagram (1.2.1). We denote $\phi_k = \phi|\colon \sqcup_{i\in f^{-1}k} S_i \to \sqcup_{j\in g^{-1}k} Q_j$. Denote $Q^k = \bigsqcup_{j\in g^{-1}k} Q_j$. Taking all that into account, (1.1.4) becomes

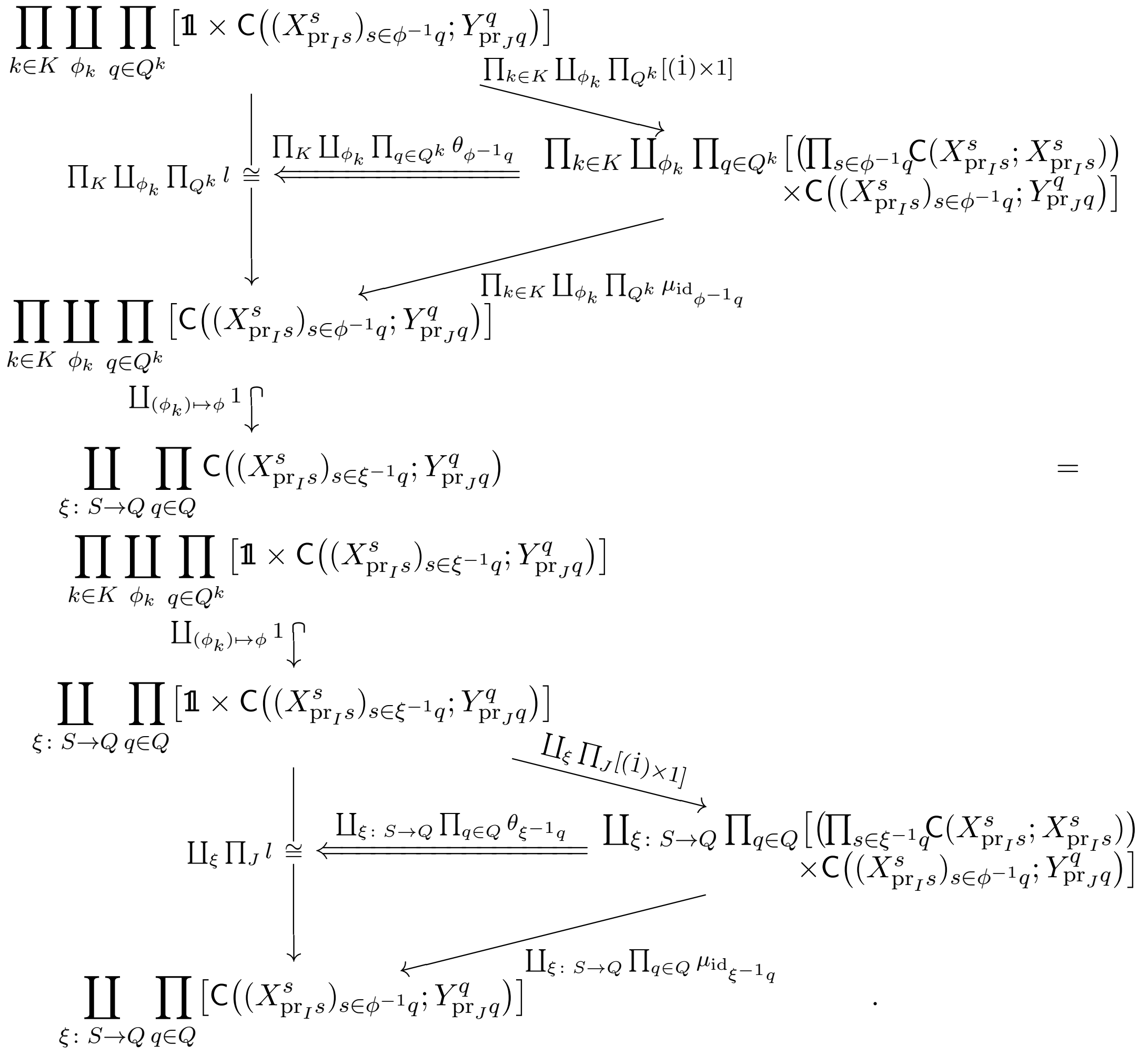

Similarly for right unitors equation (1.1.5) takes the form

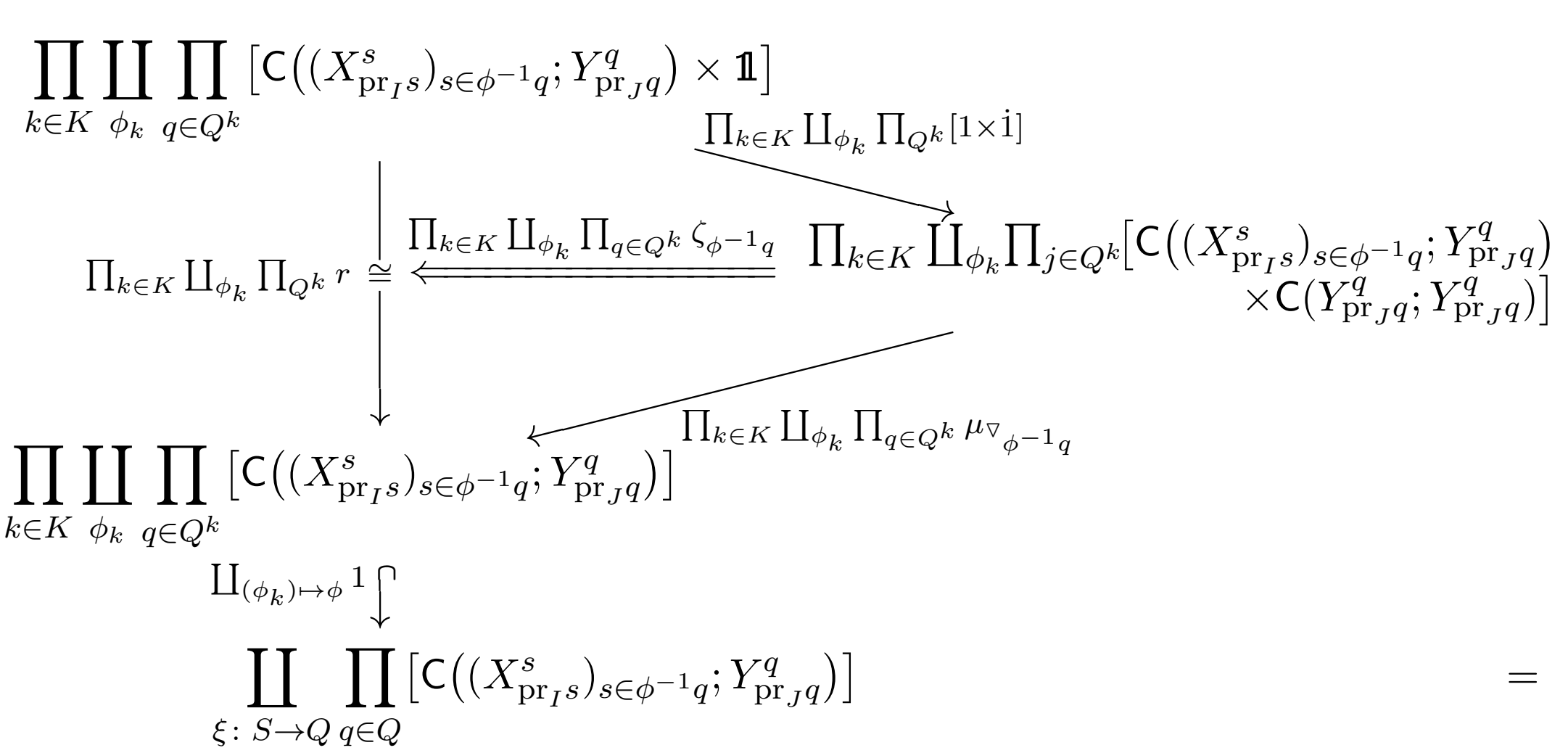

$$
\begin{array}{ccc}
\prod_{k\in K}\coprod_{\phi_k}\prod_{q\in Q^k}\bigl[\mathsf{C}\bigl((X^s_{\mathrm{pr}_I s})_{s\in\phi^{-1}q};Y^q_{\mathrm{pr}_J q}\bigr)\times\mathbb{1}\bigr] & & \\
\Big\downarrow{\scriptstyle \coprod_{(\phi_k)\mapsto\phi}1} & & \\
\coprod_{\xi\colon S\to Q}\prod_{q\in Q}\bigl[\mathsf{C}\bigl((X^s_{\mathrm{pr}_I s})_{s\in\xi^{-1}q};Y^q_{\mathrm{pr}_J q}\bigr)\times\mathbb{1}\bigr] & \xrightarrow{\coprod_\xi\prod_J[1\times i]} & \\
\Big\downarrow{\scriptstyle \coprod_\xi\prod_J r\ \cong} & \xLeftarrow{\coprod_{\xi\colon I\to J}\prod_{j\in J}\zeta_{\xi^{-1}j}} & \coprod_{\xi\colon S\to Q}\prod_{q\in Q}\bigl[\mathsf{C}\bigl((X^s_{\mathrm{pr}_I s})_{s\in\xi^{-1}q};Y^q_{\mathrm{pr}_J q}\bigr)\times\mathsf{C}(Y^q_{\mathrm{pr}_J q};Y^q_{\mathrm{pr}_J q})\bigr] \\
& \xleftarrow{\coprod_{\xi\colon S\to Q}\prod_{q\in Q}\mu_{\triangledown\colon\xi^{-1}q\to\mathbf{1}}} & \\
\coprod_{\xi\colon S\to Q}\prod_{q\in Q}\bigl[\mathsf{C}\bigl((X^s_{\mathrm{pr}_I s})_{s\in\xi^{-1}q};Y^q_{\mathrm{pr}_J q}\bigr)\bigr] & & .
\end{array}
$$

These equations hold by inspection.

Now let us prove (1.1.6). Assume that $A_i=(X^s_i)_{s\in S_i}$, $B_j=(Y^p_j)_{p\in P_j}$, $X^s_i,Y^p_j\in\operatorname{Ob}\mathsf{C}$, $S_i,P_j\in\operatorname{Ob}\mathcal{O}_{\mathsf{sk}}$. Then

$$
\begin{aligned}
\prod_{k\in K}\mathbb{F}\mathsf{C}\Bigl(\bigsqcup_{i\in f^{-1}k}A_i;\bigsqcup_{j\in g^{-1}k}B_j\Bigr)
&=\prod_{k\in K}\coprod_{\phi_k\colon\sqcup_{i\in f^{-1}k}S_i\to\sqcup_{j\in g^{-1}k}P_j}\prod_{p\in\sqcup_{j\in g^{-1}k}P_j}\mathsf{C}\bigl((X^s_i)_{(i,s)\in\phi_k^{-1}(j,p)};Y^p_j\bigr)\\
&\cong\coprod_{(\phi_k\colon\sqcup_{i\in f^{-1}k}S_i\to\sqcup_{j\in g^{-1}k}P_j)_{k\in K}}\prod_{k\in K}\prod_{j\in.g^{-1}k}\prod_{p\in P_j}\mathsf{C}\bigl((X^s_i)_{(i,s)\in\phi_k^{-1}(j,p)};Y^p_j\bigr).
\end{aligned}
$$

Denote $S=\coprod_{i\in I}S_i$, $P=\coprod_{j\in J}P_j$. The collection $(\phi_k)_{k\in K}$ amounts to a single mapping $\phi\colon S\to P$ such that

$$
\begin{array}{ccc}
S & \xrightarrow{\quad\phi\quad} & P\\
{\scriptstyle \mathrm{pr}_I\centerdot f=:\tilde f}\searrow & = & \swarrow{\scriptstyle \tilde g:=\mathrm{pr}_J\centerdot g}\\
& K &
\end{array}
$$

The bijective correspondence is established by $\phi_k=\phi|_{\sqcup_{i\in f^{-1}k}S_i}$. Thus,

$$
\prod_{k\in K}\mathbb{F}\mathsf{C}\Bigl(\bigsqcup_{i\in f^{-1}k}A_i;\bigsqcup_{j\in g^{-1}k}B_j\Bigr)=\coprod_{\phi\colon S\to P}^{\phi\centerdot\tilde g=\tilde f}\prod_{p\in P}\mathsf{C}\bigl((X^s_{\mathrm{pr}_I s})_{s\in\phi^{-1}p};Y^p_{\mathrm{pr}_J p}\bigr).
$$

The upper row of diagram (1.1.6) is

$$
\begin{aligned}
\otimes^{I\xrightarrow{f}K\xleftarrow{g}J}=\coprod_{\phi\mapsto\phi}1\colon\ &\coprod_{\phi\colon S\to P}^{\phi\centerdot\tilde g=\tilde f}\prod_{p\in P}\mathsf{C}\bigl((X^s_{\mathrm{pr}_I s})_{s\in\phi^{-1}p};Y^p_{\mathrm{pr}_J p}\bigr)\\
&\to\coprod_{\xi\colon S\to P}\prod_{p\in P}\mathsf{C}\bigl((X^s_{\mathrm{pr}_I s})_{s\in\xi^{-1}p};Y^p_{\mathrm{pr}_J p}\bigr)=\mathbb{F}\mathsf{C}\Bigl(\bigsqcup_{i\in I}A_i;\bigsqcup_{j\in J}B_j\Bigr).
\end{aligned}
$$

Denote $S^l=\coprod_{i\in\alpha^{-1}l}S_i$, $P^l=\coprod_{j\in\beta^{-1}l}P_j$. The composition of three other arrows of (1.1.6) is

$$
\begin{aligned}
\prod_{l\in L}\prod_{k\in h^{-1}l}\mathbb{F}\mathsf{C}\Bigl(\bigsqcup_{i\in f^{-1}k}A_i;\bigsqcup_{j\in g^{-1}k}B_j\Bigr)&=\prod_{l\in L}\coprod_{\phi^l\colon S^l\to P^l}^{\phi\centerdot\tilde\beta|=\tilde\alpha|}\prod_{p\in P^l}\mathsf{C}\bigl((X^s_{\mathrm{pr}_I s})_{s\in(\phi^l)^{-1}p};Y^p_{\mathrm{pr}_J p}\bigr)\\
&\xhookrightarrow{\prod_{l\in L}\coprod_{\phi^l\mapsto\phi^l}1}\prod_{l\in L}\coprod_{\phi^l\colon S^l\to P^l}\prod_{p\in P^l}\mathsf{C}\bigl((X^s_{\mathrm{pr}_I s})_{s\in(\phi^l)^{-1}p};Y^p_{\mathrm{pr}_J p}\bigr)=\\
\coprod_{\psi\colon S\to P}^{\psi\centerdot\tilde\beta=\tilde\alpha}\prod_{p\in P}\mathsf{C}\bigl((X^s_{\mathrm{pr}_I s})_{s\in\psi^{-1}p};Y^p_{\mathrm{pr}_J p}\bigr)&\xhookrightarrow{\coprod_{\psi\mapsto\psi}1}\coprod_{\xi\colon S\to P}\prod_{p\in P}\mathsf{C}\bigl((X^s_{\mathrm{pr}_I s})_{s\in\xi^{-1}p};Y^p_{\mathrm{pr}_J p}\bigr)
\end{aligned}
$$

$$= \mathbb{F}\mathsf{C}\big(\bigsqcup_{i\in I} A_i; \bigsqcup_{j\in J} B_j\big).$$

Here $\tilde{\alpha} = \mathrm{pr}_I \centerdot \alpha$, $\tilde{\beta} = \mathrm{pr}_J \centerdot \beta$ and $\psi\colon S \to P$ has to make the triangle

$$\begin{array}{ccc} S & \xrightarrow{\psi} & P \\ & = & \\ {\scriptstyle \mathrm{pr}_I \centerdot f \centerdot h =: \tilde{\alpha}} \searrow & L & \swarrow {\scriptstyle \tilde{\beta} := \mathrm{pr}_J \centerdot g \centerdot h} \end{array}$$

commutative. The restriction $\phi^l = \psi|_{\tilde{\alpha}^{-1}l}$ has to make the triangle

$$\begin{array}{ccc} S^l & \xrightarrow{\phi^l} & P^l \\ & = & \\ {\scriptstyle \tilde{\alpha}|} \searrow & h^{-1}l & \swarrow {\scriptstyle \tilde{\beta}|} \end{array}$$

commutative. Therefore, diagram of functors (1.1.6) commutes (up to symmetry in $\mathcal{S}et$).

## 1.3 Definition of morphisms of biprops

**1.3.1 Definition.** Let $\mathcal{P}$, $\mathcal{R}$ be biprops. A weak morphism of biprops $F\colon \mathcal{P} \to \mathcal{R}$ is a particular case of a homomorphism of symmetric monoidal bicategories (in unbiased form). It consists of the following data:

— a map of colours $\operatorname{Col} F\colon \operatorname{Col}\mathcal{P} \to \operatorname{Col}\mathcal{R}$ (which induces the map of objects

$$\operatorname{Ob} F = (\operatorname{Col} F)^*\colon (\operatorname{Col}\mathcal{P})^* \to (\operatorname{Col}\mathcal{R})^*, \quad A = (X_i)_{i\in I} \mapsto FA = (FX_i)_{i\in I}),$$

— a functor $F_{A;B}\colon \mathcal{P}(A;B) \to \mathcal{R}(FA;FB)$ for each pair of objects $A, B \in \operatorname{Ob}\mathcal{P}$,

— a natural isomorphism

$$\begin{array}{ccc} \mathcal{P}(A;B) \times \mathcal{P}(B;C) & \xrightarrow{m} & \mathcal{P}(A;C) \\ {\scriptstyle F_{A;B}\times F_{B;C}} \downarrow & \overset{F_{A;B;C}}{\nearrow\!\!\!\nearrow} & \downarrow {\scriptstyle F_{A;C}} \\ \mathcal{R}(FA;FB) \times \mathcal{R}(FB;FC) & \xrightarrow{m} & \mathcal{R}(FA;FC) \end{array}$$

— a natural isomorphism (a 2-isomorphism of $\mathcal{R}$)

$$\begin{array}{ccc} & \mathcal{P}(A;A) & \\ {\scriptstyle \mathrm{i}_A} \nearrow & \Downarrow {\scriptstyle F_A} & \searrow {\scriptstyle F_{A;A}} \\ \mathbb{1} & \xrightarrow[\mathrm{i}_{FA}]{} & \mathcal{R}(FA;FA) \end{array}$$

such that

$$\begin{array}{ccc} \prod_{l\in L} \mathcal{P}\big(\bigsqcup_{i\in f^{-1}l} A_i; \bigsqcup_{j\in g^{-1}l} B_j\big) & \xrightarrow{\otimes^{I \xrightarrow{f} L \xleftarrow{g} J}} & \mathcal{P}\big(\bigsqcup_{i\in I} A_i; \bigsqcup_{j\in J} B_j\big) \\ {\scriptstyle \prod_{l\in L} F_{\sqcup_{i\in f^{-1}l} A_i; \sqcup_{j\in g^{-1}l} B_j}} \downarrow & = & \downarrow {\scriptstyle F_{\sqcup_{i\in I} A_i; \sqcup_{j\in J} B_j}} \\ \prod_{l\in L} \mathcal{R}\big(\bigsqcup_{i\in f^{-1}l} FA_i; \bigsqcup_{j\in g^{-1}l} FB_j\big) & \xrightarrow{\otimes^{I \xrightarrow{f} L \xleftarrow{g} J}} & \mathcal{R}\big(\bigsqcup_{i\in I} FA_i; \bigsqcup_{j\in J} FB_j\big) \end{array} \tag{1.3.1}$$

is a commutative diagram of functors; the equation holds

$$\prod_{l\in L}\mathcal{P}\Big(\bigsqcup_{i\in f^{-1}l} A_i;\ \bigsqcup_{j\in g^{-1}l} B_j\Big) \qquad \prod_{l\in L}\mathcal{P}\Big(\bigsqcup_{j\in g^{-1}l} B_j;\ \bigsqcup_{k\in h^{-1}l} C_k\Big) \tag{1.3.2}$$

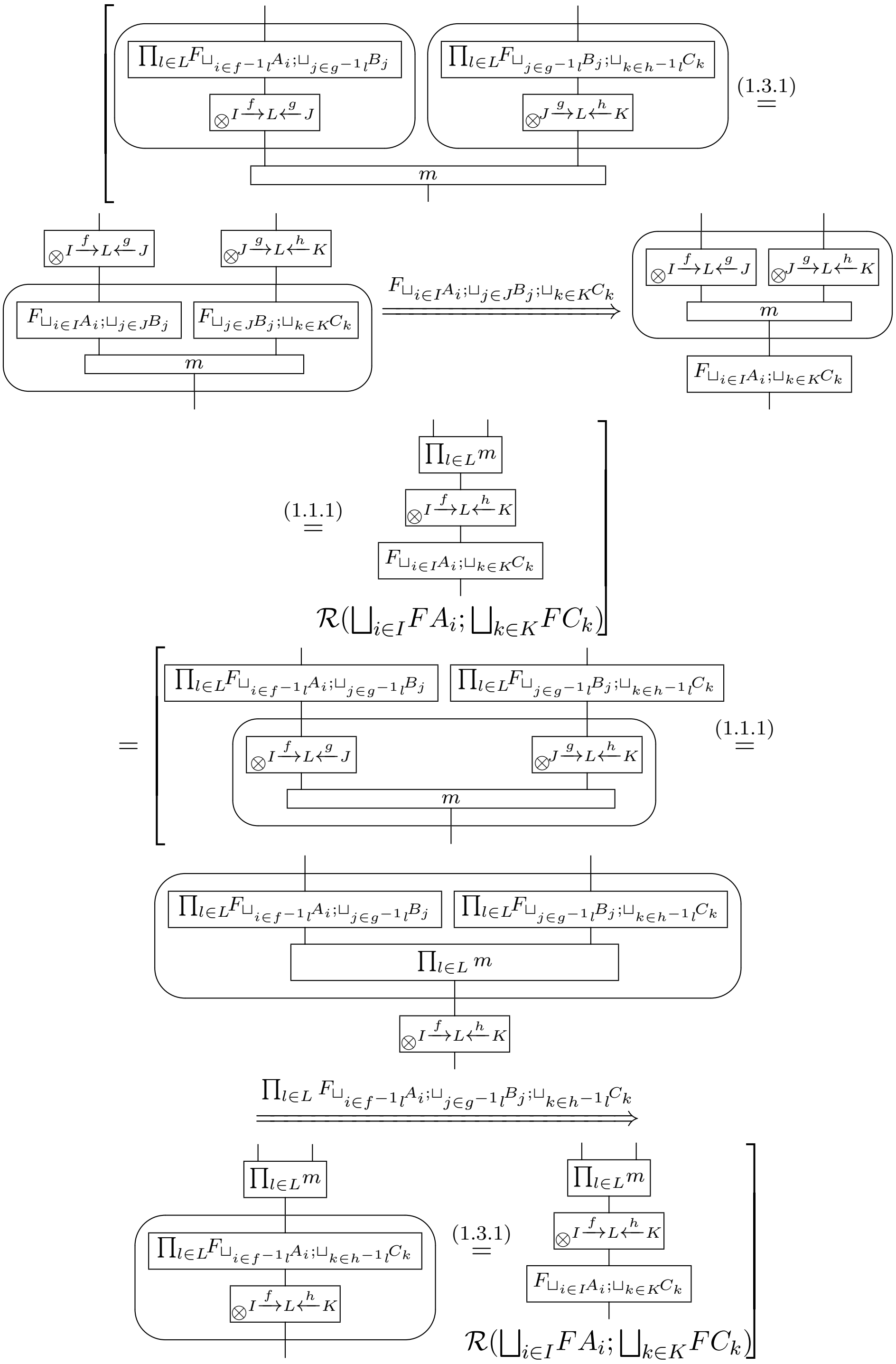

Another (pasting) form of this equation is given at Figure 2. Compatibility with compositions equation at Figure 3 must hold. The same equation but in Bartlett's notation:

$$\mathcal{P}\big((X_i)_{i\in I};(Y_j)_{j\in J}\big)\times\mathcal{P}\big((Y_j)_{j\in J};(Z_k)_{k\in K}\big)\times\mathcal{P}\big((Z_k)_{k\in K};(W_l)_{l\in L}\big) \tag{1.3.3}$$

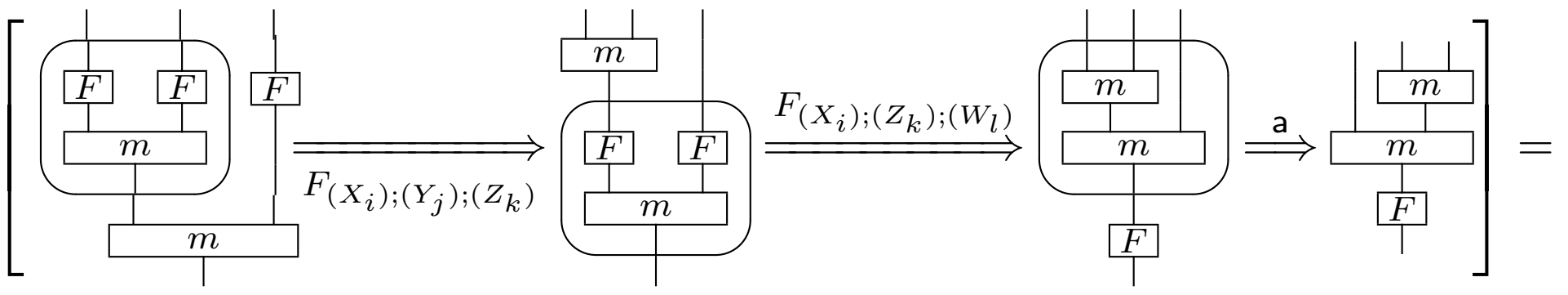

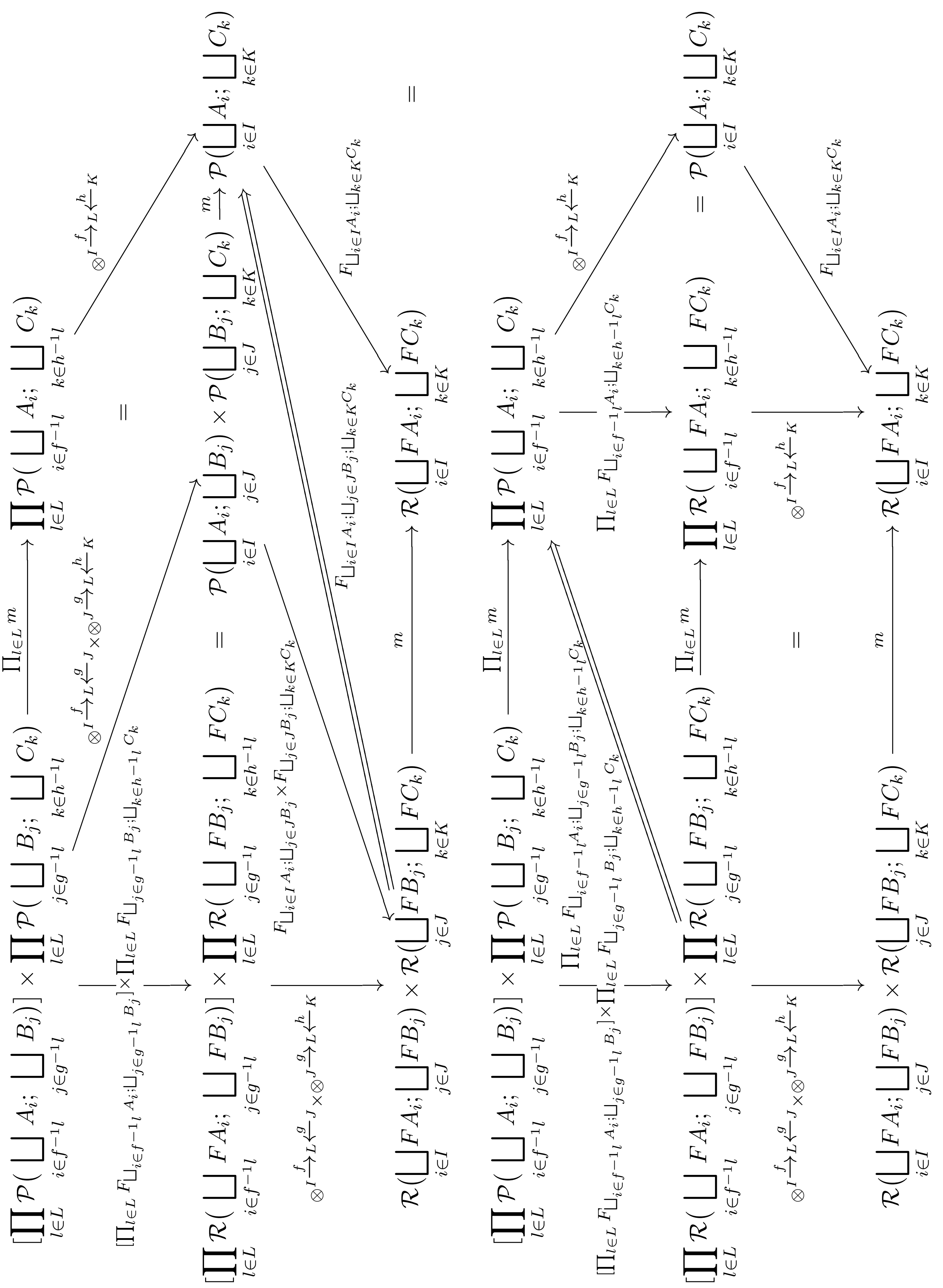

Figure 2: Coherence of a morphism of biprops with the tensor multiplication

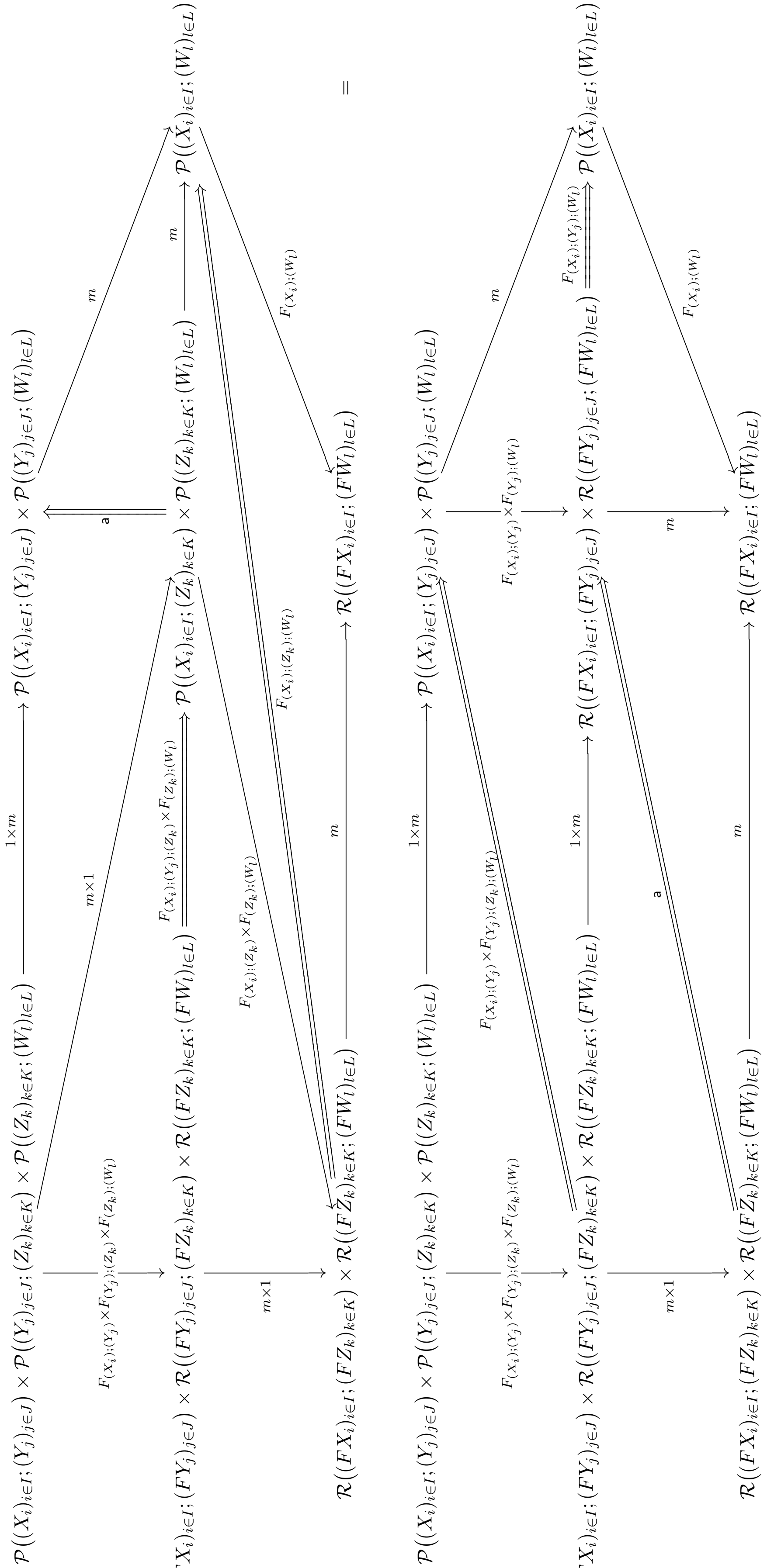

Figure 3: Coherence of a morphism of biprops with the composition

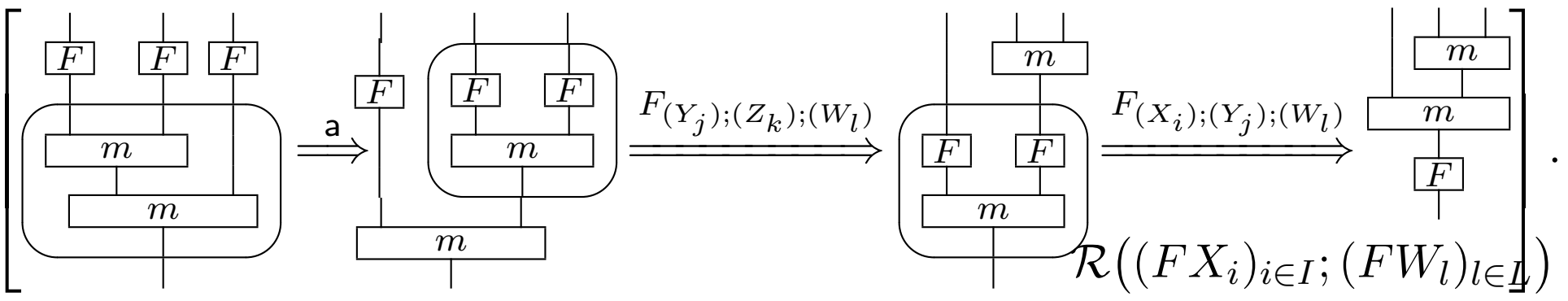

Finally, equations for unitors have to hold:

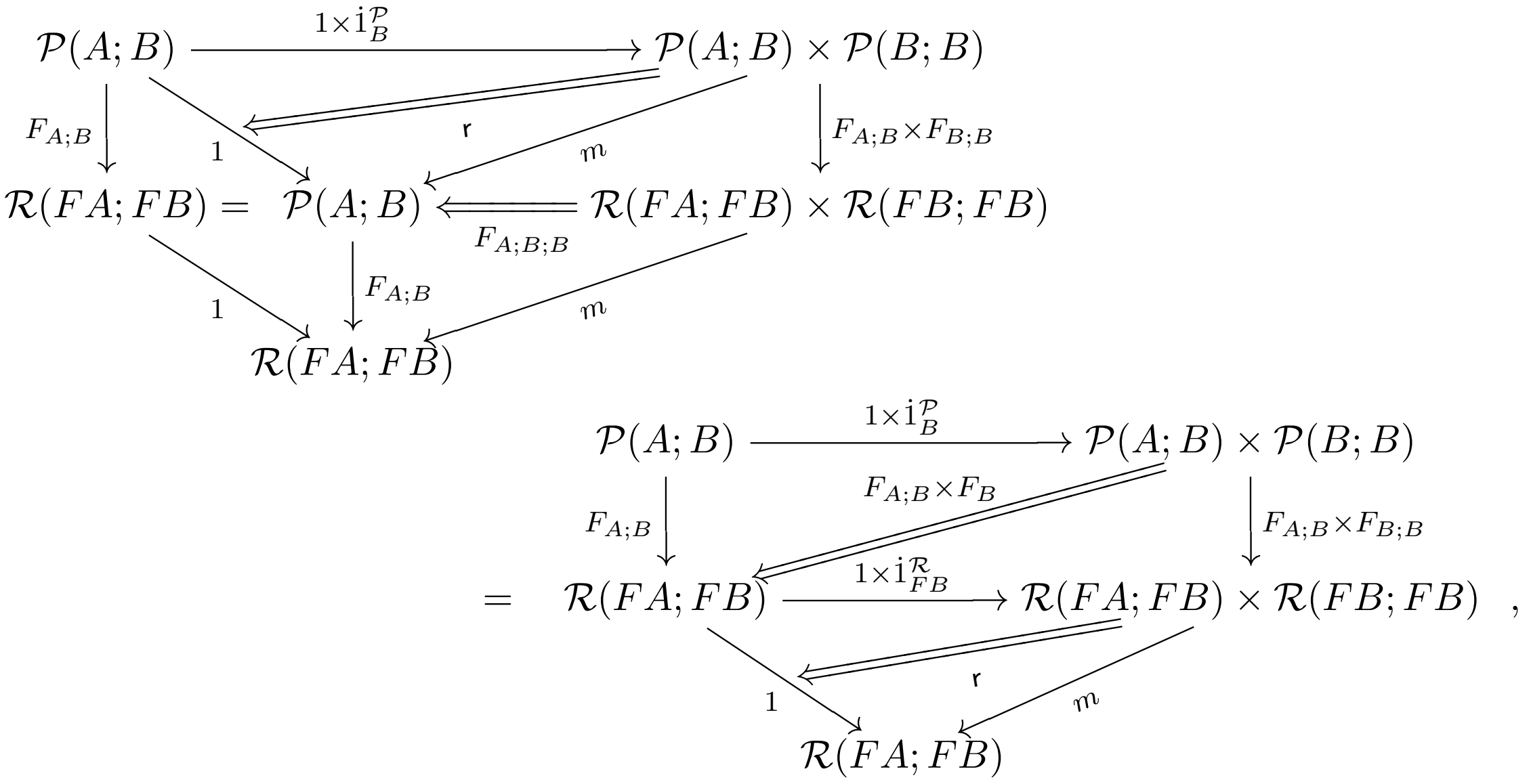

(1.3.4)

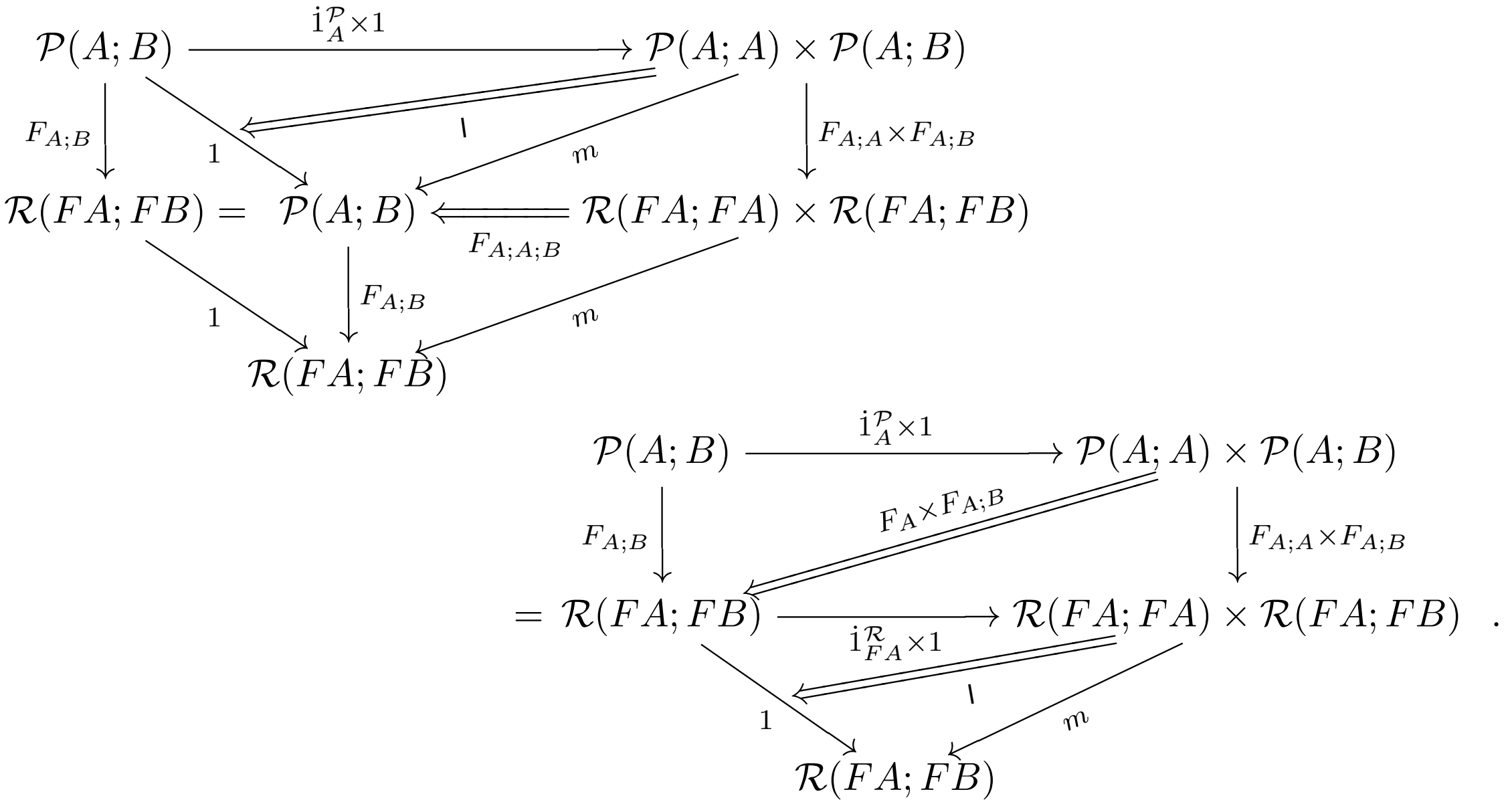

(1.3.5)

Having two morphisms of biprops $F\colon \mathcal{P}\to\mathcal{R}$, $G\colon\mathcal{R}\to\mathcal{Q}$ we define their composition $K\colon\mathcal{P}\to\mathcal{Q}$ with the components

— $\operatorname{Col} K = \big(\operatorname{Col}\mathcal{P}\xrightarrow{\operatorname{Col} F}\operatorname{Col}\mathcal{R}\xrightarrow{\operatorname{Col} G}\operatorname{Col}\mathcal{Q}\big)$;

— $K_{A;B} = \big[\mathcal{P}(A;B)\xrightarrow{F_{A;B}}\mathcal{R}(FA;FB)\xrightarrow{G_{FA;FB}}\mathcal{Q}(GFA;GFB)\big]$;

— $K_{A;B;C} = \chi =$

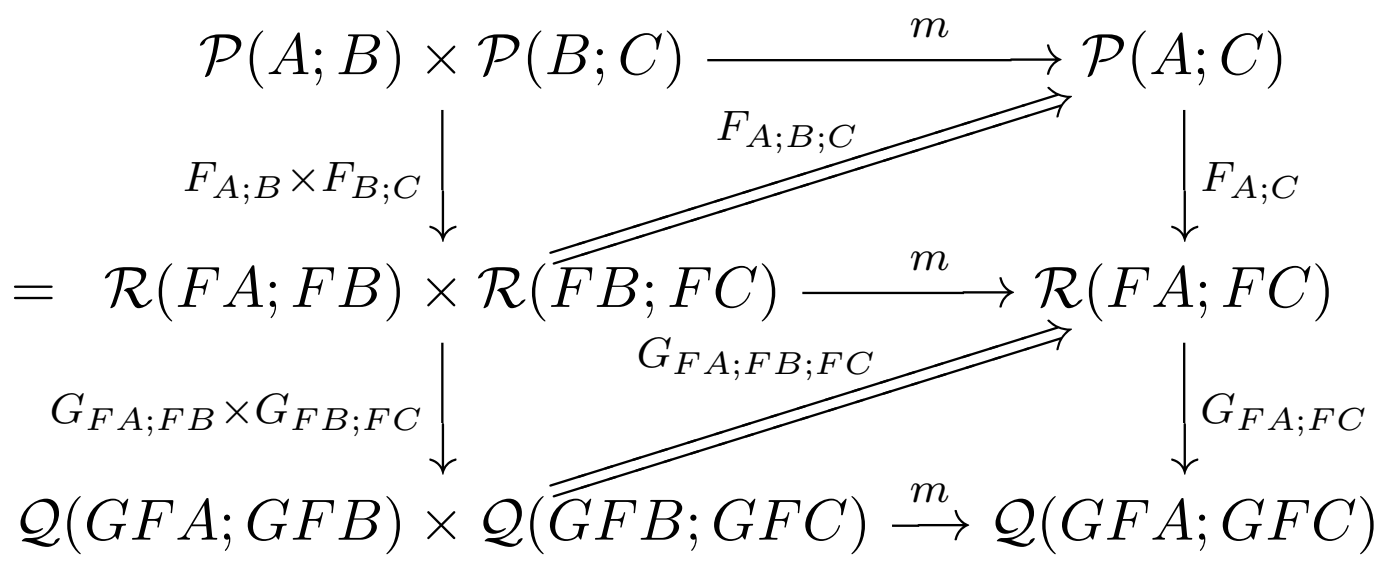

;

$$
\text{— } K_A = \kappa = \begin{array}{ccc}
\mathbb{1} & \xrightarrow{1_A} & \mathcal{P}(A;A) \\
\| & \overset{F_A}{\swarrow\!\!\swarrow} & \downarrow F_{A;A} \\
\mathbb{1} & \xrightarrow{1_{FA}} & \mathcal{R}(FA;FA) \\
\| & \overset{G_{FA}}{\swarrow\!\!\swarrow} & \downarrow G_{FA;FA} \\
\mathbb{1} & \xrightarrow{1_{GFA}} & \mathcal{Q}(GFA;GFA)
\end{array} .
$$

Diagram (1.3.1) for $K$ is a pasting of two diagrams (1.3.1) for $F$ and for $G$ stacked one on top of the other. Therefore, it commutes. Stacking one commutative cube at Figure 2 for $F$ on top of another for $G$, we prove that equation at Figure 2 holds for $(K,\chi)$. Stacking one commutative cube at Figure 3 for $F$ on top of another for $G$, we prove that equation at Figure 3 holds for $(K,\chi)$. Similarly, stacking one commutative prism (1.3.4) for $F$ on top of another for $G$, we get commutative prism (1.3.4) for $(K,\kappa)$. Likewise, stacking one commutative prism (1.3.5) for $F$ on top of another for $G$, we get commutative prism (1.3.5) for $(K,\kappa)$. Thus we have defined a composition law for morphisms of biprops.

The identity morphism $\mathrm{Id}\colon \mathcal{P} \to \mathcal{P}$ has identity map $\mathrm{Col}\,\mathrm{Id} = \mathrm{id}\colon \mathrm{Col}\,\mathcal{P} \to \mathrm{Col}\,\mathcal{P}$, identity functors $\mathrm{Id}_{A;B} = \mathrm{Id}\colon \mathcal{P}(A;B) \to \mathcal{P}(A;B)$, identity transformations $\mathrm{Id}_{A;B;C} = \mathrm{id}$ and $\mathrm{Id}_A = \mathrm{id}$. The identity morphisms are left and right units with respect to composition law for morphisms of biprops.

The composition law is associative. In fact, given $F\colon \mathcal{P} \to \mathcal{R}$, $G\colon \mathcal{R} \to \mathcal{Q}$, $E\colon \mathcal{Q} \to \mathcal{T}$ compose in two ways to $K\colon \mathcal{P} \to \mathcal{T}$ with

$$
\text{— } K_{A;B;C} = \begin{array}{ccc}
\mathcal{P}(A;B) \times \mathcal{P}(B;C) & \xrightarrow{\quad m \quad} & \mathcal{P}(A;C) \\
{\scriptstyle F_{A;B}\times F_{B;C}}\downarrow & \overset{F_{A;B;C}}{\nearrow\!\!\nearrow} & \downarrow {\scriptstyle F_{A;C}} \\
\mathcal{R}(FA;FB) \times \mathcal{R}(FB;FC) & \xrightarrow{\quad m \quad} & \mathcal{R}(FA;FC) \\
{\scriptstyle G_{FA;FB}\times G_{FB;FC}}\downarrow & \overset{G_{FA;FB;FC}}{\nearrow\!\!\nearrow} & \downarrow {\scriptstyle G_{FA;FC}} \\
\mathcal{Q}(GFA;GFB) \times \mathcal{Q}(GFB;GFC) & \xrightarrow{\quad m \quad} & \mathcal{Q}(GFA;GFC) \\
{\scriptstyle E_{GFA;GFB}\times E_{GFB;GFC}}\downarrow & \overset{E_{GFA;GFB;GFC}}{\nearrow\!\!\nearrow} & \downarrow {\scriptstyle E_{GFA;GFC}} \\
\mathcal{T}(EGFA;EGFB) \times \mathcal{T}(EGFB;EGFC) & \xrightarrow{m} & \mathcal{Q}(EGFA;EGFC)
\end{array} ;
$$

$$
\text{— } K_A = \begin{array}{ccc}
\mathbb{1} & \xrightarrow{1_A} & \mathcal{P}(A;A) \\
\| & \overset{F_A}{\swarrow\!\!\swarrow} & \downarrow F_{A;A} \\
\mathbb{1} & \xrightarrow{1_{FA}} & \mathcal{R}(FA;FA) \\
\| & \overset{G_{FA}}{\swarrow\!\!\swarrow} & \downarrow G_{FA;FA} \\
\mathbb{1} & \xrightarrow{1_{GFA}} & \mathcal{Q}(GFA;GFA) \\
\| & \overset{E_{GFA}}{\swarrow\!\!\swarrow} & \downarrow E_{GFA;GFA} \\
\mathbb{1} & \xrightarrow{1_{EGFA}} & \mathcal{T}(EGFA;EGFA)
\end{array} .
$$

Summing up, there is a category of biprops and their morphisms.

## 1.4 Biprop envelope of a symmetric weak multifunctor

Let $F\colon \mathsf{C} \to \mathsf{D}$ be a symmetric weak multifunctor. We produce from it a morphism of biprops $\mathbb{F}F\colon \mathbb{F}\mathsf{C} \to \mathbb{F}\mathsf{D}$ with $\mathrm{Col}\,\mathbb{F}F = \mathrm{Ob}\,F$, with functors

$$
(\mathbb{F}F)_{(X_i)_{i\in I};(Y_j)_{j\in J}} = \coprod_{\phi\colon I\to J\in\mathcal{S}_{\mathsf{sk}}} \prod_{j\in J} F_{(X_i)_{i\in\phi^{-1}j};Y_j} : \coprod_{\phi\colon I\to J\in\mathcal{S}_{\mathsf{sk}}} \prod_{j\in J} \mathsf{C}\big((X_i)_{i\in\phi^{-1}j};Y_j\big)
$$

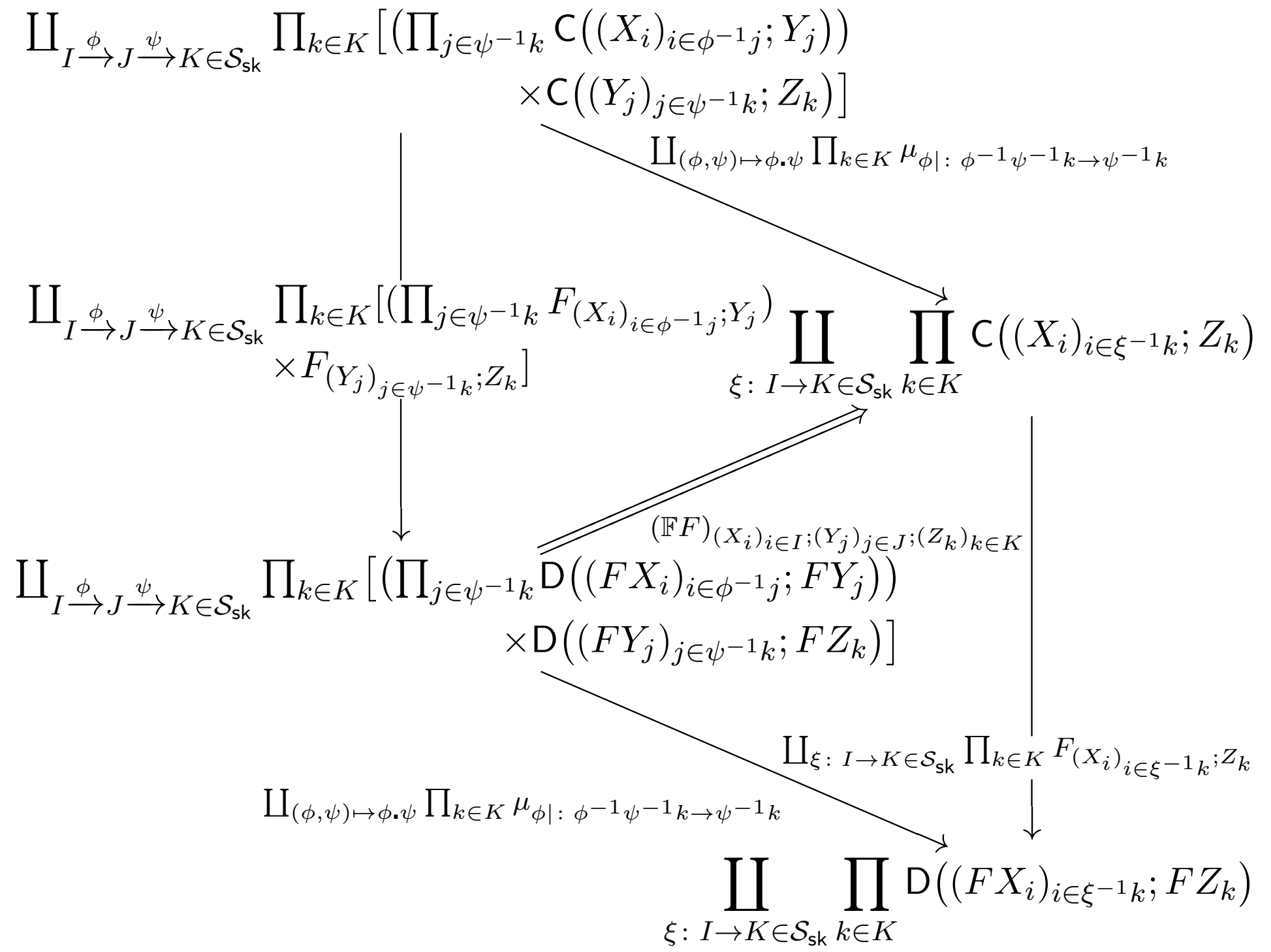

defined as whiskering, the transformation

$$\coprod_{I \xrightarrow{\phi} J \xrightarrow{\psi} K} \prod_{k\in K} F^{\phi|:\, \phi^{-1}\psi^{-1}k\to\psi^{-1}k}_{(X_i)_{i\in\phi^{-1}\psi^{-1}k};(Y_j)_{j\in\psi^{-1}k};Z_k},$$

followed by the functor

$$\coprod_{(\phi,\psi)\mapsto\phi\centerdot\psi} : \coprod_{I \xrightarrow{\phi} J \xrightarrow{\psi} K} \prod_{k\in K} \mathsf{D}\big((FX_i)_{i\in\phi^{-1}\psi^{-1}k}; FZ_k\big) \to \coprod_{\xi:\, I\to K} \prod_{k\in K} \mathsf{D}\big((FX_i)_{i\in\xi^{-1}k}; FZ_k\big),$$

Figure 4: Natural isomorphism for $\mathbb{F}F$

$$\to \coprod_{\phi:\, I\to J\in\mathcal{S}_{\mathsf{sk}}} \prod_{j\in J} \mathsf{D}\big((FX_i)_{i\in\phi^{-1}j}; FY_j\big).$$

natural isomorphism given on Figure 4, natural isomorphism

$(\mathbb{F}F)_{(X_i)_{i\in I}} =$

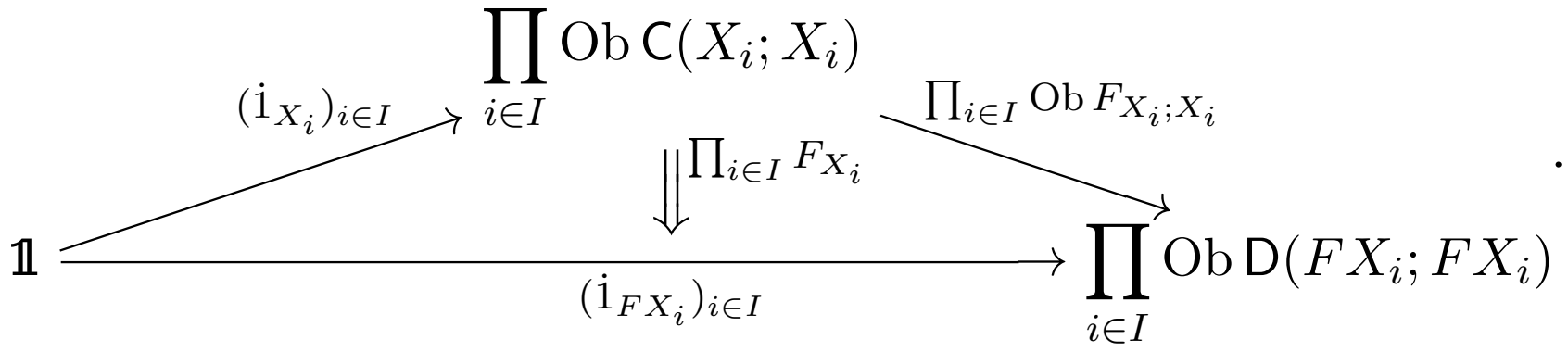

.

**1.4.1 Proposition.** *So defined* $\mathbb{F}F\colon \mathbb{F}\mathsf{C} \to \mathbb{F}\mathsf{D}$ *is a morphism of biprops.*

*Proof.* Let us prove equalities (1.3.1)–(1.3.5). Fix some maps $I \xrightarrow{f} L \xleftarrow{g} J \in \mathcal{S}_{\mathsf{sk}}$. Consider objects $A_i = (X_i^s)_{s\in S_i}$, $B_j = (Y_j^q)_{q\in Q_j}$. Introduce ordered sets $S = \bigsqcup_{i\in I} S_i$, $Q = \bigsqcup_{j\in J} Q_j$. A map $\phi : S \to Q \in \mathcal{S}_{\mathsf{sk}}$ such that diagram (1.2.1) (with $L$ in place of $K$) commutes amounts to family of maps $\phi_l\colon \bigsqcup_{i\in f^{-1}l} S_i \to \bigsqcup_{j\in g^{-1}l} Q_j \in \mathcal{S}_{\mathsf{sk}}$, $l \in L$. We have

$$\prod_{l\in L} \mathbb{F}\mathsf{C}\big(\bigsqcup_{i\in f^{-1}l} A_i; \bigsqcup_{j\in g^{-1}l} B_j\big) = \prod_{l\in L} \mathbb{F}\mathsf{C}\Big(\big((X_i^s)_{s\in S_i}\big)_{i\in f^{-1}l}; \big((Y_j^q)_{q\in Q_j}\big)_{j\in g^{-1}l}\Big)$$

$$= \prod_{l\in L} \bigsqcup_{\phi_l:\, \sqcup_{i\in f^{-1}l} S_i \to \sqcup_{j\in g^{-1}l} Q_j} \prod_{q\in\sqcup_{j\in g^{-1}l} Q_j} \mathsf{C}\big((X^s_{\mathrm{pr}_I s})_{s\in\phi_l^{-1}q}; Y^q_{\mathrm{pr}_J q}\big)$$

$$= \coprod_{\phi\colon S\to Q}^{\phi\centerdot \mathrm{pr}_J\centerdot g=\mathrm{pr}_I\centerdot f} \prod_{l\in L} \prod_{q\in(\mathrm{pr}_J g)^{-1}l} \mathsf{C}\big((X^s_{\mathrm{pr}_I s})_{s\in\phi^{-1}q}; Y^q_{\mathrm{pr}_J q}\big).$$

Thus equation (1.3.1) takes the form of an identity

$$\begin{array}{ccc}
\displaystyle\coprod_{\phi\colon S\to Q}^{\phi\centerdot \mathrm{pr}_J\centerdot g=\mathrm{pr}_I\centerdot f} \prod_{l\in L} \prod_{q\in(\mathrm{pr}_J g)^{-1}l} \mathsf{C}\big((X^s_{\mathrm{pr}_I s})_{s\in\phi^{-1}q}; Y^q_{\mathrm{pr}_J q}\big) & \xrightarrow{\coprod_{\phi\mapsto\phi} 1} & \displaystyle\coprod_{\xi\colon S\to Q} \prod_{q\in Q} \mathsf{C}\big((X^s_{\mathrm{pr}_I s})_{s\in\xi^{-1}q}; Y^q_{\mathrm{pr}_J q}\big) \\
\Big\downarrow {\scriptstyle \coprod_{\phi\colon S\to Q}^{\phi\centerdot \mathrm{pr}_J\centerdot g=\mathrm{pr}_I\centerdot f} \prod_{l\in L}\prod_{q\in(\mathrm{pr}_J g)^{-1}l} F_{(X^s_{\mathrm{pr}_I s})_{s\in\phi^{-1}q}; Y^q_{\mathrm{pr}_J q}}} & = & \Big\downarrow {\scriptstyle \coprod_{\xi\colon S\to Q}\prod_{q\in Q} F_{(X^s_{\mathrm{pr}_I s})_{s\in\xi^{-1}q}; Y^q_{\mathrm{pr}_J q}}} \\
\displaystyle\coprod_{\phi\colon S\to Q}^{\phi\centerdot \mathrm{pr}_J\centerdot g=\mathrm{pr}_I\centerdot f} \prod_{l\in L} \prod_{q\in(\mathrm{pr}_J g)^{-1}l} \mathsf{D}\big((FX^s_{\mathrm{pr}_I s})_{s\in\phi^{-1}q}; FY^q_{\mathrm{pr}_J q}\big) & \xrightarrow{\coprod_{\phi\mapsto\phi} 1} & \displaystyle\coprod_{\xi\colon S\to Q} \prod_{q\in Q} \mathsf{D}\big((FX^s_{\mathrm{pr}_I s})_{s\in\xi^{-1}q}; FY^q_{\mathrm{pr}_J q}\big)
\end{array}$$

Now let us prove equality (1.3.2). Write relevant objects as $A_i = (X^s_i)_{s\in S_i}$, $B_j = (Y^q_j)_{q\in Q_j}$, $C_k = (Z^p_k)_{p\in P_k}$. Denote $S = \bigsqcup_{i\in I} S_i$, $Q = \bigsqcup_{j\in J} Q_j$, $P = \bigsqcup_{k\in K} P_k$. We consider pairs of maps $S \xrightarrow{\phi} Q \xrightarrow{\psi} P$ such that diagram (1.2.3) commutes. Up to certain isomorphisms the left hand side of equation (1.3.2) equals

$$\coprod_{(\phi,\psi)\mapsto\phi\centerdot\psi} \left( \coprod_{I\xrightarrow{f}L\xleftarrow{g}J}^{\substack{\phi\centerdot\mathrm{pr}_J\centerdot g=\mathrm{pr}_I\centerdot f \\ \psi\centerdot\mathrm{pr}_K\centerdot h=\mathrm{pr}_J\centerdot g}} \prod_{p\in P} F^{\phi|\colon(\phi\psi)^{-1}p\to\psi^{-1}p}_{(X^s_{\mathrm{pr}_I s})_{s\in(\phi\psi)^{-1}p};(Y^q_{\mathrm{pr}_J q})_{q\in\psi^{-1}p};Z^p_{\mathrm{pr}_K p}} \right),$$

while the right hand side (up to certain isomorphisms) is

$$\coprod_{(\phi,\psi)\mapsto\phi\centerdot\psi} \left( \coprod_{I\xrightarrow{f}L\xleftarrow{g}J}^{\substack{\phi\centerdot\mathrm{pr}_J\centerdot g=\mathrm{pr}_I\centerdot f \\ \psi\centerdot\mathrm{pr}_K\centerdot h=\mathrm{pr}_J\centerdot g}} \prod_{l\in L} \prod_{p\in(\mathrm{pr}_K h)^{-1}l} F^{\phi|\colon(\phi\psi)^{-1}p\to\psi^{-1}p}_{(X^s_{\mathrm{pr}_I s})_{s\in(\phi\psi)^{-1}p};(Y^q_{\mathrm{pr}_J q})_{q\in\psi^{-1}p};Z^p_{\mathrm{pr}_K p}} \right).$$

These expressions are identified, hence equation (1.3.2) holds true.

Equation (1.3.3) follows from the disjoint sum over $I \xrightarrow{\phi} J \xrightarrow{\psi} K \xrightarrow{\xi} L \in \mathcal{S}_{\mathsf{sk}}$ of products over $l \in L$ of equations [Lyu25, Figure 2]. Equation (1.3.4) follows from the disjoint sum over $\phi\colon I \to J \in \mathcal{S}_{\mathsf{sk}}$ of products over $j \in J$ of equations [Lyu25, (1.1.16)]. Equation (1.3.5) follows from the disjoint sum over $\psi\colon J \to K \in \mathcal{S}_{\mathsf{sk}}$ of products over $k \in K$ of equations [Lyu25, (1.1.17)]. □

Let $F\colon \mathsf{C} \to \mathsf{D}$, $G\colon \mathsf{D} \to \mathsf{E}$ be symmetric weak multifunctors. Their composition $K = F\centerdot G\colon \mathsf{C} \to \mathsf{E}$ [Lyu25, (1.1.21) and around] has $\mathrm{Ob}\,K = \big(\mathrm{Ob}\,\mathsf{C} \xrightarrow{\mathrm{Ob}\,F} \mathrm{Ob}\,\mathsf{D} \xrightarrow{\mathrm{Ob}\,G} \mathrm{Ob}\,\mathsf{E}\big)$, $K_{A;Y} = F_{A;Y} \centerdot G_{FA;FY}$, where $A = (X_i)_{i\in I}$, $K^\phi$ is a pasting of $F^\phi$ and $G^\phi$, and $K_X$ is a pasting of $F_X$ and $G_{FX}$. We claim that $\mathbb{F}(F\centerdot G) = (\mathbb{F}F)\centerdot(\mathbb{F}G)$. In fact,

$$\mathrm{Col}\,\mathbb{F}(F\centerdot G) = \big(\mathrm{Ob}\,\mathsf{C} \xrightarrow{\mathrm{Ob}\,F} \mathrm{Ob}\,\mathsf{D} \xrightarrow{\mathrm{Ob}\,G} \mathrm{Ob}\,\mathsf{E}\big) = \mathrm{Col}(\mathbb{F}F)\centerdot\mathrm{Col}(\mathbb{F}G),$$

$$\begin{aligned}
\big[(\mathbb{F}F)\centerdot(\mathbb{F}G)\big]_{(X_i)_{i\in I};(Y_j)_{j\in J}} &= (\mathbb{F}F)_{(X_i)_{i\in I};(Y_j)_{j\in J}} \centerdot (\mathbb{F}G)_{(FX_i)_{i\in I};(FY_j)_{j\in J}} \\
&= \Big( \coprod_{\phi\colon I\to J} \prod_{j\in J} F_{(X_i)_{i\in\phi^{-1}j};Y_j} \Big) \centerdot \coprod_{\phi\colon I\to J} \prod_{j\in J} G_{(FX_i)_{i\in\phi^{-1}j};FY_j} \\
&= \coprod_{\phi\colon I\to J} \prod_{j\in J} \big[ F_{(X_i)_{i\in\phi^{-1}j};Y_j} \centerdot G_{(FX_i)_{i\in\phi^{-1}j};FY_j} \big] = \big[\mathbb{F}(F\centerdot G)\big]_{(X_i)_{i\in I};(Y_j)_{j\in J}},
\end{aligned}$$

$$\big[(\mathbb{F}F)\centerdot(\mathbb{F}G)\big]_{(X_i)_{i\in I};(Y_j)_{j\in J};(Z_k)_{k\in K}} =$$

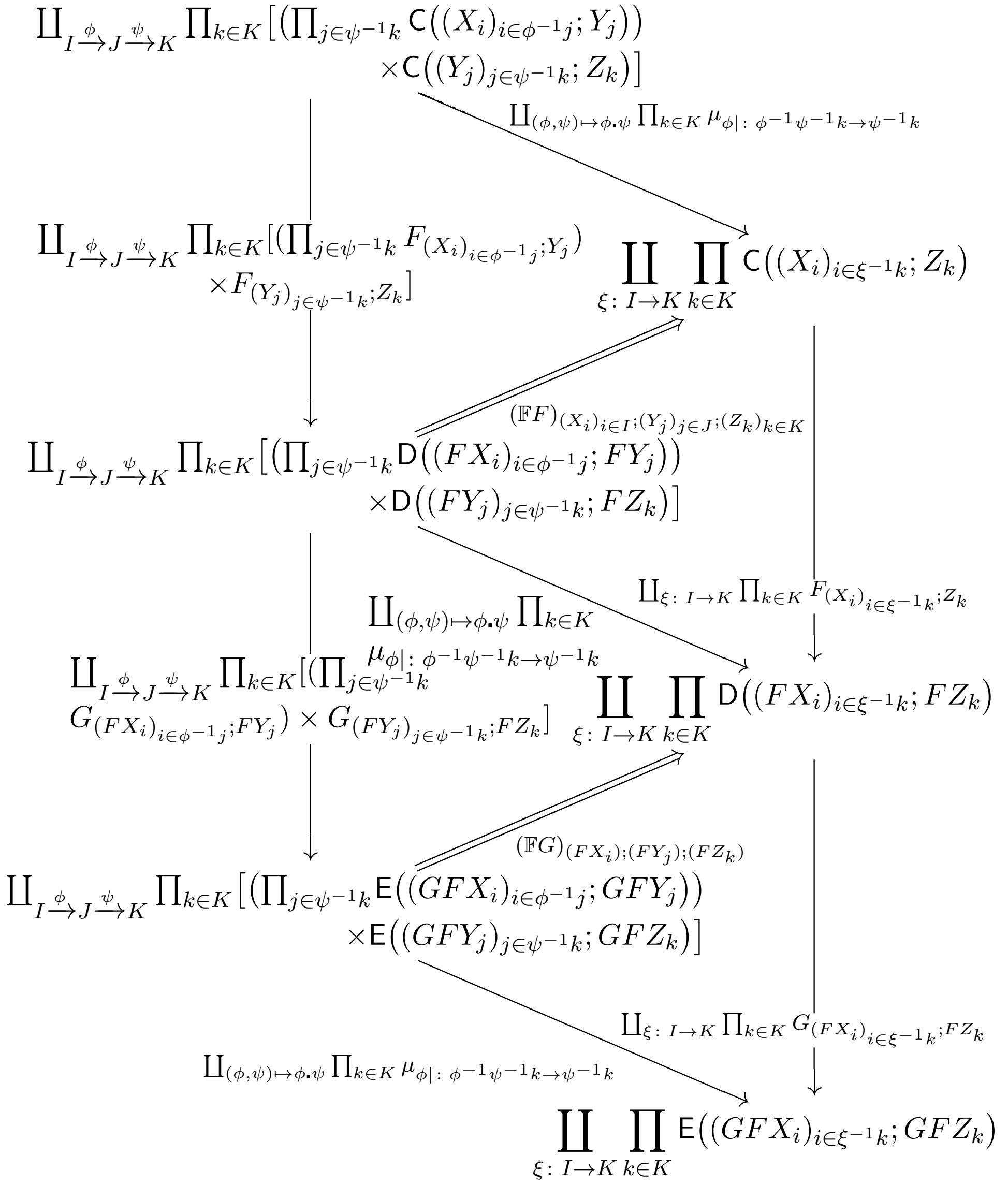

$$= \boxed{\dfrac{\coprod_{(\phi,\psi)\mapsto\phi\centerdot\psi}\Big(\coprod_{I\xrightarrow{\phi}J\xrightarrow{\psi}K}\prod_{k\in K}F^{\phi|:\,\phi^{-1}\psi^{-1}k\to\psi^{-1}k}_{(X_i)_{i\in\phi^{-1}\psi^{-1}k};(Y_j)_{j\in\psi^{-1}k};Z_k}\Big)}{\coprod_{(\phi,\psi)\mapsto\phi\centerdot\psi}\Big(\coprod_{I\xrightarrow{\phi}J\xrightarrow{\psi}K}\prod_{k\in K}G^{\phi|:\,\phi^{-1}\psi^{-1}k\to\psi^{-1}k}_{(FX_i)_{i\in\phi^{-1}\psi^{-1}k};(FY_j)_{j\in\psi^{-1}k};FZ_k}\Big)}}$$

$$= \coprod_{(\phi,\psi)\mapsto\phi\centerdot\psi}\Big(\coprod_{I\xrightarrow{\phi}J\xrightarrow{\psi}K}\prod_{k\in K}(F\centerdot G)^{\phi|:\,\phi^{-1}\psi^{-1}k\to\psi^{-1}k}_{(X_i)_{i\in\phi^{-1}\psi^{-1}k};(Y_j)_{j\in\psi^{-1}k};Z_k}\Big)$$

$$= \big[\mathbb{F}(F\centerdot G)\big]_{(X_i)_{i\in I};(Y_j)_{j\in J};(Z_k)_{k\in K}}.$$

Figure 5: The composition is preserved by $\mathbb{F}$

The proof that the composition is preserved by $\mathbb{F}$ is completed on the previous page.

$\mathbb{F}$ takes the identity multifunctor $\mathrm{Id}\colon \mathsf{C} \to \mathsf{C}$ to the identity morphism of biprops $\mathrm{Id}\colon \mathbb{F}\mathsf{C} \to \mathbb{F}\mathsf{C}$. Indeed, $\mathrm{Col}\,\mathbb{F}\,\mathrm{Id} = \mathrm{id}$, $(\mathbb{F}\,\mathrm{Id})_{(X_i)_{i\in I};(Y_j)_{j\in J}} = \mathrm{Id}$, $(\mathbb{F}\,\mathrm{Id})_{(X_i)_{i\in I};(Y_j)_{j\in J};(Z_k)_{k\in K}} = \mathrm{id}$, $(\mathbb{F}\,\mathrm{Id})_{(X_i)_{i\in I}} = \mathrm{id}$.

Therefore, $\mathbb{F}$ is a functor from the category of symmetric weak multicategories and symmetric weak multifunctors to the category of biprops and their morphisms.

# 2 Actions of symmetric groups on a biprop

Similarly to Leinster [Lei03, Lemma A.2.2] we derive two actions of the symmetric group on hom-categories of a biprop, rather than including the action in the definition of a biprop. The first by permuting the arguments, the second by permuting the outputs. These actions are obtained from the composition functors and from units. The symmetric group acts by functors. Naturally, these actions are not strict and the composition is preserved only up to isomorphism of functors (Proposition 2.1.2), which satisfies the non-abelian cocycle identity (Propositions 2.1.5 and 2.1.6).

## 2.1 Symmetric group action functors

Let $\beta\colon J \to K \in \mathcal{S}_{\mathsf{sk}}$ be a bijection. Let $(Y_j)_{j\in J}$, $(Z_k)_{k\in K}$, $(W_i)_{i\in I}$ be families of colours of a biprop $\mathcal{P}$ such that $Z_k = Y_{\beta^{-1}k}$. Similarly to [Lei03, Lemma A.2.2] and to [Lyu23, § A.1] define functors

$$\begin{aligned} l_\beta = \bigl\{\mathcal{P}\bigl((Z_k)_{k\in K};(W_i)_{i\in I}\bigr) &\xrightarrow{(\mathrm{i}_{Z_k})_{k\in K}\times 1} \bigl[\prod_{k\in K}\mathcal{P}(Y_{\beta^{-1}k};Z_k)\bigr]\times\mathcal{P}\bigl((Z_k)_{k\in K};(W_i)_{i\in I}\bigr) \\ &\xrightarrow{\otimes^{J\xrightarrow{\beta}K\xleftarrow{1}K}\times 1} \mathcal{P}\bigl((Y_j)_{j\in J};(Z_k)_{k\in K}\bigr)\times\mathcal{P}\bigl((Z_k)_{k\in K};(W_i)_{i\in I}\bigr) \xrightarrow{m} \mathcal{P}\bigl((Y_j)_{j\in J};(W_i)_{i\in I}\bigr)\bigr\}, 7 \end{aligned}$$

$$\begin{aligned} r_\beta = \bigl\{\mathcal{P}\bigl((W_i)_{i\in I};(Y_j)_{j\in J}\bigr) &\xrightarrow{1\times(\mathrm{i}_{Z_k})_{k\in K}} \mathcal{P}\bigl((W_i)_{i\in I};(Y_j)_{j\in J}\bigr)\times\prod_{k\in K}\mathcal{P}(Y_{\beta^{-1}k};Z_k) \\ &\xrightarrow{1\times\otimes^{J\xrightarrow{\beta}K\xleftarrow{1}K}} \mathcal{P}\bigl((W_i)_{i\in I};(Y_j)_{j\in J}\bigr)\times\mathcal{P}\bigl((Y_j)_{j\in J};(Z_k)_{k\in K}\bigr) \xrightarrow{m} \mathcal{P}\bigl((W_i)_{i\in I};(Z_k)_{k\in K}\bigr)\bigr\}. \end{aligned}$$

**2.1.1 Remark.** Let $\beta = 1\colon K \to K$. Then due to (1.1.3) there are isomorphisms $\mathsf{l}\colon l_{1_K} \to \mathrm{Id}$, $\mathsf{r}\colon r_{1_K} \to \mathrm{Id}$.

**2.1.2 Proposition.** *Assume that both $\alpha$ and $\beta$ are bijections from $\mathcal{S}_{\mathsf{sk}}$, $\gamma = (I \xrightarrow{\alpha} J \xrightarrow{\beta} K)$. Then there are isomorphisms*

$$\begin{aligned} \phi^l_{\alpha,\beta}\colon l_\gamma &\to \bigl[\mathcal{P}\bigl((Y_{\beta^{-1}k})_{k\in K};(W_l)_{l\in L}\bigr) \xrightarrow{l_\beta} \mathcal{P}\bigl((Y_j)_{j\in J};(W_l)_{l\in L}\bigr) \xrightarrow{l_\alpha} \mathcal{P}\bigl((Y_{\alpha i})_{i\in I};(W_l)_{l\in L}\bigr)\bigr], \\ \phi^r_{\alpha,\beta}\colon r_\gamma &\to \bigl[\mathcal{P}\bigl((W_l)_{l\in L};(Y_{\alpha i})_{i\in I}\bigr) \xrightarrow{r_\alpha} \mathcal{P}\bigl((W_l)_{l\in L};(Y_j)_{j\in J}\bigr) \xrightarrow{r_\beta} \mathcal{P}\bigl((W_l)_{l\in L};(Y_{\beta^{-1}k})_{k\in K}\bigr)\bigr]. \end{aligned}$$

*Proof.* Consider $X_i = Y_{\alpha i}$, hence, $Y_j = X_{\alpha^{-1}j}$ as well as $Z_k = Y_{\beta^{-1}k}$. The natural transformation

$(\phi^l_{\alpha,\beta})^{-1}$ is the composition

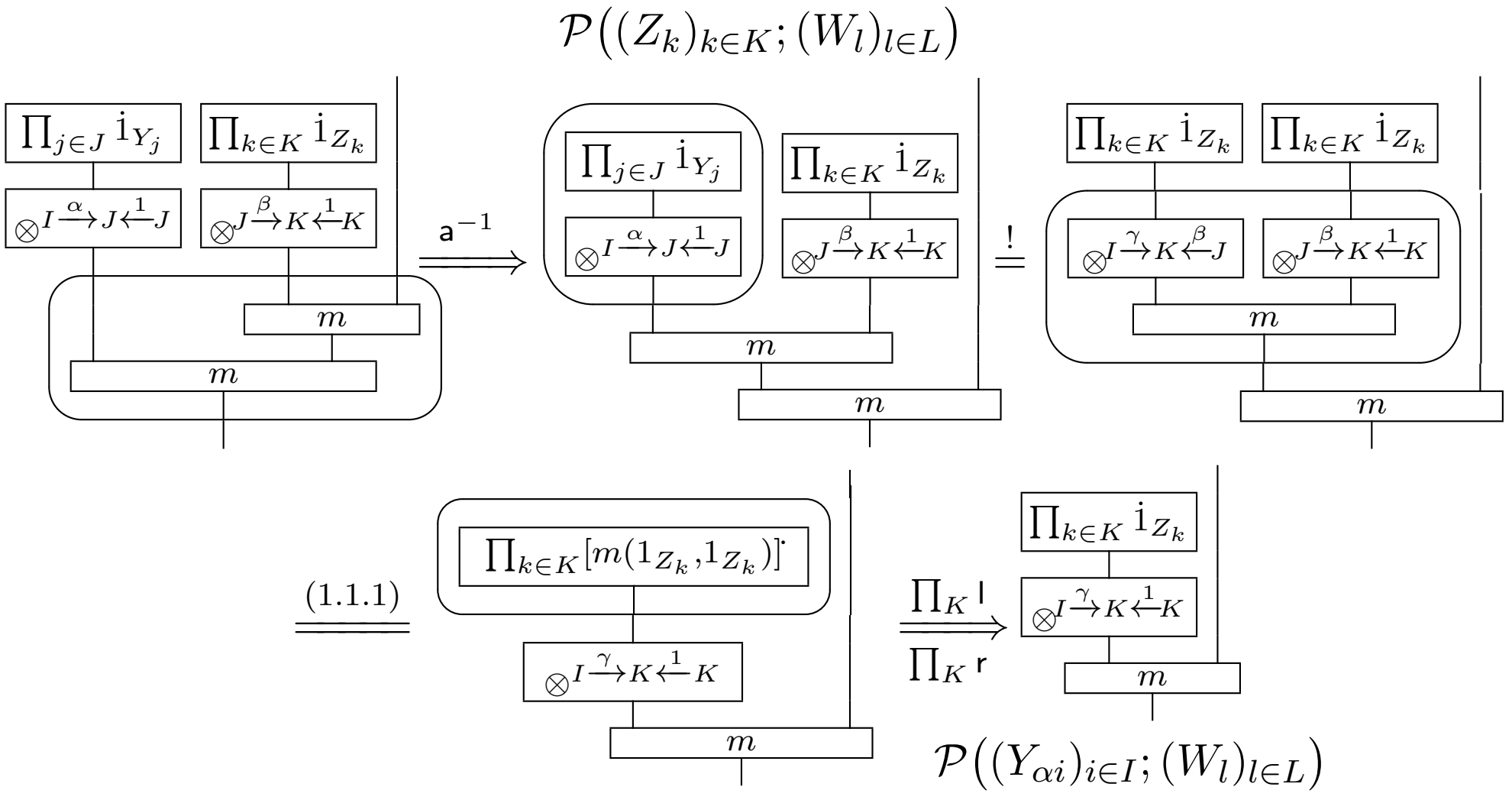

Equality marked ! follows from axiom (1.1.6) applied to namely,

$$\otimes^{I\xrightarrow{\alpha}J\xleftarrow{1}J} = \Big[\prod_{j\in J}\mathcal{P}(X_{\alpha^{-1}j};Y_j) \cong \prod_{k\in K}\mathcal{P}(X_{\gamma^{-1}k};Z_k) = \prod_{k\in K}\mathcal{P}(X_{\gamma^{-1}k};Y_{\beta^{-1}k}) \xrightarrow{\otimes^{I\xrightarrow{\gamma}K\xleftarrow{\beta}J}} \mathcal{P}\big((X_i)_{i\in I};(Y_j)_{j\in J}\big)\Big]. \quad (2.1.1)$$

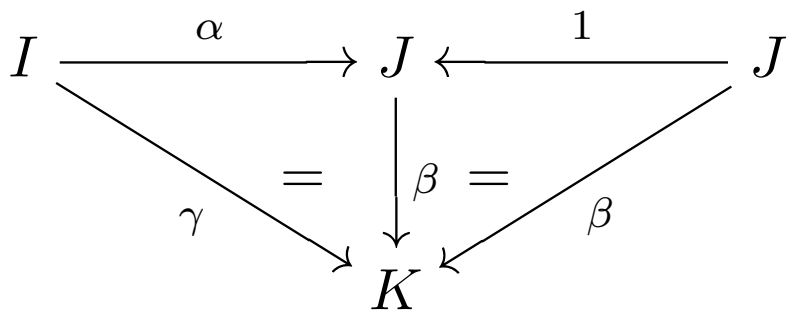

The natural transformation $(\phi^r_{\alpha,\beta})^{-1}$ is the composition

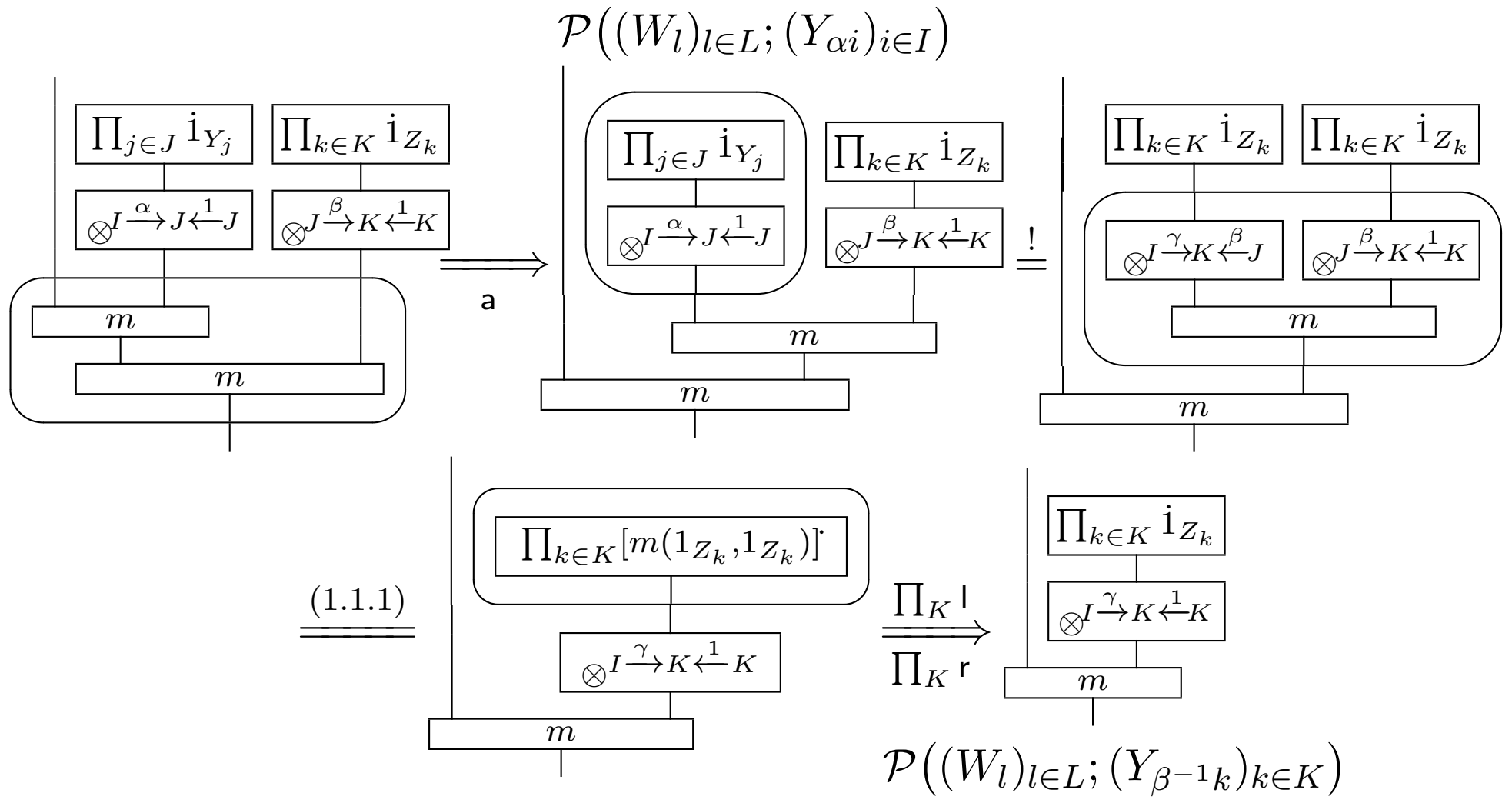

Equality marked ! follows from (2.1.1). ☐

**2.1.3 Corollary.** *For any bijection $\beta$ the functors $l_\beta$, $r_\beta$ are equivalences.*

*Proof.* Follows from Proposition 2.1.2 and Remark 2.1.1. ☐

**2.1.4 Proposition.** *Assume that both $\alpha\colon I \to J$ and $\beta\colon K \to L$ are bijections from $\mathcal{S}_{\mathsf{sk}}$. Then there is an isomorphism*

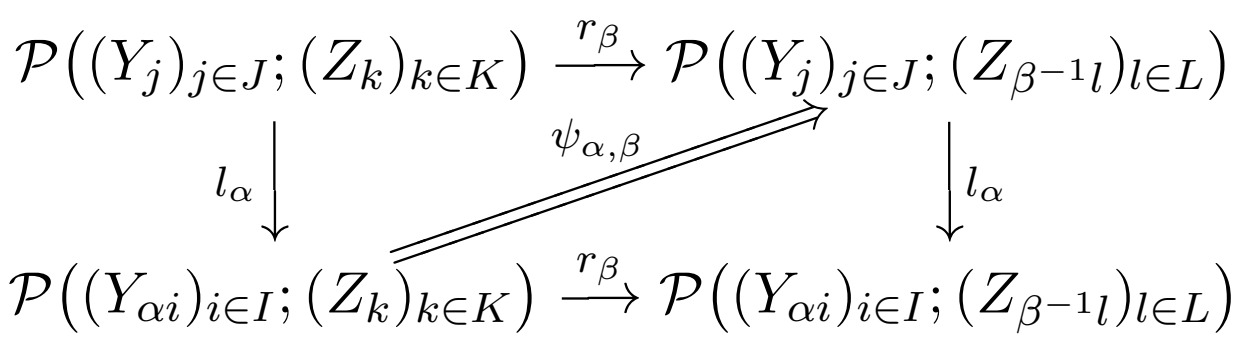

*Proof.* Denote $X_i = Y_{\alpha i}$, $W_l = Z_{\beta^{-1}l}$. Then there is the associativity isomorphism

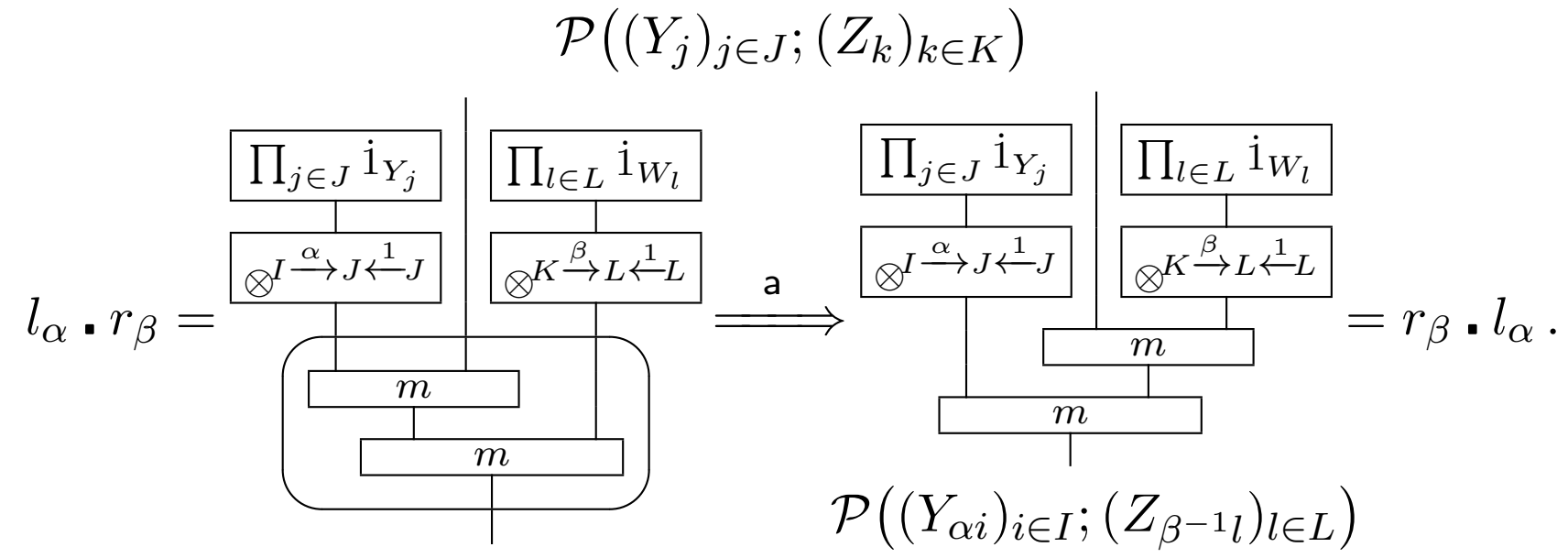

□

Below we prove the coherence condition asserting that the two ways of decomposing a left reindexing along a triple composite of bijections into three left reindexings agree.

**2.1.5 Proposition.** *Let $\alpha,\beta,\gamma$ from $I \xrightarrow{\alpha} J \xrightarrow{\beta} K \xrightarrow{\gamma} L \in \mathcal{S}_{\mathsf{sk}}$ be bijections. Let $(X_i)_{i\in I}$, $(Y_j)_{j\in J}$, $(Z_k)_{k\in K}$, $(W_l)_{l\in L}$, $(U_n)_{n\in N}$ be families of colours of $\mathcal{P}$ such that $X_i = Y_{\alpha i}$, $Y_j = Z_{\beta j}$, $Z_k = W_{\gamma k}$ for all $i \in I$, $j \in J$, $k \in K$. Then the non-abelian cocycle identity holds for $\phi^l$:*

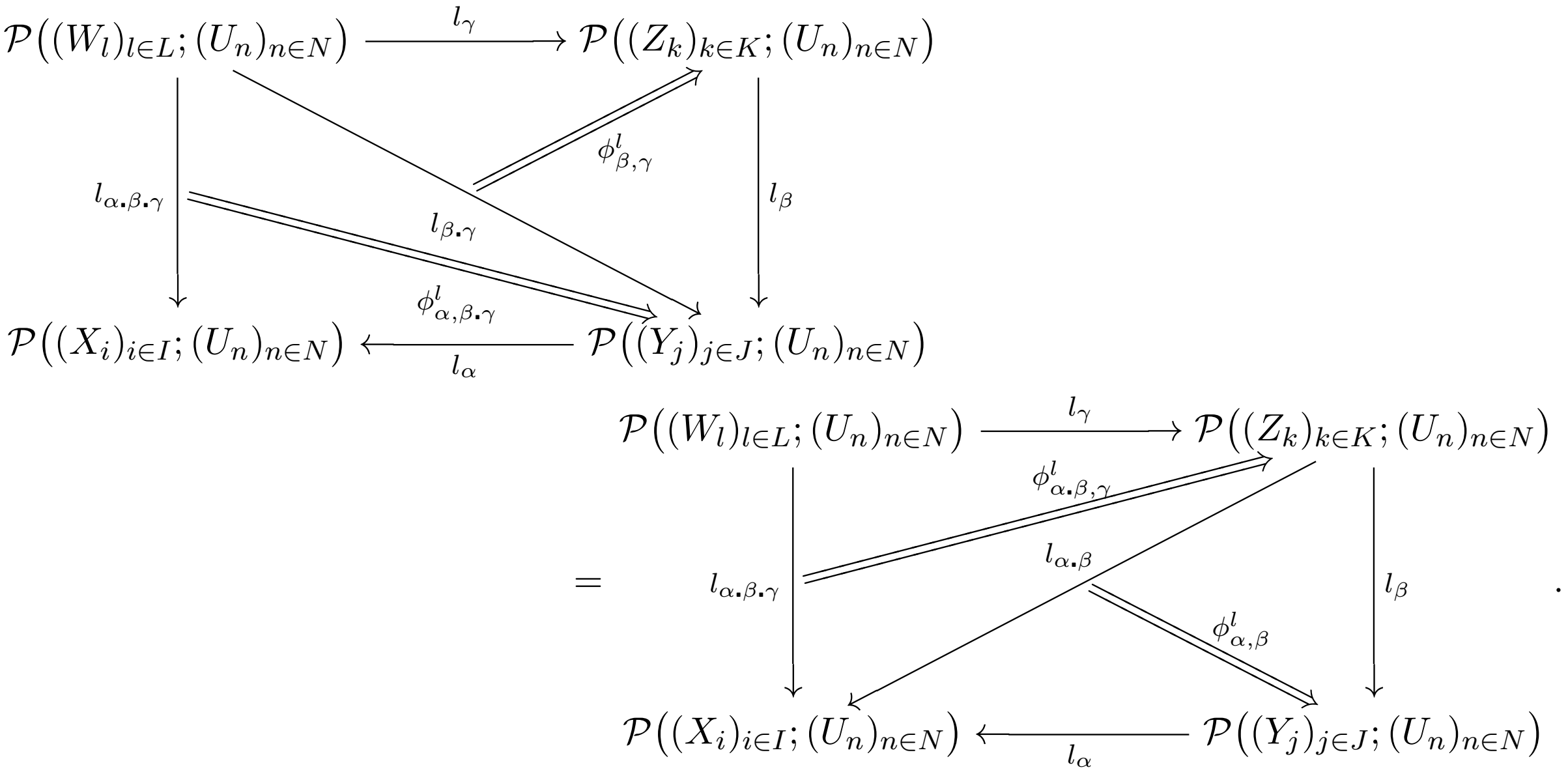

*Proof.* The left hand side is the composition

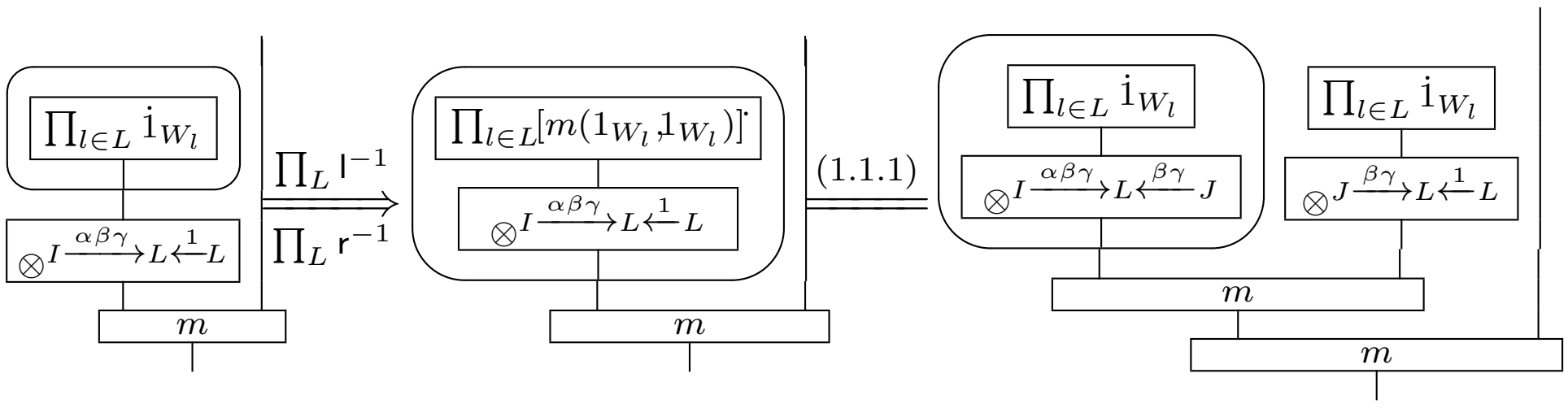

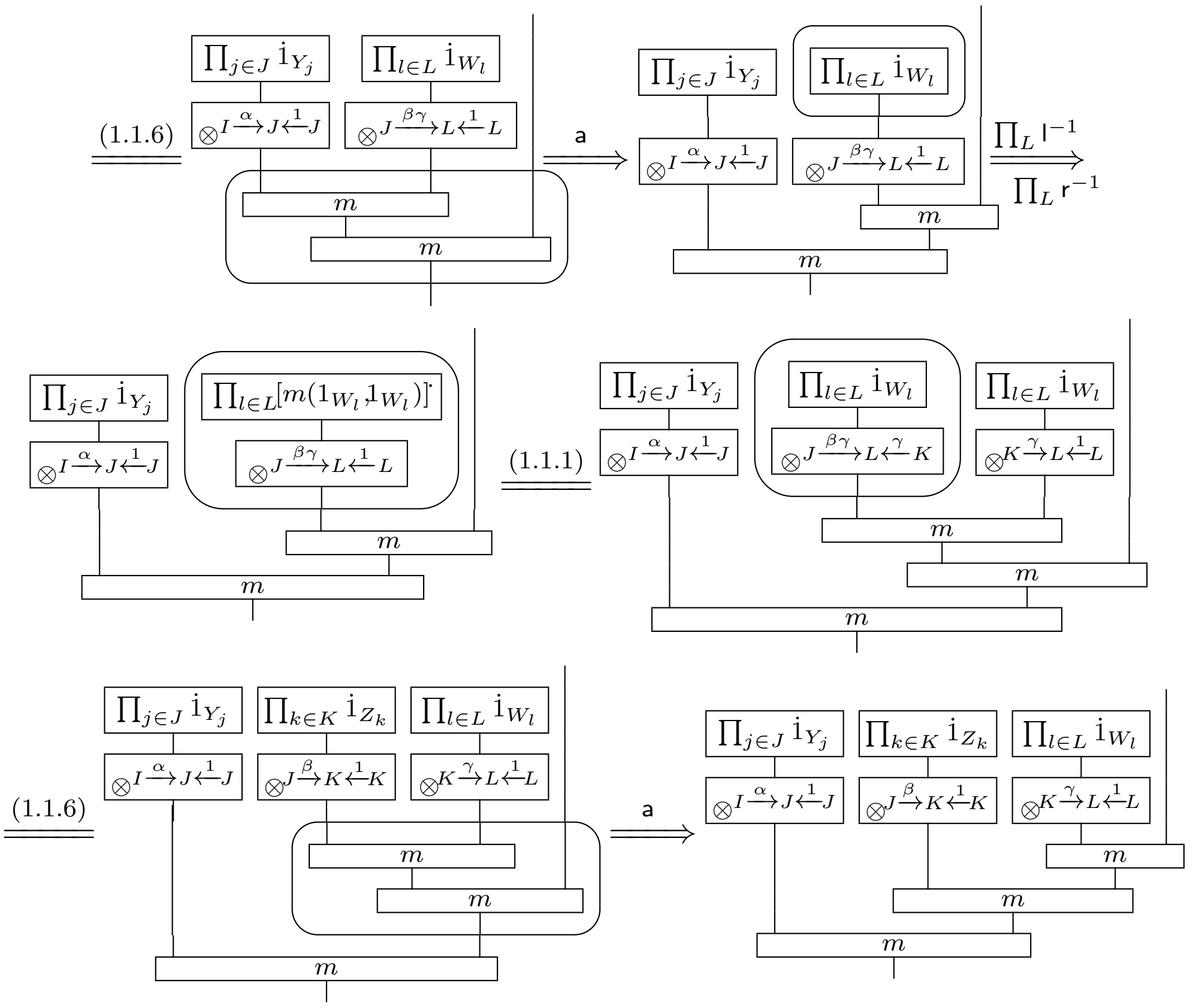

For the first arrow (and other similar arrows) we may use either $\prod_L \mathsf{l}^{-1}$ or $\prod_L \mathsf{r}^{-1}$. The resulting transformations are equal as [Lyu25, Lemma 2.1.4] shows. The right hand side is the composition

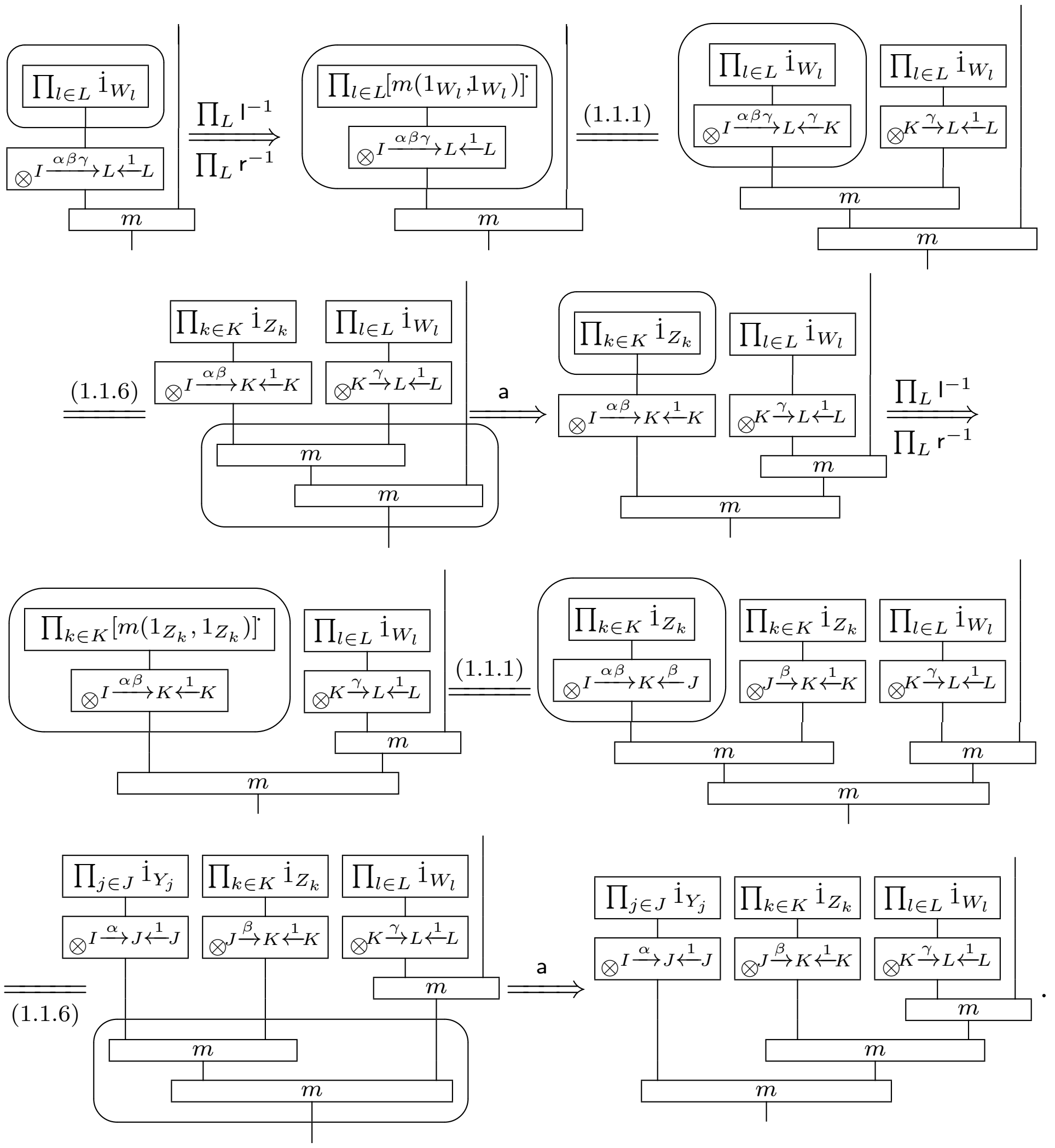

The equality follows from identity (1.1.8) and from associativity pentagon for the associators $\mathsf{a}$. $\square$

Below we prove the coherence condition asserting that the two ways of decomposing a right reindexing along a triple composite of bijections into three right reindexings agree.

**2.1.6 Proposition.** *Let $\alpha, \beta, \gamma$ from $I \xrightarrow{\alpha} J \xrightarrow{\beta} K \xrightarrow{\gamma} L \in \mathcal{S}_{\mathsf{sk}}$ be bijections. Let $(X_i)_{i\in I}$, $(Y_j)_{j\in J}$, $(Z_k)_{k\in K}$, $(W_l)_{l\in L}$, $(U_n)_{n\in N}$ be families of colours of $\mathcal{P}$ such that $X_i = Y_{\alpha i}$, $Y_j = Z_{\beta j}$, $Z_k = W_{\gamma k}$ for all $i \in I$, $j \in J$, $k \in K$. Then the non-abelian cocycle identity holds for $\phi^r$:*

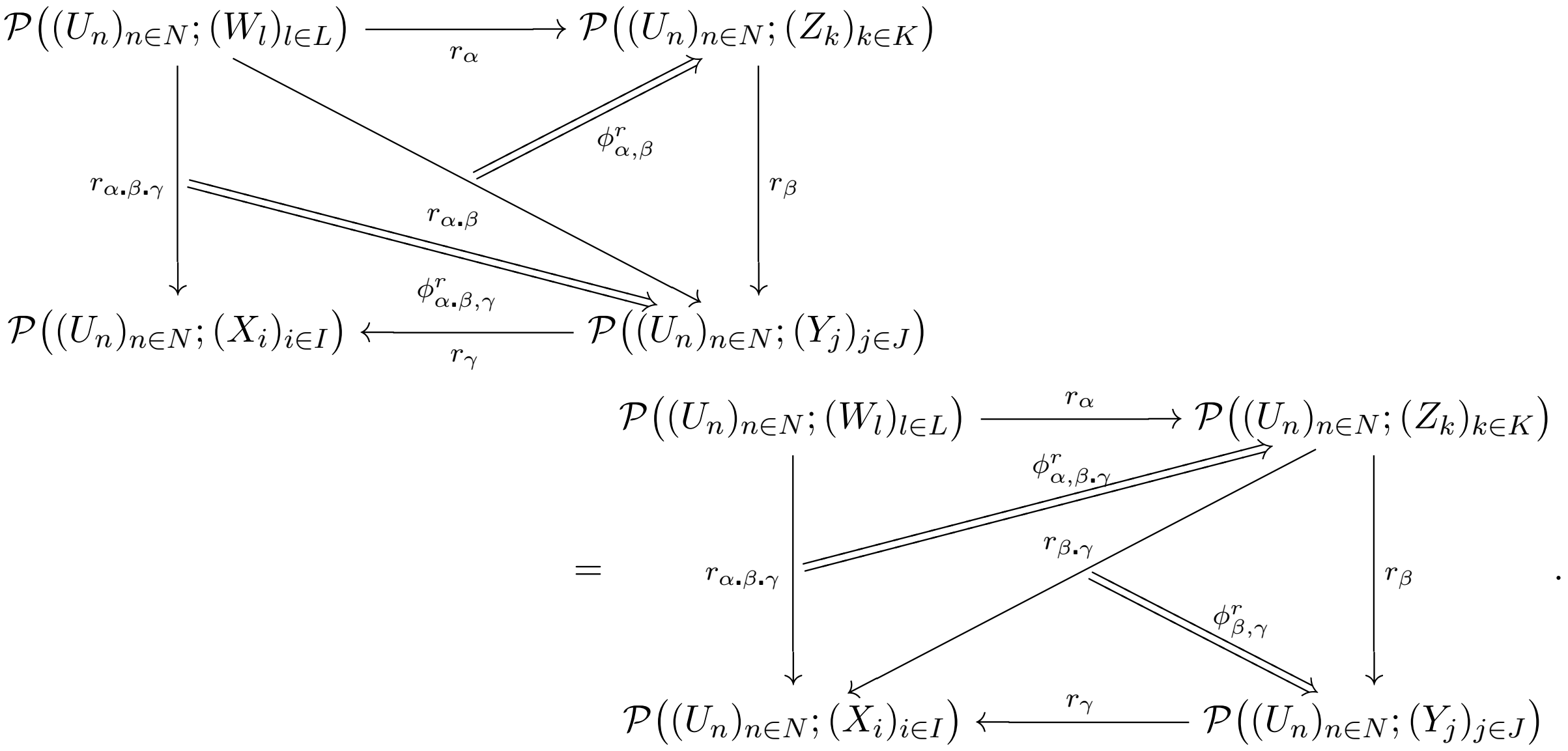

*Proof.* The left hand side is the composition

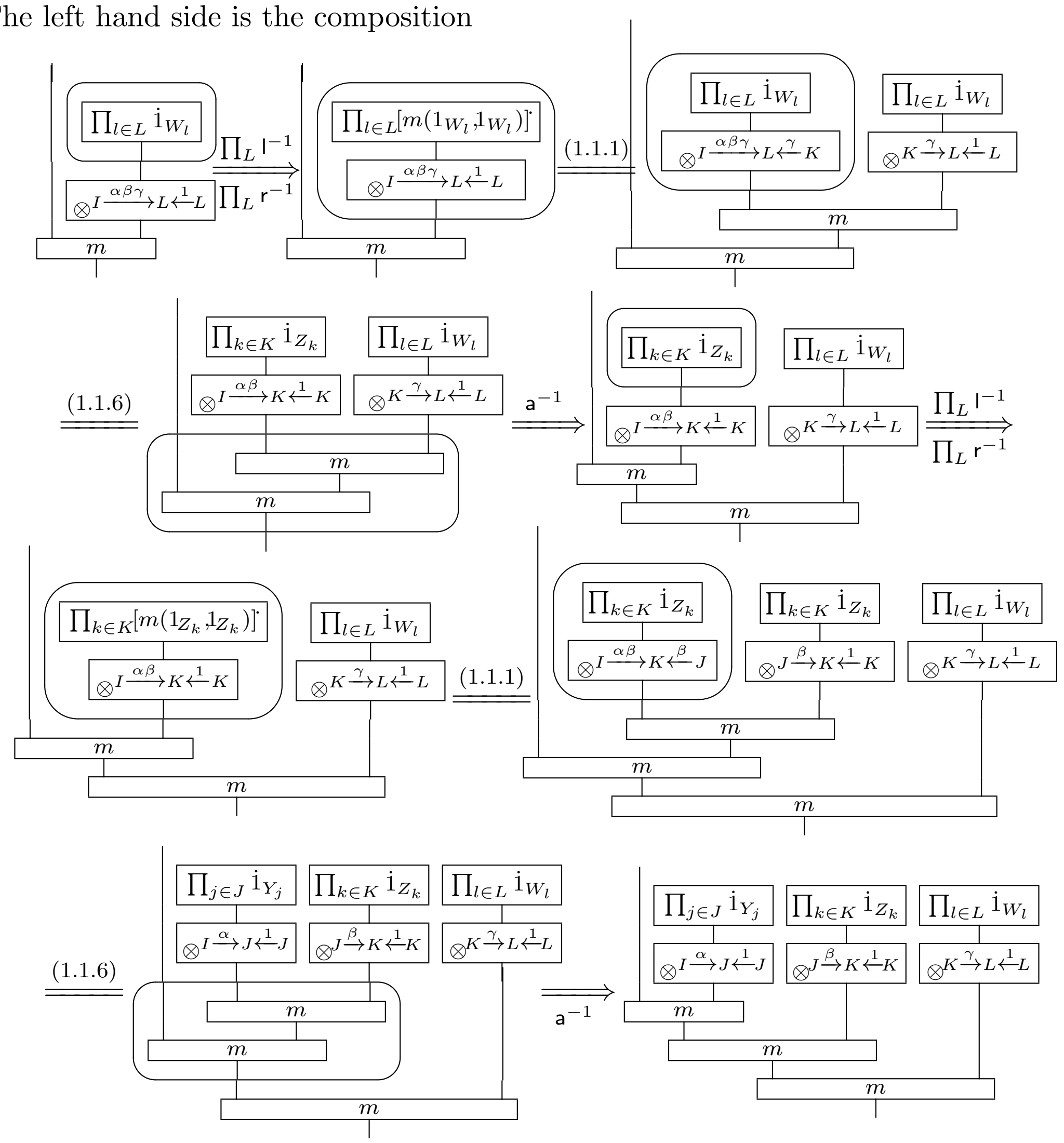

The right hand side is the composition

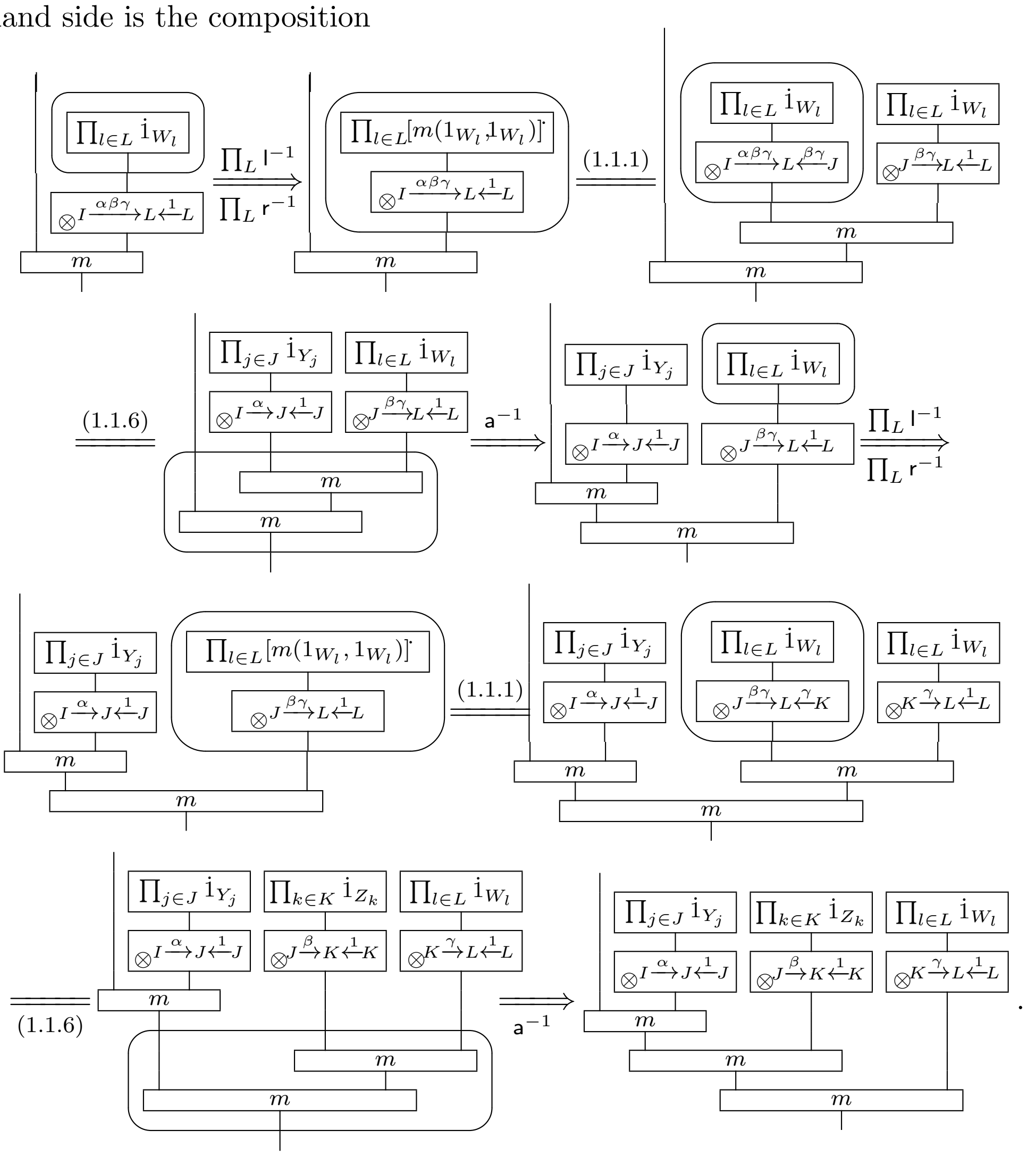

The equality follows from identity (1.1.8) and from associativity pentagon for the associators $\mathsf{a}$. $\square$

**2.1.7 Remark.** Let $\beta\colon J \to K \in \mathcal{S}_{\mathsf{sk}}$ be a bijection. Let $(X_i)_{i\in I}$, $(Y_j)_{j\in J}$, $(Z_k)_{k\in K}$, $(W_n)_{n\in N}$ be families of colours of a biprop $\mathcal{P}$ such that $Z_k = Y_{\beta^{-1}k}$. Then there is the associativity isomorphism

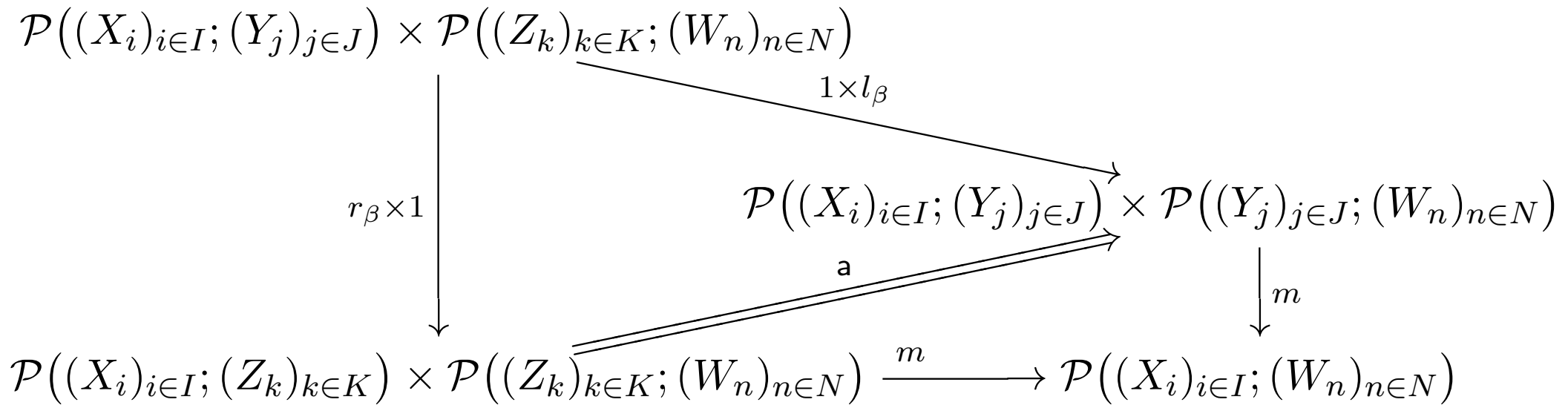

*Proof.* Indeed, this is the isomorphism

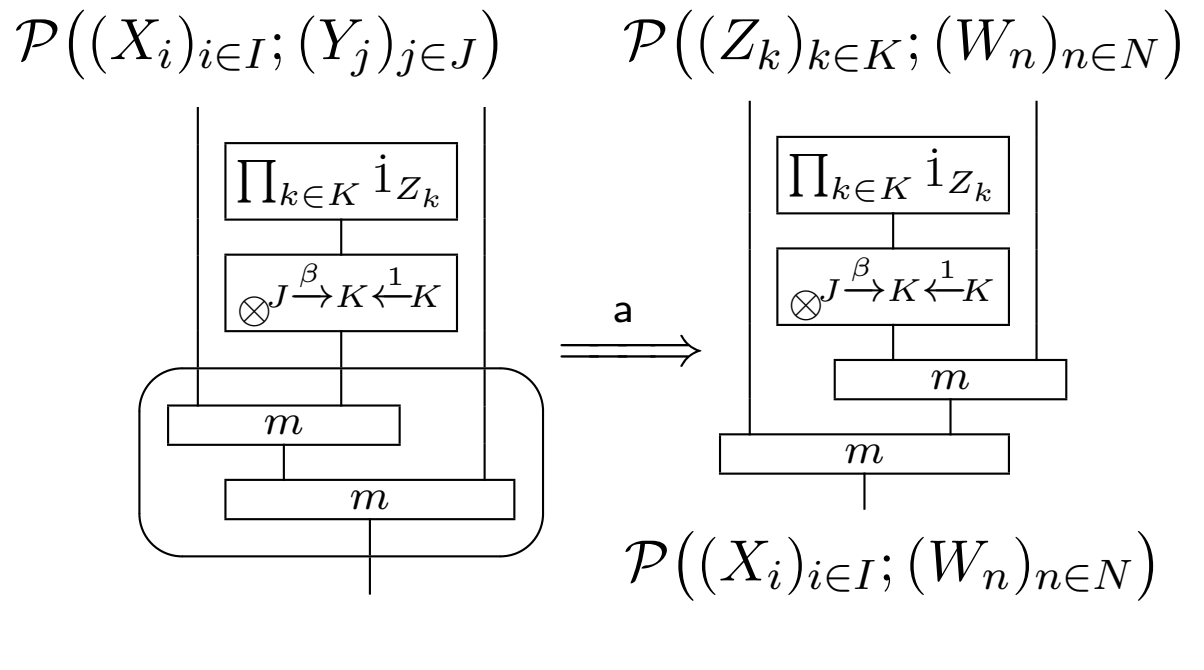

$\square$

**2.1.8 Remark.** Let $\beta\colon N \to I \in \mathcal{S}_{\mathsf{sk}}$ be a bijection. Let $(X_i)_{i\in I}$, $(Y_j)_{j\in J}$, $(Z_k)_{k\in K}$, $(W_n)_{n\in N}$ be families of colours of a biprop $\mathcal{P}$ such that $X_i = W_{\beta^{-1}i}$. Then there is the associativity isomorphism

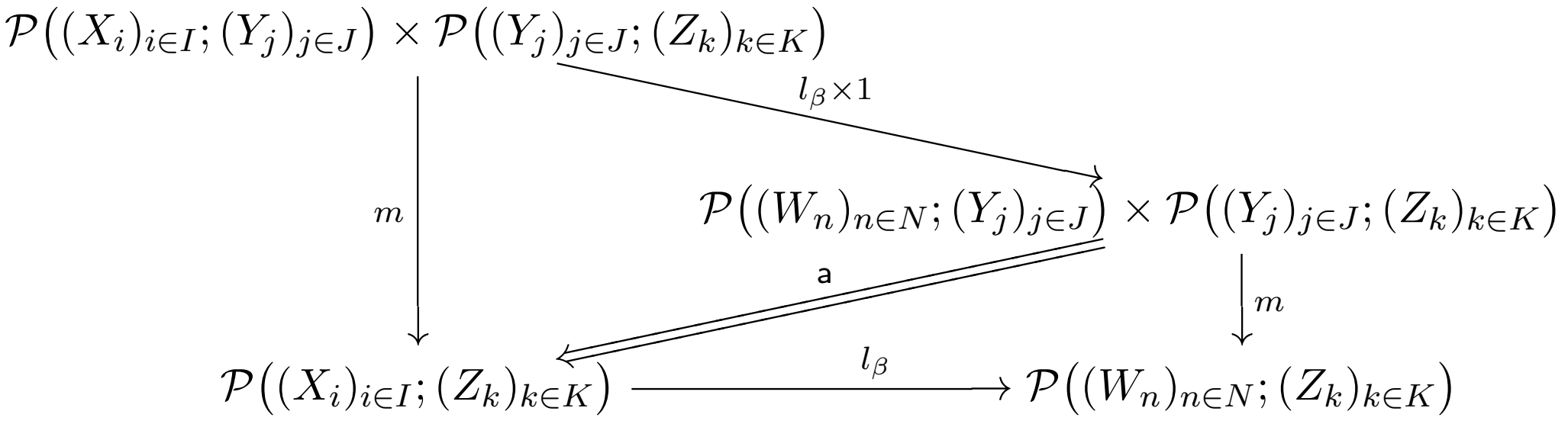

*Proof.* Indeed, this is the isomorphism

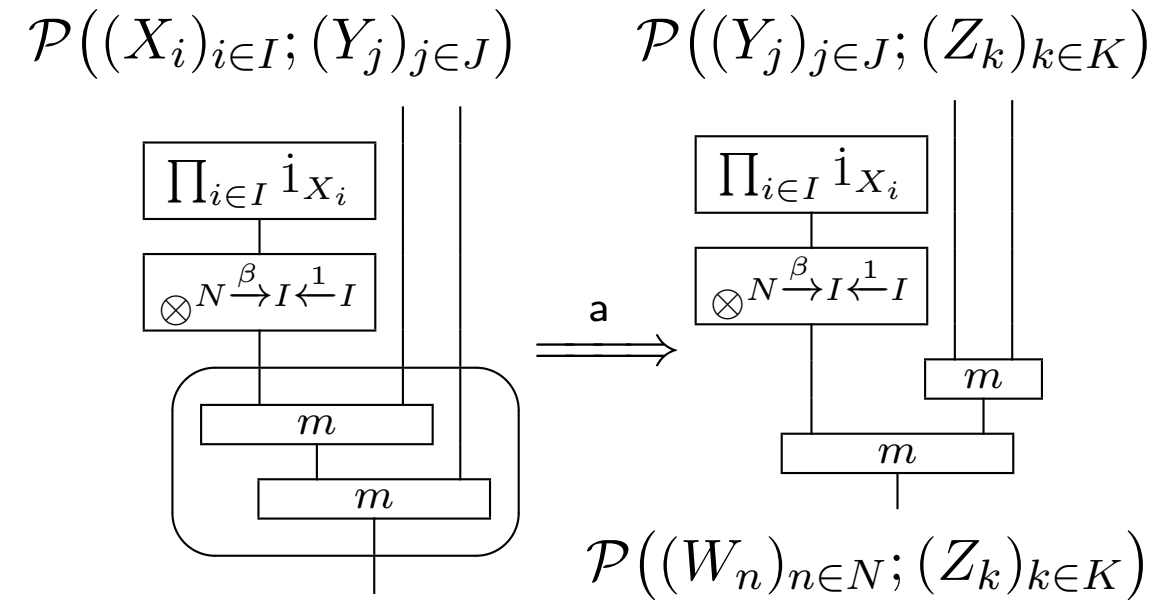

□

**2.1.9 Remark.** Let $\beta\colon K \to N \in \mathcal{S}_{\mathsf{sk}}$ be a bijection. Let $(X_i)_{i\in I}$, $(Y_j)_{j\in J}$, $(Z_k)_{k\in K}$, $(W_n)_{n\in N}$ be families of colours of a biprop $\mathcal{P}$ such that $W_n = Z_{\beta^{-1}n}$. Then there is the associativity isomorphism

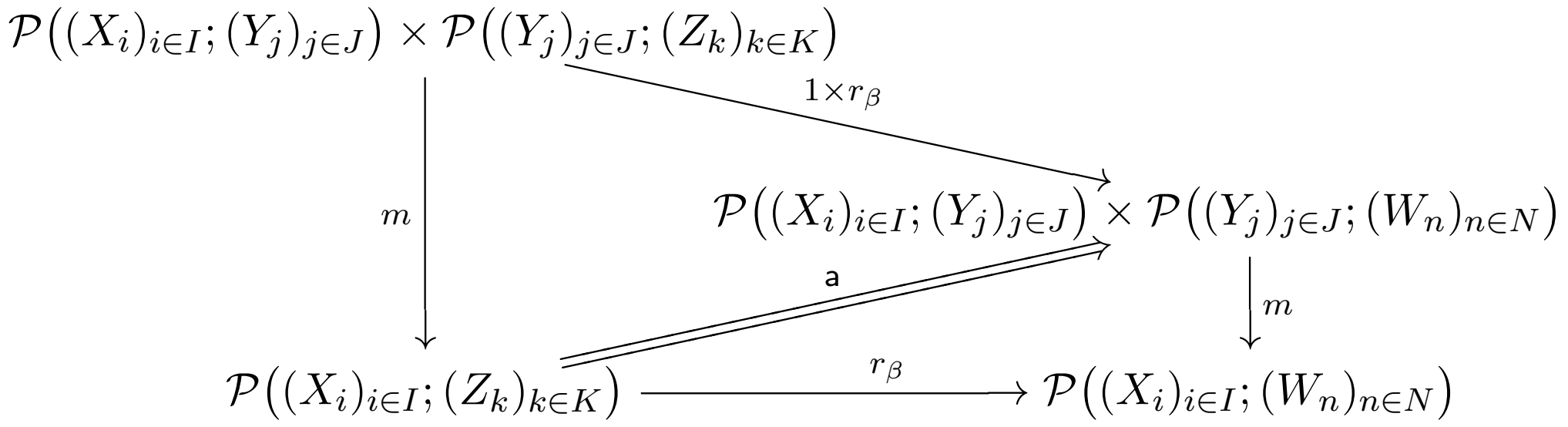

*Proof.* Indeed, this is the isomorphism

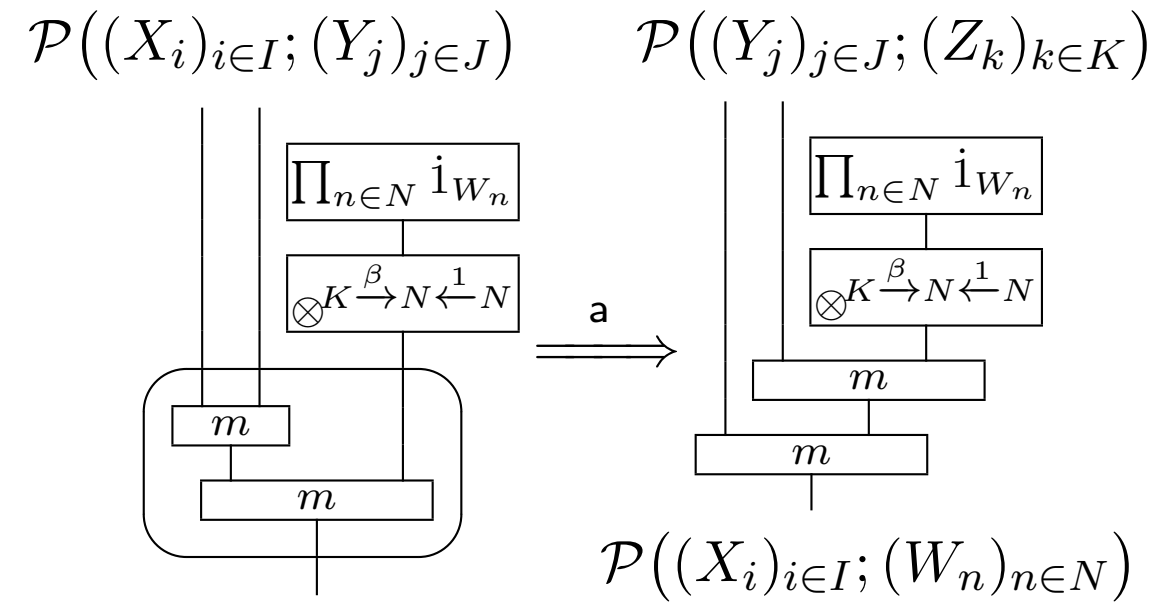

□

## 2.2 Strong natural transformations

Strong is a synonym for weak more often used nowadays. Let $F^k, G^l\colon \mathcal{P} \to \mathcal{R}$, $k \in K$, $l \in L$ be morphisms of biprops. A *(strong) natural transformation* $t\colon (F^k)_{k\in K} \to (G^l)_{l\in L}\colon \mathcal{P} \to \mathcal{R}$ consists of

— objects $t_X \in \operatorname{Ob} \mathcal{R}\big((F^k X)_{k\in K}; (G^l X)_{l\in L}\big)$, $X \in \operatorname{Col} \mathcal{P}$;

— natural isomorphisms

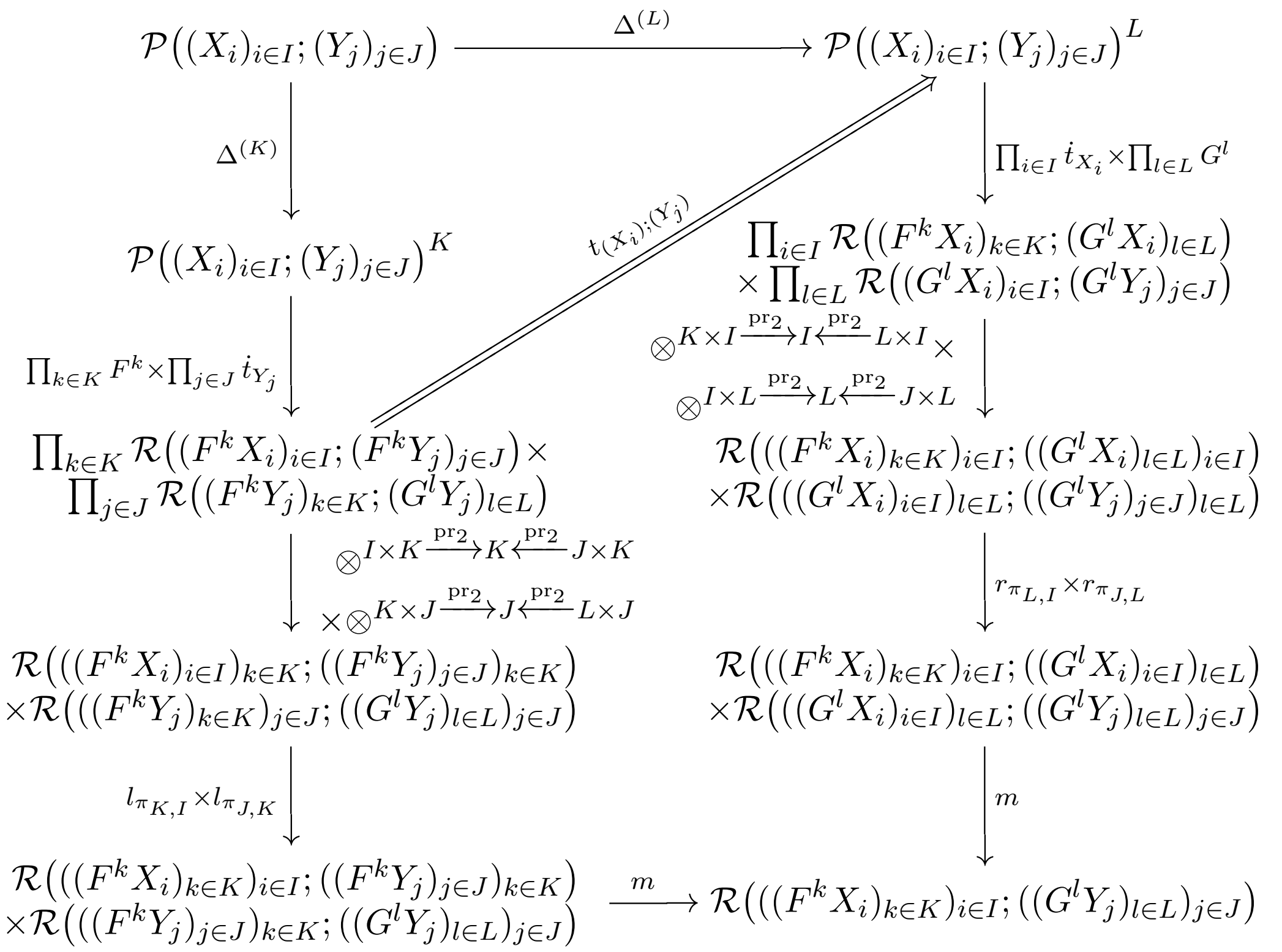

(2.2.1)

Here $\Delta^{(L)}\colon x\mapsto(x,x,\dots,x)$ is the diagonal functor. In Bartlett's notation we write

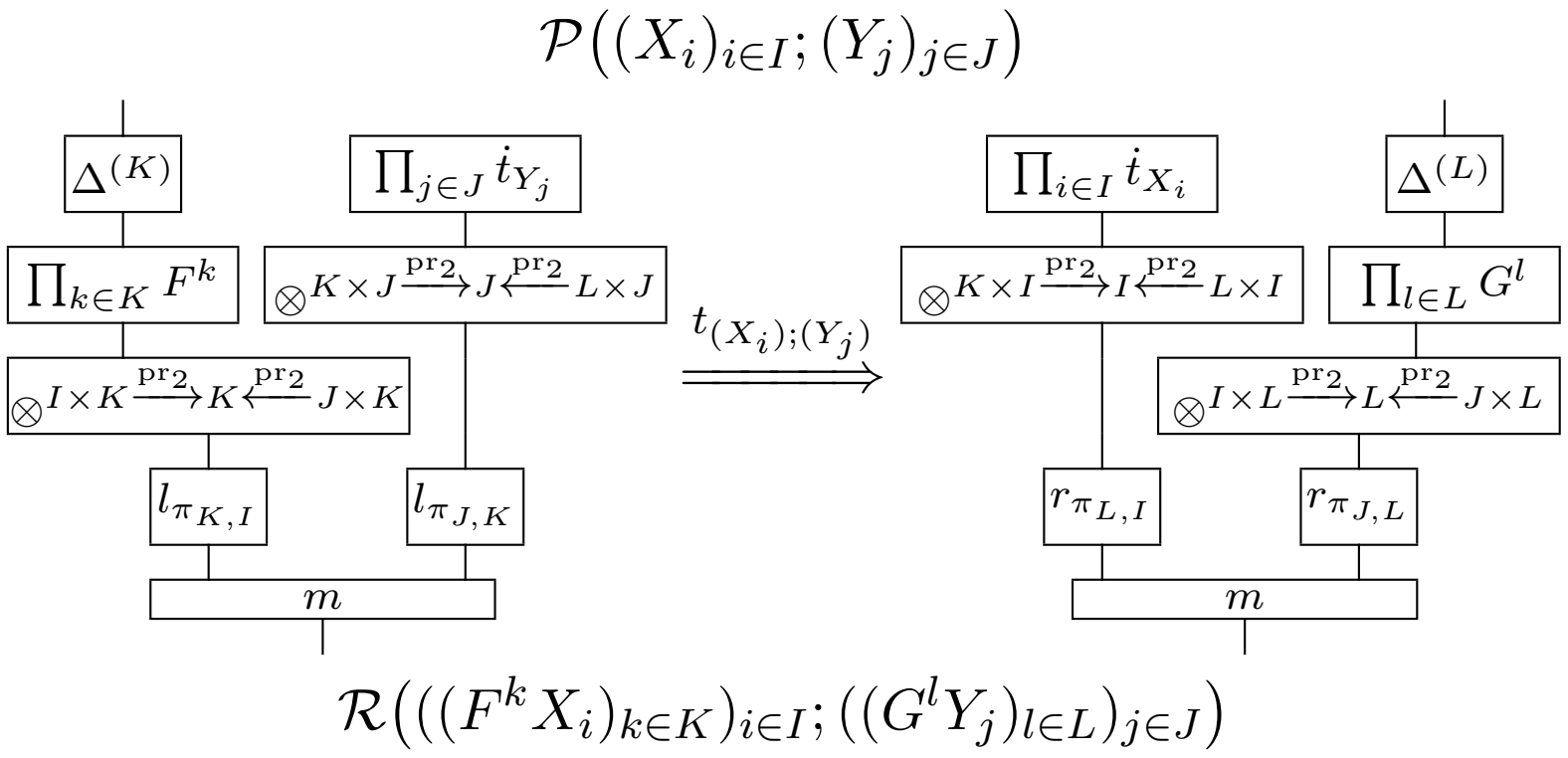

These transformations have to agree with compositions in the sense of the following equation, where $X_i, Y_j, Z_n$ are colours of $\mathcal{P}$, thus, $A=(X_i)_{i\in I}$, $B=(Y_j)_{j\in J}$, $C=(Z_n)_{n\in N}$ are objects:

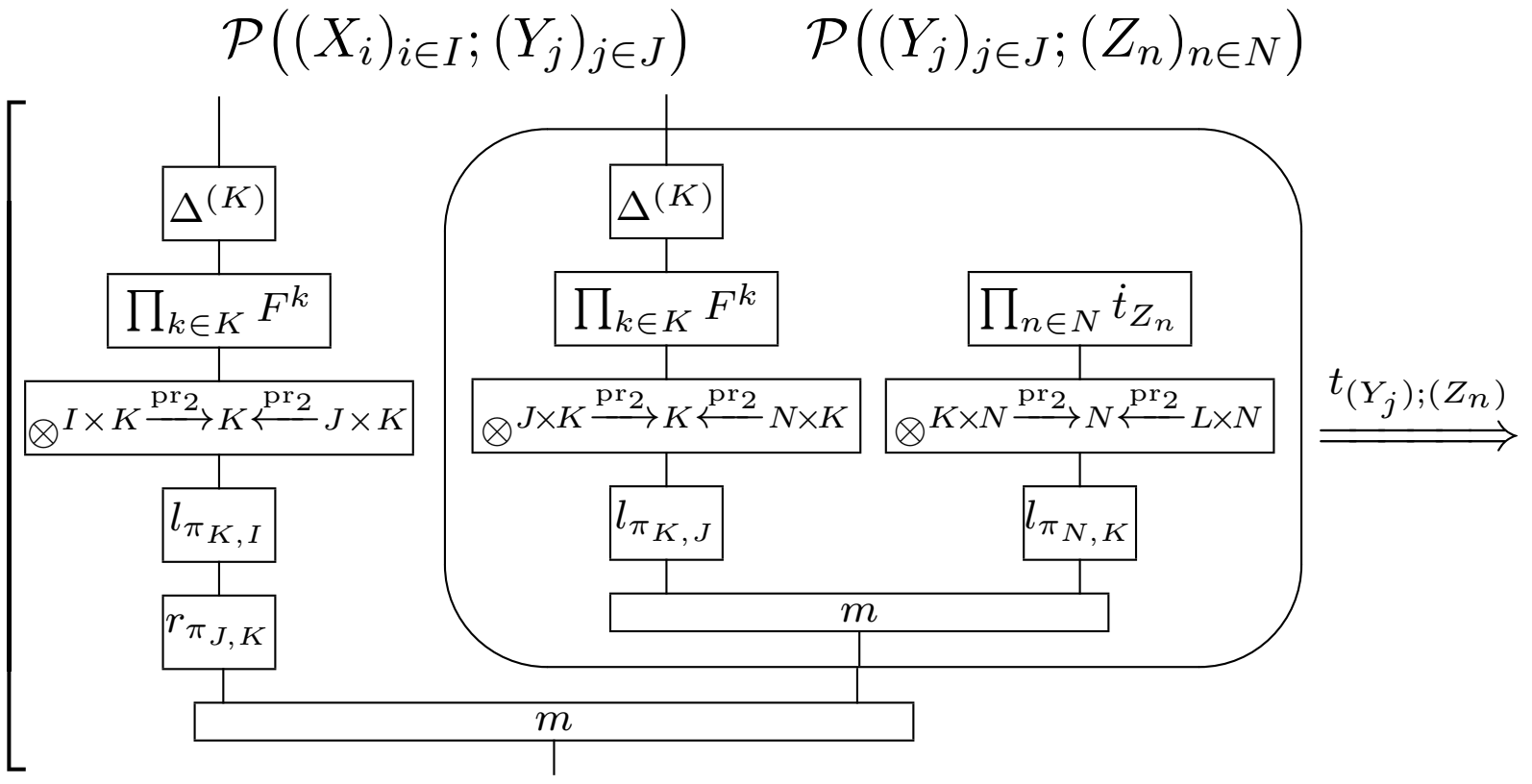

(2.2.2)

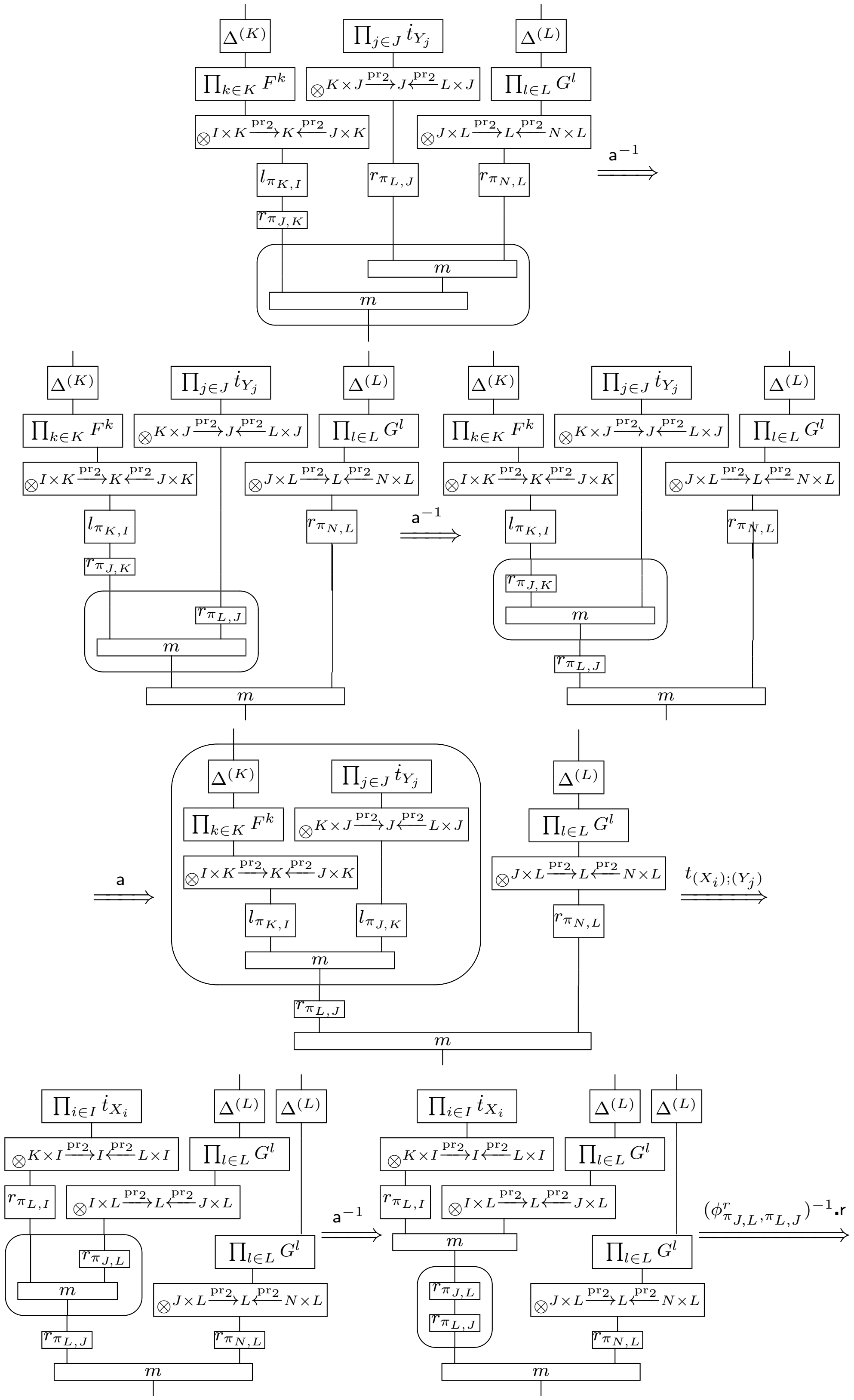
a−1
a−1
a
t(Xi);(Yj)
a−1
(φr πJ,L,πL,J)−1.r

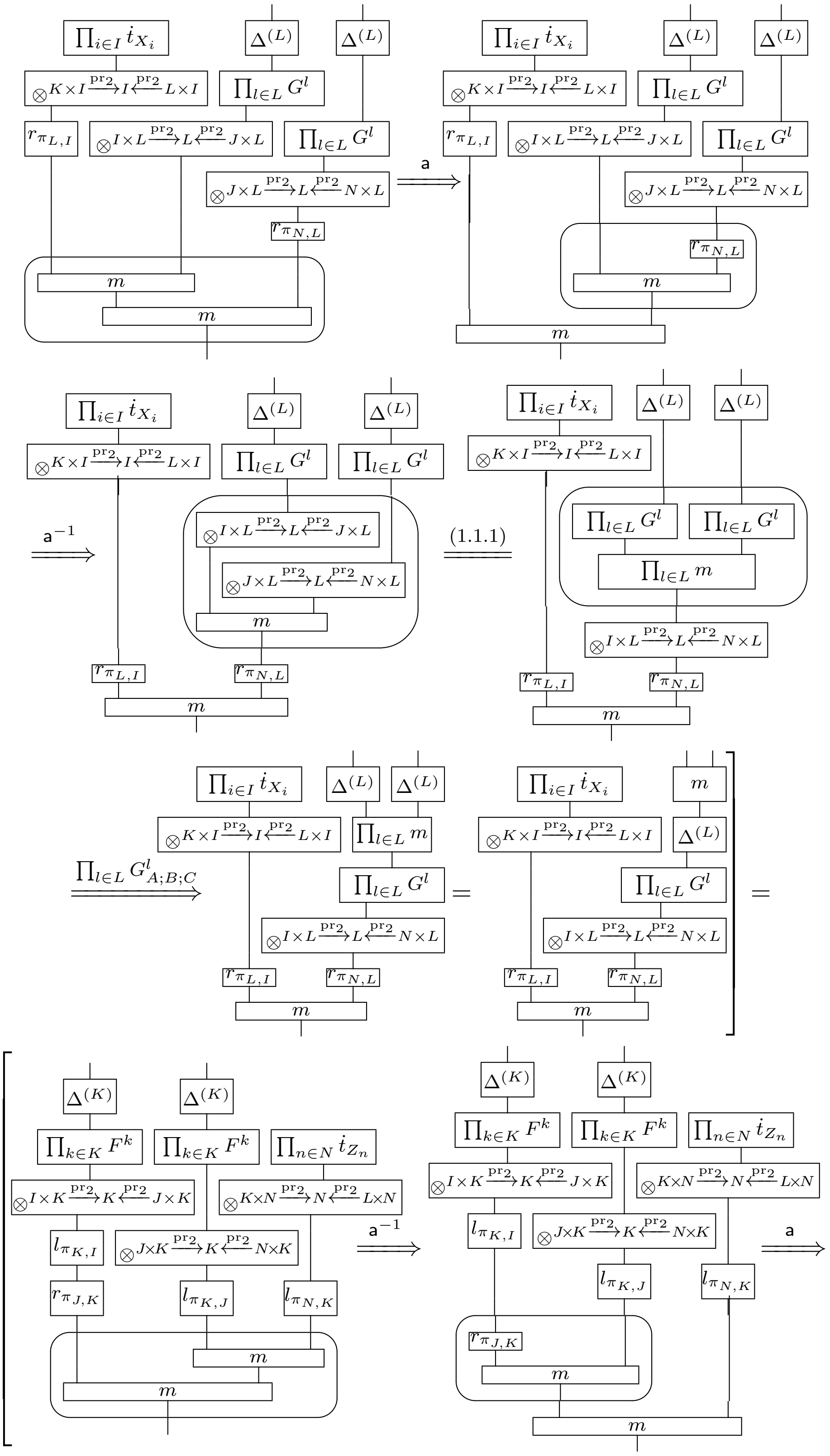
a
a^{-1}
(1.1.1)
$\prod_{l\in L} G^l_{A;B;C}$
=
=
a^{-1}
a

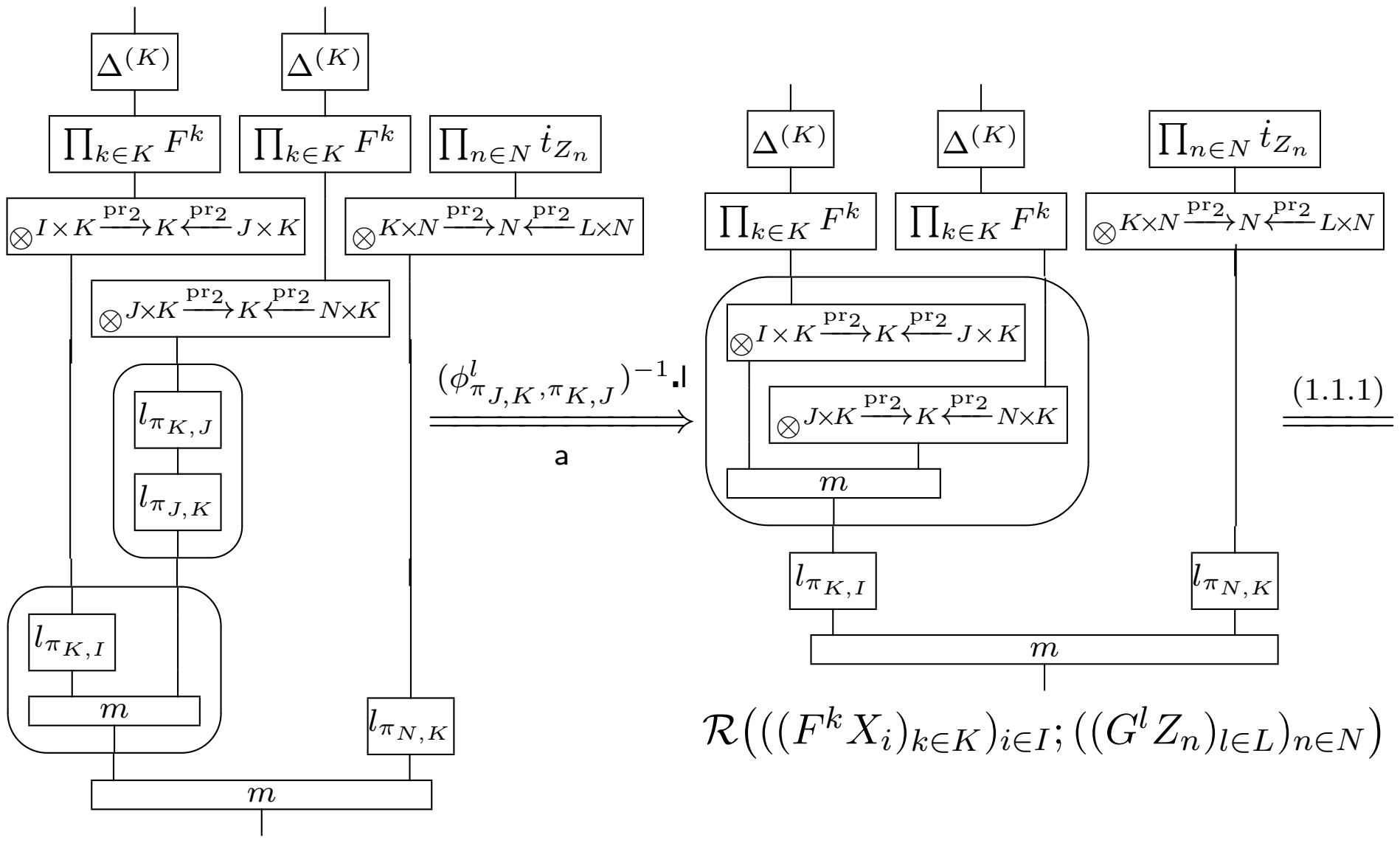

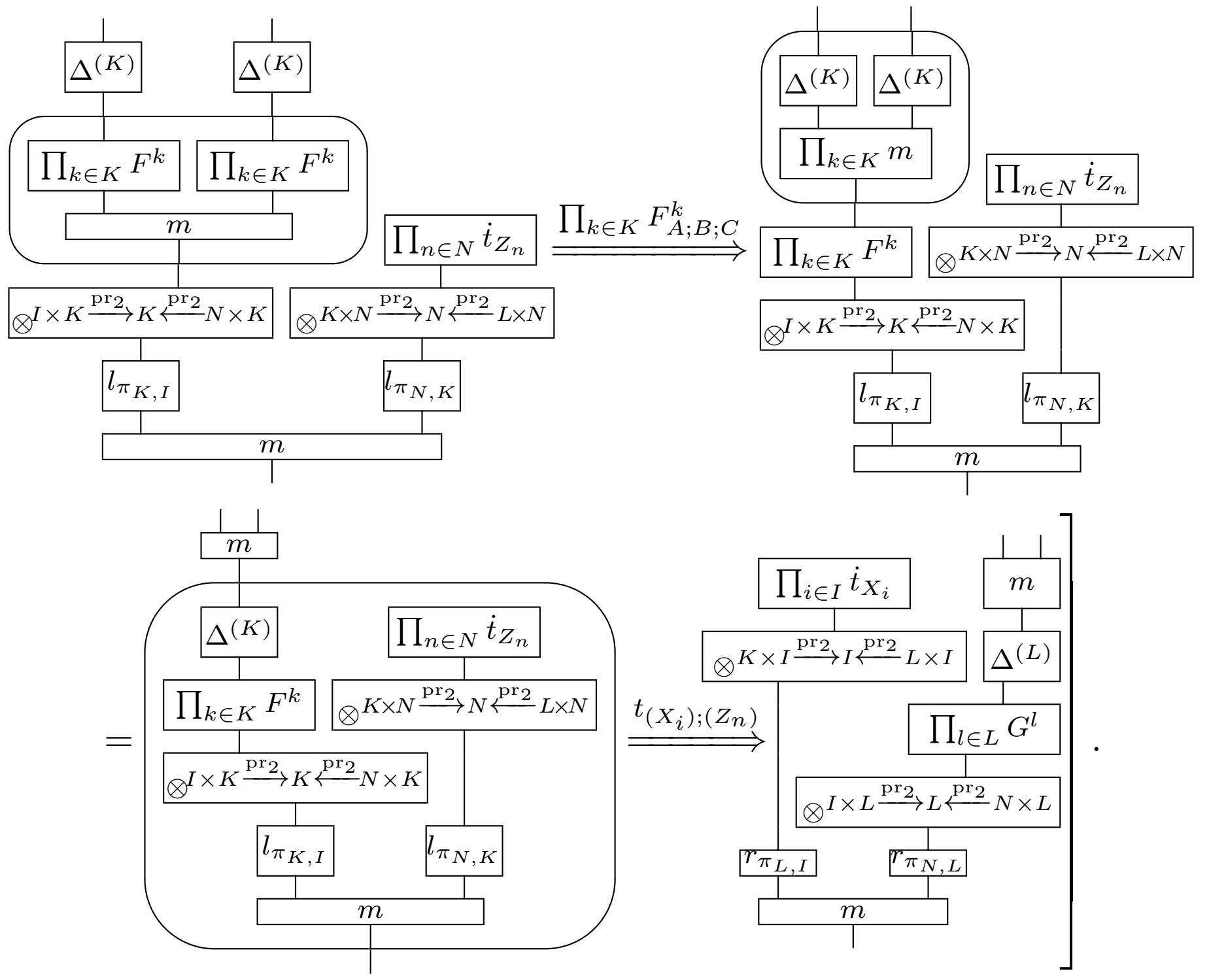

Note that

$$\big(K = K \times \mathbf{1} \underset{(12)}{\overset{\pi_{K,\mathbf{1}}}{-\!\!\|\!\!\to}} \mathbf{1} \times K = K\big) = 1_K = \big(K = \mathbf{1} \times K \underset{(12)}{\overset{\pi_{\mathbf{1},K}}{-\!\!\|\!\!\to}} K \times \mathbf{1} = K\big).$$

Using this we write another equation, which specifies the behaviour with respect to units:

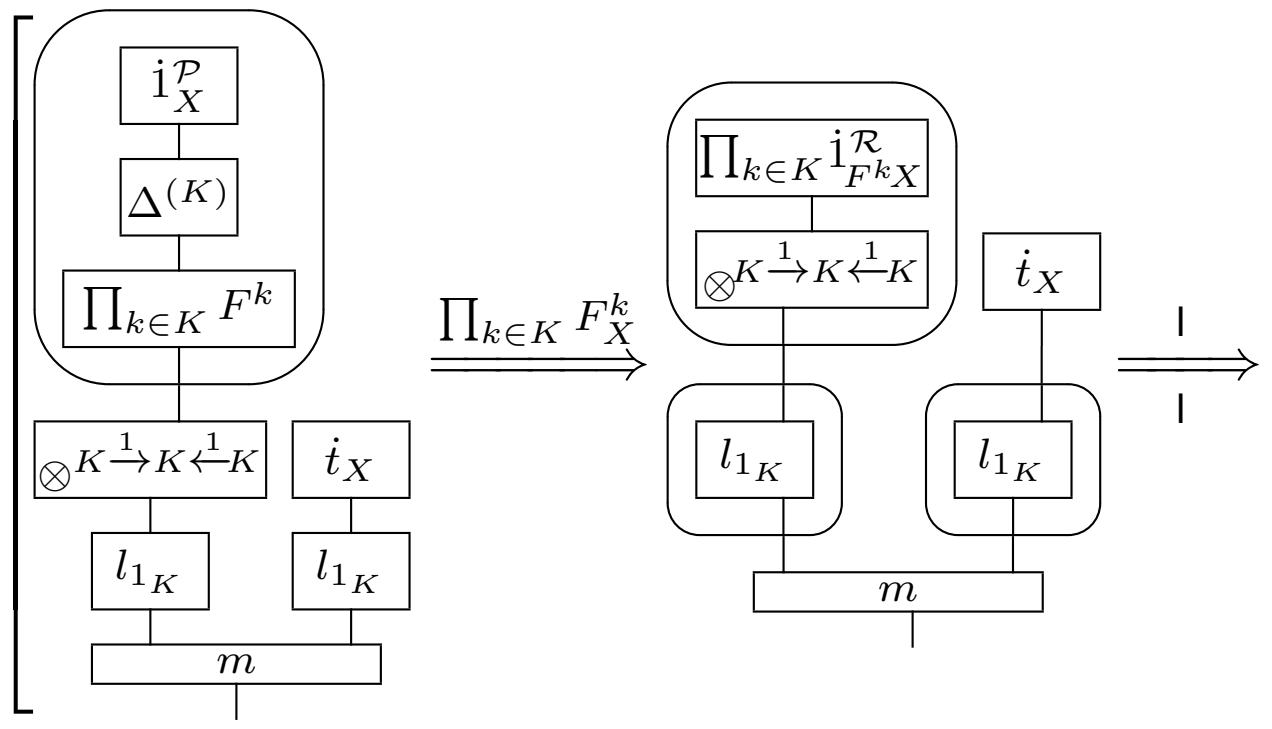

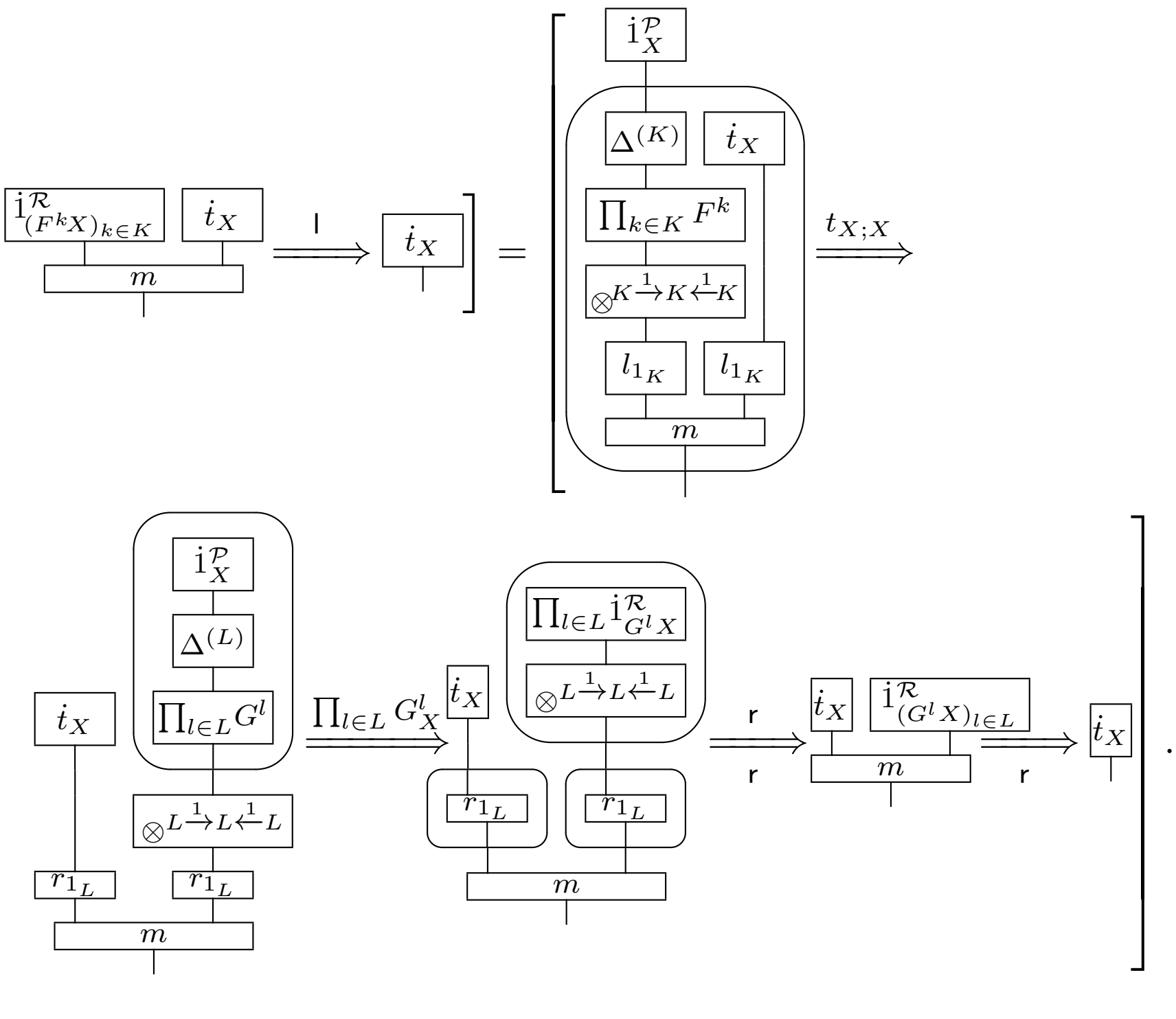

(2.2.3)

## 2.3 Modifications

**2.3.1 Definition.** Let $t, p\colon (F^k)_{k\in K} \to (G^l)_{l\in L}\colon \mathcal{P} \to \mathcal{R}$ be natural transformations of morphisms of biprops. A *modification* $f\colon t \to p\colon (F^k)_{k\in K} \to (G^l)_{l\in L}\colon \mathcal{P} \to \mathcal{R}$ is a family of

— morphisms $f_X \in \mathcal{R}((F^kX)_{k\in K}; (G^lX)_{l\in L})(t_X, p_X)$, $X \in \operatorname{Col}\mathcal{P}$

such that the following square of natural transformations commutes:

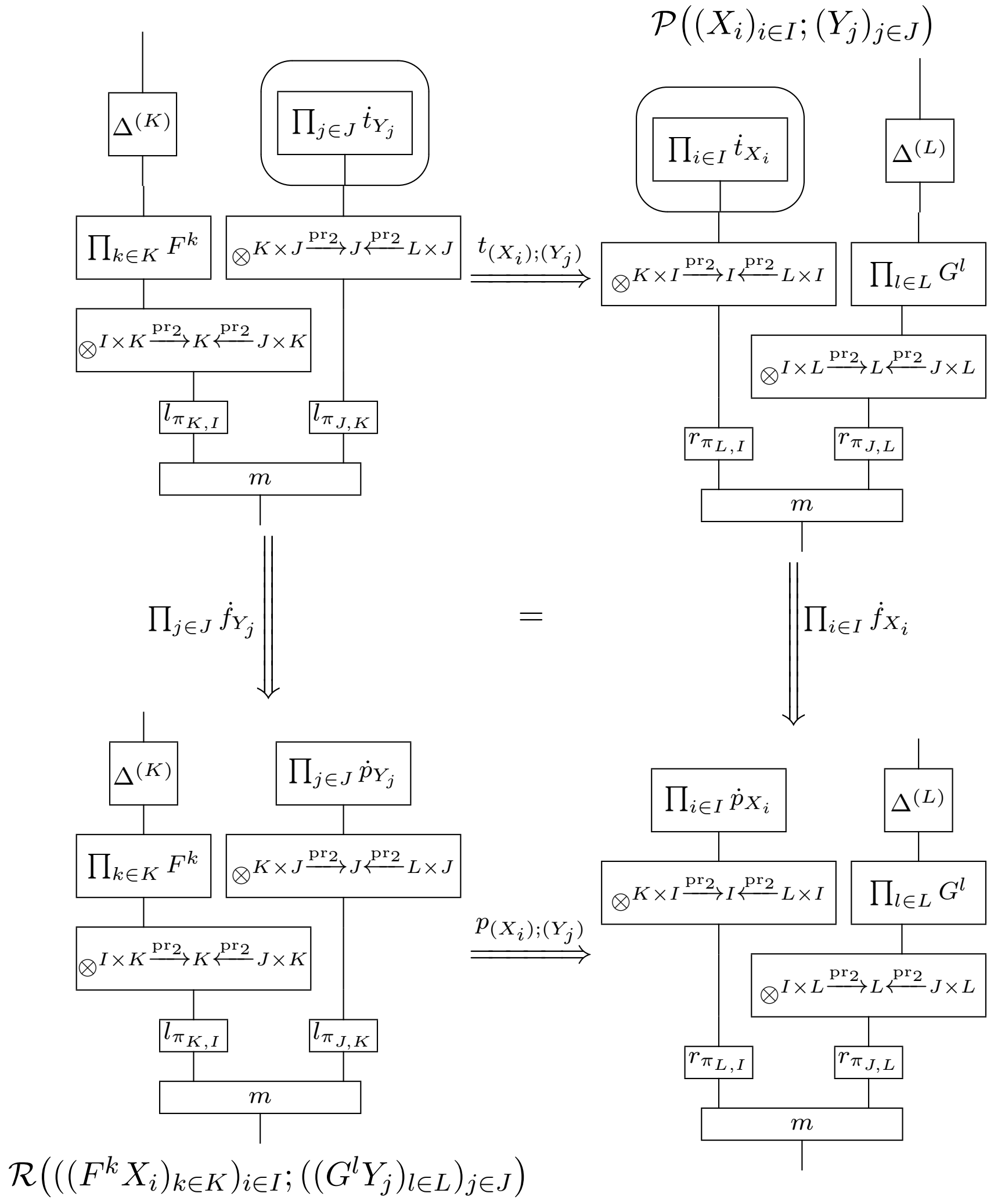

(2.3.1)

# 3 Declarations

## 3.1 Funding

Partial financial support was received from the National Centre of Competence in Research SwissMAP of the Swiss National Science Foundation (grant number 205607), Faculty of Mathematics, Computer Science and Natural Sciences of the University of Hamburg, the Isaac Newton Institute and London Mathematical Society, the School of Mathematics at the University of Edinburgh and from Università degli studi di Milano.

## 3.2 Conflicts of interest/Competing interests

The author has no relevant financial or non-financial interests to disclose.

Institute of Mathematics, NAS Ukraine, 3 Tereshchenkivska st., Kyiv, 01024, Ukraine